\documentclass{article}

\usepackage{latexsym}
\usepackage{amsfonts}
\usepackage{amssymb}
\usepackage{amsmath}
\usepackage{amscd}
\usepackage{eucal}
\usepackage{amsthm}

\usepackage{latexsym}
\usepackage{amsfonts}
\usepackage{amssymb}
\usepackage{amsmath}
\usepackage{amscd}
\usepackage{eucal}
\usepackage[all, knot]{xy}
\xyoption{arc}
\xyoption{rotate}
\xyoption{dvips}
\usepackage{amsthm}
\usepackage[dvips]{pict2e}

\usepackage{verbatim} 
\usepackage[dvips]{graphicx}

\newcommand{\G}{\mathbb G}

\newtheorem{Def}{Definition}[section]
\newtheorem{thm}[Def]{Theorem}
\newtheorem{lem}[Def]{Lemma}
\newtheorem{prop}[Def]{Proposition}
\newtheorem{cor}[Def]{Corollary}

\newtheoremstyle{example}{\topsep}{\topsep}%
     {}
     {}
     {\bfseries}
     {.}
     {8pt}
     {\thmname{#1}\thmnumber{ #2}\thmnote{ #3}}

   \theoremstyle{example}
   
\newtheorem{rem}[Def]{Remark}

\newcommand{\lra}{\ensuremath{\longrightarrow}}
\newcommand{\cat}[1]{\ensuremath{\mbox{\bfseries {\upshape {#1}}}}}
\newcommand{\cl}[1]{\ensuremath{\mathcal {#1}}}

\newcommand{\ed}{\end{document}}

\newcommand{\numarabic}{\renewcommand{\labelenumi}{\arabic{enumi}.}}

\newcommand{\sunit}{\setlength{\unitlength}{1mm}}
\newcommand{\setunit}[1]{\setlength{\unitlength}{#1}}

\newcommand{\triarrow}{
\sunit
\begin{picture}(24,0)(-0.5,-1)
\put(0,0){\line(1,0){20}}
\put(19.8,-0.2){\makebox(0,0)[l]{\ttr}}
\end{picture}}

\newcommand{\degenarrow}{
\begin{picture}(0,0)
\qbezier[20](0,0)(-3,-2.25)(-6,-4.5)
\put(-6,-4.5){\vector(-12,-9){0.1}}
\end{picture}}

\newcommand{\psinv}[6]{\xy
(#1,0)*+{#3}="x";
(#2,0)*+{#5}="y";
{\ar@<.7ex>^{#4} "x"; "y"};
{\ar@<.7ex>^{#6} "y"; "x"};
\endxy
}

\newcommand{\ttr}{\ensuremath{\triangleright}}

\newcommand{\bs}{\boldsymbol}

\begin{document}

\newtheoremstyle{example}{\topsep}{\topsep}%
     {}
     {}
     {\bfseries}
     {.}
     {2pt}
     {\thmname{#1}\thmnumber{ #2}\thmnote{ #3}}

   \theoremstyle{example}
   \newtheorem{nota}{Notation}
   \newtheorem{Defi}{Definition}

\title{The periodic table of $n$-categories for low dimensions II: degenerate tricategories} 
\author{Eugenia Cheng\\Department of Mathematics, Universit\'e de Nice Sophia-Antipolis \\E-mail: eugenia@math.unice.fr\\ \\Nick Gurski\\Department of Mathematics, Yale University\\E-mail: michaeln.gurski@yale.edu}
\date{June 2007}
\maketitle

%


\begin{abstract}

We continue the project begun in \cite{cg1} by examining degenerate tricategories and comparing them with the structures predicted by the Periodic table.  For triply degenerate tricategories we exhibit a triequivalence with the partially discrete tricategory of commutative monoids.  For the doubly degenerate case we explain how to construct a braided monoidal category from a given doubly degenerate category, but show that this does not induce a straightforward comparison between \cat{BrMonCat} and \cat{Tricat}.  We show how to alter the natural structure of \cat{Tricat} in two different ways to provide a comparison, but show that only the more brutal alteration yields an equivalence.  Finally we study degenerate tricategories in order to give the first fully algebraic definition of monoidal bicategories and the full tricategory structure \cat{MonBicat}.  
\end{abstract}

\maketitle
\tableofcontents

\section*{Introduction}
\addcontentsline{toc}{section}{Introduction}

This work is a continuation of the work begun in \cite{cg1}, studying the ``Periodic Table'' of $n$-categories proposed by Baez and Dolan \cite{bd3}.  The idea of the Periodic Table is to study ``degenerate'' $n$-categories, that is, $n$-categories in which the lowest dimensions are trivial.  For small $n$ this is supposed to yield well-known algebraic structures such as commutative monoids or braided monoidal categories; this helps us understand some specific part of the whole $n$-category via better-known algebraic structures, and also helps us to try to predict what $n$-categories should look like for higher $n$.

More precisely, the idea of degeneracy is as follows.  Consider an $n$-category in which the lowest non-trivial dimension is the $k$th dimension, that is, there is only one cell of each dimension lower than $k$.  We call this a ``$k$-degenerate $n$-category''.  We can then perform a ``dimension shift'' and consider the $k$-cells of the old $n$-category to be 0-cells of a new $(n-k)$-category, as shown in the schematic diagram in Figure~\ref{dshift}.

\begin{figure}\label{shift}
\caption{Dimension-shift for $k$-fold degenerate $n$-categories} \label{dshift}
\setunit{1mm}
\begin{tabular}{|lcc|} \hline &&\\[-4pt]
{\bfseries ``old" $n$-category} & \triarrow & {\bfseries ``new" $(n-k)$-category} \\[8pt] \hline && \\[-4pt]
\hspace{1em}$\left.\begin{array}{l}
\mbox{0-cells} \\ \mbox{1-cells} \\ \hspace{1em}\vdots \\[4pt] \makebox(13,0)[r]{$(k-1)$-cells}
\end{array}\right\}$ trivial &&\\[32pt]
\hspace{1.4em}$k$-cells & \triarrow & 0-cells \\[6pt]
$(k+1)$-cells & \triarrow & 1-cells \\
\hspace{2.2em} \vdots & \vdots & \vdots \\[2pt]
\hspace{1.4em}$n$-cells & \triarrow & $(n-k)$-cells \\[4pt] \hline
\end{tabular}
\end{figure}

This yields a ``new'' $(n-k)$-category, but it will always have some special extra structure: the $k$-cells of the old $n$-category have $k$ different compositions defined on them (along bounding cells of each lower dimension), so the 0-cells of the ``new'' $(n-k)$-category must have $k$ multiplications defined on them, interacting via the interchange laws from the old $n$-category.  Likewise every cell of higher dimension will have $k$ ``extra'' multiplications defined on them as well as composition along bounding cells.  

In \cite{bd3}, Baez and Dolan {\it define} a ``$k$-tuply monoidal $(n-k)$-category'' to be a $k$-degenerate $n$-category, but {\it a priori} it should be an $(n-k)$-category with $k$-monoidal structures on it, interacting via coherent pseudo-invertible cells.  A direct definition has not yet been made for general $n$ and $k$.  (Balteanu et al \cite{bfsv1} study a lax version of this, where the monoidal structures interact via non-invertible cells; this gives different structures not studied in the present work.)

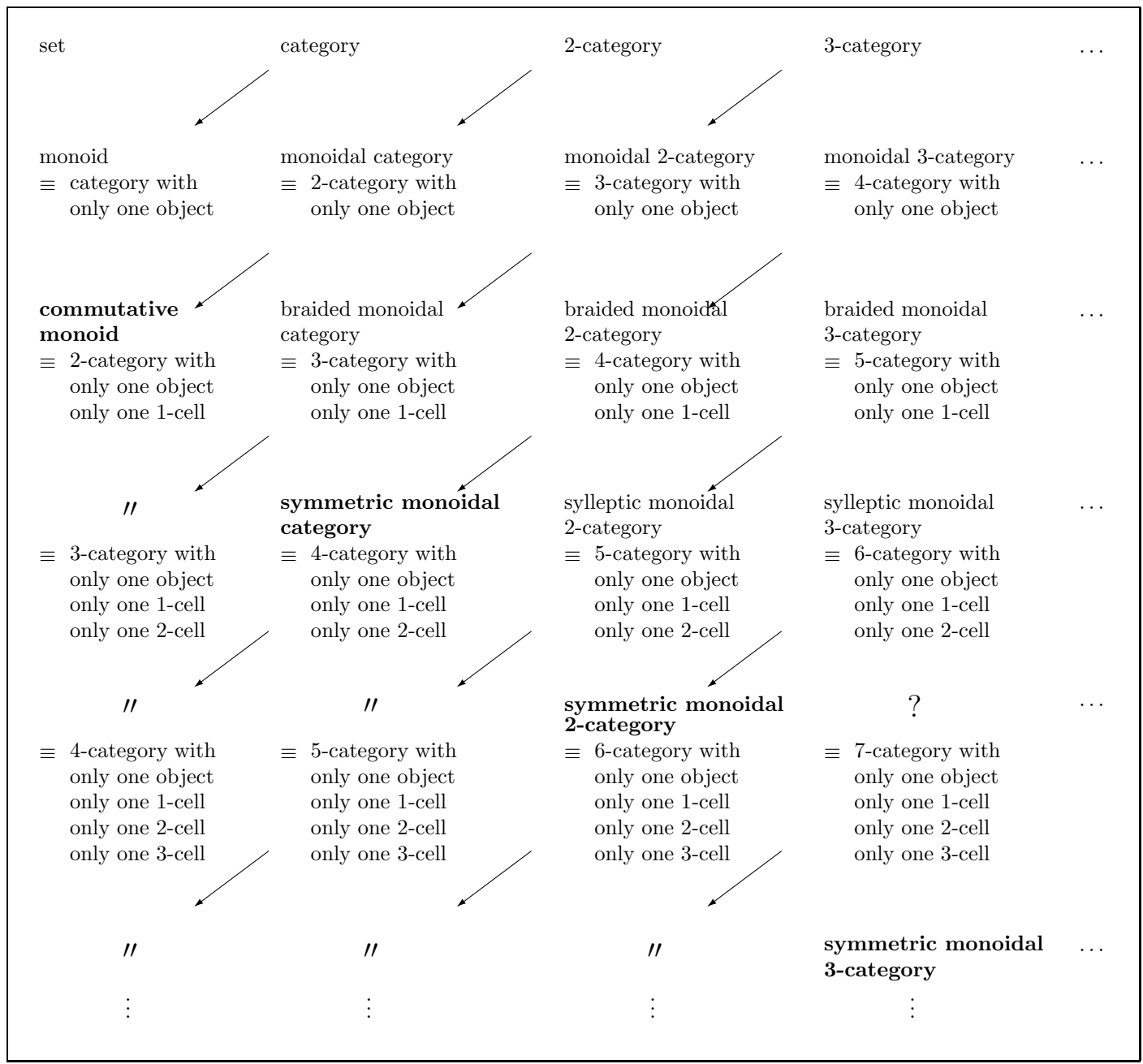
\begin{figure}
\caption{The hypothesised Periodic Table of $n$-categories} \label{ptable}

\setunit{2mm}  \nopagebreak[4]
\begin{picture}(80,90)(16,-8)

\put(0,37.5){
\renewcommand{\tabcolsep}{1.5em}

\begin{tabular}{|lllll|} \hline &&&&\\
set & category & 2-category & 3-category & \makebox(0,0)[r]{$\cdots$}  \\[40pt]
monoid & monoidal category & monoidal 2-category & monoidal 3-category & \makebox(0,0)[r]{$\cdots$} \\
\makebox(2.5,0.5)[l]{$\equiv$}category with & \makebox(2.5,0.5)[l]{$\equiv$}2-category with & \makebox(2.5,0.5)[l]{$\equiv$}3-category with & \makebox(2.5,0.5)[l]{$\equiv$}4-category with &\\ 
\makebox(2.5,0.5)[l]{}only one object & \makebox(2.5,0.5)[l]{}only one object & \makebox(2.5,0.5)[l]{}only one object & \makebox(2.5,0.5)[l]{}only one object &\\[35pt]
{\bfseries commutative}& braided monoidal & braided monoidal & braided monoidal & \makebox(0,0)[r]{$\cdots$}  \\
{\bfseries monoid}& category & 2-category & 3-category &\\
\makebox(2.5,0.5)[l]{$\equiv$}2-category with & \makebox(2.5,0.5)[l]{$\equiv$}3-category with & \makebox(2.5,0.5)[l]{$\equiv$}4-category with & \makebox(2.5,0.5)[l]{$\equiv$}5-category with &\\ 
\makebox(2.5,0.5)[l]{}only one object & \makebox(2.5,0.5)[l]{}only one object & \makebox(2.5,0.5)[l]{}only one object & \makebox(2.5,0.5)[l]{}only one object &\\
\makebox(2.5,0.5)[l]{}only one 1-cell & \makebox(2.5,0.5)[l]{}only one 1-cell & \makebox(2.5,0.5)[l]{}only one 1-cell & \makebox(2.5,0.5)[l]{}only one 1-cell &\\[30pt]
\makebox(8,0)[r]{\large $\prime\prime$}& {{\bfseries symmetric monoidal}} & sylleptic monoidal & {sylleptic monoidal} & \makebox(0,0)[r]{$\cdots$}  \\
& {\bfseries category} & 2-category & 3-category &\\
\makebox(2.5,0.5)[l]{$\equiv$}3-category with & \makebox(2.5,0.5)[l]{$\equiv$}4-category with & \makebox(2.5,0.5)[l]{$\equiv$}5-category with & \makebox(2.5,0.5)[l]{$\equiv$}6-category with &\\ 
\makebox(2.5,0.5)[l]{}only one object & \makebox(2.5,0.5)[l]{}only one object & \makebox(2.5,0.5)[l]{}only one object & \makebox(2.5,0.5)[l]{}only one object &\\
\makebox(2.5,0.5)[l]{}only one 1-cell & \makebox(2.5,0.5)[l]{}only one 1-cell & \makebox(2.5,0.5)[l]{}only one 1-cell & \makebox(2.5,0.5)[l]{}only one 1-cell &\\
\makebox(2.5,0.5)[l]{}only one 2-cell & \makebox(2.5,0.5)[l]{}only one 2-cell & \makebox(2.5,0.5)[l]{}only one 2-cell & \makebox(2.5,0.5)[l]{}only one 2-cell &\\ [20pt]
\makebox(8,0)[r]{\large $\prime\prime$}& \makebox(8,0)[r]{\large $\prime\prime$} & \makebox(16,0)[l]{{\bfseries symmetric monoidal}} & \makebox(8,0)[r]{\Large ?}& \makebox(0,0)[r]{$\cdots$}  \\
& & {\bfseries 2-category} & &\\
\makebox(2.5,0.5)[l]{$\equiv$}4-category with & \makebox(2.5,0.5)[l]{$\equiv$}5-category with & \makebox(2.5,0.5)[l]{$\equiv$}6-category with & \makebox(2.5,0.5)[l]{$\equiv$}7-category with &\\ 
\makebox(2.5,0.5)[l]{}only one object & \makebox(2.5,0.5)[l]{}only one object & \makebox(2.5,0.5)[l]{}only one object & \makebox(2.5,0.5)[l]{}only one object &\\
\makebox(2.5,0.5)[l]{}only one 1-cell & \makebox(2.5,0.5)[l]{}only one 1-cell & \makebox(2.5,0.5)[l]{}only one 1-cell & \makebox(2.5,0.5)[l]{}only one 1-cell &\\
\makebox(2.5,0.5)[l]{}only one 2-cell & \makebox(2.5,0.5)[l]{}only one 2-cell & \makebox(2.5,0.5)[l]{}only one 2-cell & \makebox(2.5,0.5)[l]{}only one 2-cell &\\ 
\makebox(2.5,0.5)[l]{}only one 3-cell & \makebox(2.5,0.5)[l]{}only one 3-cell & \makebox(2.5,0.5)[l]{}only one 3-cell & \makebox(2.5,0.5)[l]{}only one 3-cell &\\[30pt] 
\makebox(8,0)[r]{\large $\prime\prime$}& \makebox(8,0)[r]{\large $\prime\prime$}& \makebox(8,0)[r]{\large $\prime\prime$}& {\bfseries symmetric monoidal}& \makebox(0,0)[r]{$\cdots$} \\
&&& {\bfseries 3-category} &\\
\makebox(7.5,0)[r]{$\vdots$} & \makebox(7.5,0)[r]{$\vdots$}
& \makebox(7.5,0)[r]{$\vdots$} & \makebox(7.5,0)[r]{$\vdots$} & \\
&&&& \\
&&&& \\
 
\hline
\end{tabular}}

\put(22,76){\degenarrow}
\put(22,61){\degenarrow}
\put(22,46){\degenarrow}
\put(22,30){\degenarrow}
\put(22,12){\degenarrow}

\put(43.5,76){\degenarrow}
\put(43.5,61){\degenarrow}
\put(43.5,46){\degenarrow}
\put(43.5,30){\degenarrow}
\put(43.5,12){\degenarrow}

\put(64,76){\degenarrow}
\put(64,61){\degenarrow}
\put(64,46){\degenarrow}
\put(64,30){\degenarrow}
\put(64,12){\degenarrow}


\end{picture}
\end{figure}

The Periodic Table seeks to answer the question: exactly what sort of $(n-k)$-category structure does the degeneracy process produce? Figure~\ref{ptable} shows the first few columns of the hypothesised Periodic Table: the $(n,k)$th entry predicts what a $k$-degenerate $n$-category ``is''. (In this table we follow Baez and Dolan and omit the word ``weak" understanding that all the $n$-categories in consideration are weak.)

One consequence of the present work is that although $k$-tuply monoidal $(n-k)$-categories and $k$-degenerate $n$-categories are related, we see that the relationship is not straightforward.  So in fact we need to consider three possible structures for each $n$ and $k$:

\begin{itemize}
\item $k$-degenerate $n$-categories
\item $k$-tuply monoidal $(n-k)$-categories
\item the $(n,k)$th entry of the Periodic Table.
\end{itemize}

In \cite{cg1} we examined the top left hand corner of the table, that is, degenerate categories and degenerate bicategories.  We found that we had to be careful about the exact meaning of ``is''.  The main problem is the presence of some unwanted extra structure in the ``new" $(n-k)$-categories in the form of distinguished elements, arising from the structure constraints in the original $n$-categories --- a specified $k$-cell structure constraint in the ``old" $n$-category will appear as a distinguished 0-cell in the ``new" $(n-k)$-category under the dimension-shift depicted in Figure~\ref{dshift}. (For $n=2$ this phenomenon is mentioned by Leinster in \cite{lei4} and was further described in a talk \cite{lei10}.)

This problem becomes worse when considering functors, transformations, modifications, and so on, as we will discuss in the next section.

\subsection{Totalities of structures}

Broadly speaking we have two different aims:

\begin{enumerate}
\item Theoretical: to make precise statements about the claims of the Periodic Table by examining the {\it totalities} of the structures involved, that is, not just the degenerate $n$-categories but also all the higher morphisms between them.

\item Practical: to find the structures predicted by the Periodic Table in a way that somehow naturally arises from degenerate tricategories and their morphisms.
\end{enumerate}

The point of (2) is that in practice we may simply want to know that a given doubly degenerate tricategory is a braided monoidal category, or that a given functor is a braided monoidal functor, for example, without needing to know if the {\it theory} of doubly degenerate tricategories corresponds to the theory of braided monoidal categories.  The motivating example discussed in \cite{bd3} is the degenerate $n$-category of ``manifolds with corners embedded in $n$-cubes''; work towards constructing such a structure appears in \cite{bl1} and \cite{cg2}.

In this work we see that although the tricategories and functors behave as expected, the higher morphisms are much more general than the ones we want.  Moreover, for (1) we see that the overall dimensions of the totalities do not match up.  On the one hand we have $k$-degenerate $n$-categories, which naturally organise themselves into an $(n+1)$-category---the full sub-$(n+1)$-category of \cat{nCat}; by contrast, the structure predicted by the Periodic Table is an $(n-k)$-category with extra structure, and these organise themselves into an $(n-k+1)$-category---the full sub-$(n-k+1)$-category of \cat{(n-k)Cat}.  In order to compare an $(n+1)$-category with an $(n-k+1)$-category we either need to remove some dimensions from the former or add some to the latter.  

The most obvious thing to do is add dimensions to the latter in the form of higher identity cells.  However, we quickly see that this does not yield an equivalence of $(n+1)$-categories because the $(n+1)$-cells of \cat{nCat} are far from trivial. 

So instead we try to reduce the dimensions of \cat{nCat}.  We cannot in general apply a simple truncation to $j$-dimensions as this will not result in a $j$-category.  Besides, we would also like to restrict our higher morphisms in order to achieve a better comparison with the structures given in the Periodic Table---{\it a priori} our higher morphisms are too general.

So we perform a construction analogous to the construction of ``icons'' \cite{lp1}.  The idea of icons is to organise bicategories into a {\it bicategory} rather than a tricategory, by discarding the modifications, selecting only those transformations that have all their components the identity, and altering their composition to ensure closure.  

For the $(n+1)$-category of $n$-categories we can try to perform an analogous ``collapse'' to obtain a $j$-category of $n$-categories, for any $2 \leq j \leq n$.  We discard all cells of \cat{nCat} of dimension greater than $j$, and for dimensions up to $j$ we select only those special cases where the components in the lowest $(n-j)$ dimensions are the identity; we then redefine composition to force closure.  More explicitly, recall that for an $m$-cell in \cat{nCat} the data is essentially:

\begin{itemize}
\item for all 0-cells an $(m-1)$-cell
\item for all 1-cells an $m$-cell
\item for all 2-cells an $(m+1)$-cell
\item for all 3-cells an $(m+2)$-cell
\item $\vdots$
\item for all $(n-m-1)$-cells an $(n-2)$-cell
\item for all $(n-m)$-cells an $(n-1)$-cell
\item for all $(n-m+1)$-cells an $n$-cell
\end{itemize}
We can force the first $(n-j+1)$ of these to be the identity, so that the first non-trivial piece of data is ``for all $(n-j+1)$-cells a $\ldots$''.  This automatically forces all the morphisms of dimension higher than $j$ to be the identity, and once we have redefined composition to force closure, we have a $j$-category of $n$-categories.  Then, using $j=n-k+1$ and restricting to the $k$-degenerate case, we find we have a naturally arising $(n-k+1)$-category of $k$-degenerate $n$-categories.


So to compare with braided monoidal categories we can make a {\it bicategory} of tricategories \cite{gg1} with
\begin{itemize}
\item 0-cells: tricategories
\item 1-cells: functors between them
\item 2-cells: lax transformations whose 1- and 2-cell components are identities, but 3-cell components are non-trivial. 
\end{itemize}
Taking the full sub-bicategory whose 0-cells are the doubly degenerate tricategories, this does give rise to braided monoidal categories, braided monoidal functors, and monoidal transformations; however the correspondence is not at all straightforward and moreover we do not get a biequivalence of bicategories.

For monoidal bicategories we make a {\it tricategory} of tricategories with
\begin{itemize}
\item 0-cells: tricategories
\item 1-cells: functors between them
\item 2-cells: lax transformations whose 1-cell components are identities, but 2- and 3-cell components are non-trivial
\item 3-cells: lax modifications whose 2-cell components are identities, but 3-cell components are non-trivial.
\end{itemize}
Taking the full sub-tricategory whose 0-cells are the degenerate tricategories, this gives rise to monoidal bicategories, monoidal functors, monoidal transformations and monoidal modifications; in fact this is how we define them, as we will later discuss.

Another approach would be to restrict the degenerate tricategories much more.  We could adopt the philosophy that for $k$-degenerate $n$-categories, when a piece of data says ``for all $m$-cells a $\ldots$'', then if $m$ is one of the trivial dimensions $m<k$ this data should be trivial as well.

This approach seems less natural from the point of view of $n$-categories, and might exclude some examples; it produces an equivalence of structures for the case of braided monoidal categories that seems almost tautological.  Essentially we have made a functor
	\[\cat{BrMonCat}_4 \lra \cat{Tricat}\]
where $\cat{BrMonCat}_4$ is the partially discrete tetracategory formed by adding higher identities to \cat{BrMonCat}; we can fairly easily construct such a functor by choosing all extra structure to be identities.  We might try to look for the essential image of this functor, to find a sub-tetracategory of \cat{Tricat} that is equivalent to $\cat{BrMonCat}_4$, but it is currently much too hard for us to work with tetracategories in this way.  Instead, we can do something very crude---take the precise image of this functor (which will be trivial in dimensions 3 and 4) and somehow force it to be a bicategory.  It is not clear what this achieves.

\subsection{Results}

The main results of \cite{cg1} can be summed up as follows.  (Here we write ``degenerate'' for ``1-degenerate'', and ``doubly degenerate'' for ``2-degenerate'', although in general we also use ``degenerate'' for any level of degeneracy.)

\begin{itemize}

\item Comparing each degenerate category with the monoid formed by its 1-cells, we exhibit an equivalence of categories of these structures, but not a biequivalence of bicategories.

\item Comparing each doubly degenerate bicategory with the commutative monoid formed by its 2-cells, we exhibit a biequivalence of bicategories of these structures, but not an equivalence of categories or a triequivalence of tricategories.

\item Comparing each degenerate bicategory with the monoidal category formed by its 1-, 2-, and 3-cells, we exhibit an equivalence of categories of these structures, but not a biequivalence of bicategories or a triequivalence of tricategories.

\end{itemize}

In the present work we proceed to the next dimension and study degenerate tricategories.  We use the fully algebraic definition of tricategory given in \cite{gur1}; this is based on the definition given in \cite{gps1} which is not fully algebraic.  The results can be summed up as follows, but cannot be stated quite so succinctly.

\begin{itemize}

\item Comparing each triply degenerate tricategory with the commutative monoid formed by its 3-cells, we exhibit a triequivalence of tricategories of these structures, but not an equivalence of categories, a biequivalence of bicategories, or a tetra-equivalence of tetra-categories.

\item We show how doubly degenerate tricategories give rise to braided monoidal categories, and similarly functors.  For comparisons we use the bicategory of tricategories described above.  We exhibit comparison functors in both directions, but no equivalence except in the ``tautological'' case (see above).

\item A degenerate tricategory gives, by definition, a monoidal bicategory formed by its 1-cells and 2-cells.  The totality of monoidal bicategories has not previously been understood; here we use the tricategory of bicategories described above, and use this to define a tricategory \cat{MonBicat} of monoidal bicategories, in which the higher-dimensional structure is not directly inherited from \cat{Tricat}. 

\end{itemize}

These results might be thought of as being significantly ``worse'' than the results for degenerate categories and bicategories, in the sense that it is much harder to make any precise statements about the structures in question being equivalent.  Indeed it is a general principle of $n$-category theory that there is a large (perhaps disproportionate) difference between $n=2$ and $n=3$: bicategories exhibit a level of coherence that does not generalise to higher dimensions.  Two of the main examples of this are:

\begin{enumerate}

\item Bicategories and weak functors form a category, although ``morally'' they should only be expected to form a tricategory.  However for $n \geq 3$ weak $n$-categories and weak functors cannot form a category as composition is not strictly associative or unital, at least in the algebraic case.  Note that bicategories, weak functors and weak transformations do not form a bicategory; in general $n$-categories should not form any coherent structure in fewer than $(n+1)$-dimensions without restricting the higher morphisms between them and altering their composition as described above.

\item Every bicategory is biequivalent to a (strict) 2-category.  However, it is not the case that every tricategory is triequivalent to a strict 3-category, and thus the general result for $n \geq 3$ is expected to be more subtle.  The obstruction for tricategories is exactly what produces a braiding rather than a symmetry for doubly degenerate tricategories (considered as monoidal categories), and is also what enables weak 3-categories to model homotopy 3-types where strict 3-categories cannot.  Another way of expressing this difference in coherence is that in a bicategory every diagram of constraints commutes, whereas in a tricategory there are diagrams of constraints that might not commute.

\end{enumerate}

The second remark above shows the importance of understanding the relationship between doubly degenerate tricategories and braided monoidal categories.  We will show that although every doubly degenerate tricategory does give rise to a braided monoidal category of its 2-cells and 3-cells, the relationship is not straightforward.  The process of producing the braiding is complicated, and there is a great deal of ``extra structure'' on the resulting braided monoidal category.  The disparity is even greater for functors, transformations and modifications, meaning that it is not clear in what sense braided monoidal categories and doubly degenerate tricategories should be considered equivalent.  In fact the source of all this difficulty can be seen in the case of doubly degenerate bicategories, although in that case many of the difficulties can resolved to produce a strict Eckmann-Hilton argument.  When categorifying this to the case of doubly degenerate tricategories, we might hope to produce a categorified Eckmann-Hilton argument yielding a categorified commutativity, i.e. a braiding.  However, the difficulties involved are also ``categorified'', as we will further discuss in Section~\ref{overview}.

The organisation of the paper is as follows; it is worth noting that each section is significant for different reasons, as we will point out.  In Section~\ref{3degen} we examine triply degenerate tricategories; the significance of this section is that this is a ``stable'' case, and the results therefore have implications for the Stabilisation Hypothesis.  In Section~\ref{2degen} we examine doubly degenerate tricategories, whose significance we have discussed above; we draw the reader's attention to the extensive informal overview given in Section~\ref{overview} as the severe technical details are in danger of obscuring the important principles involved.

In Section~\ref{1degen} we examine degenerate tricategories (i.e. 1-degenerate tricategories); the main purpose of this section is to give the first full definition of algebraic monoidal bicategories, together with their functors, transformations and modifications.  Some of the large technical diagrams are deferred to the Appendix.

Finally we note that the present paper is necessarily most concerned with studying precisely what structures do arise as degenerate tricategories, since good comparisons do not naturally arise. For the purposes of correctly interpreting the Periodic Table, it seems likely that a direct definition of ``$k$-tuply monoidal'' higher category will be more fruitful.  The case of doubly degenerate tricategories shows us that a $k$-degenerate $n$-category does not give rise to $k$ monoidal structures on the associated $(n-k)$-category in a straightforward way; however it is possible that $k$-tuply monoidal $(n-k)$-categories defined directly could more naturally yield the desired entries in the Periodic Table.

\section{Triply degenerate tricategories}\label{3degen}

In this section, we will study triply degenerate tricategories and the higher morphisms between them---functors, transformations, modifications and perturbations.  By the Periodic Table, triply degenerate tricategories are expected to be commutative monoids; by results of \cite{cg1} we now expect them to be commutative monoids {\it equipped with some distinguished invertible elements} arising from the structure constraints in the tricategory.  The process of finding how many such elements there are is highly technical and not particularly enlightening; we simply examine the data and axioms for a tricategory and calculate which constraints determine the others in the triply degenerate case.  The importance of these results is not in the exact number of distinguished invertible elements, but rather in the fact that there are any at all, and more than in the bicategory case.  We expect $n$-degenerate $n$-categories to have increasing numbers of distinguished invertible elements as $n$ increases, and thus for the precise algebraic situation to become more and more intractible in a somewhat uninteresting way.

The other important part of this result is to examine whether the higher morphisms between triply degenerate tricategories rectify the situation---if any higher morphisms essentially ignore the distinguished invertible elements already specified, then we can still have a structure equivalent to commutative monoids.  For doubly degenerate bicategories, this happened at the transformation level; for triply degenerate tricategories, this happens at the modification level.  As expected from results of \cite{cg1}, the top level morphisms, that is the perturbations, destroy the possibility of an equivalence on the level of tetracategories.

Throughout this section we use results of \cite{cg1} to characterise the (single) doubly degenerate hom-bicategory of a triply degenerate tricategory.

\subsection{Basic results}

The overall results for triply degenerate tricategories are as follows; we will discuss the calculations that lead to these results in the following sections.  We should also point out that the results in this section show that the higher-dimensional hypotheses we made in \cite{cg1} are incorrect. 

\numarabic

\begin{thm} \ \smallskip
\begin{enumerate}
\item  A triply degenerate tricategory $T$ is precisely a commutative monoid $X_{T}$ together with eight distinguished invertible elements $d, m, a, l, r, u, \pi, \mu$. 
\item  Extending the above correspondence, a weak functor $S \rightarrow T$ is precisely a monoid homomorphism $F:S \rightarrow T$ together with four distinguished invertible elements $m_{F}, \chi, \iota, \gamma$.
\item  Extending the above correspondence, a tritransformation $\alpha:F \rightarrow G$ is precisely the assertion that $(F,m_{F}) = (G,m_{G})$ together with distinguished invertible elements $\Pi$ and $\alpha_{T}$.
\item Extending the above correspondence, a trimodification $m:\alpha \Rightarrow \beta$ is precisely the assertion that $\alpha$ and $\beta$ are parallel.
\item Extending the above correspondence, a perturbation $\sigma:m \Rrightarrow n$ is precisely an element $\sigma$ in $T$.
\end{enumerate}
\end{thm}

\subsection{Tricategories}

In this section we perform the calculations for the triply degenerate tricategories themselves.  First we prove a useful lemma concerning adjoint equivalences.  The data for a tricategory involves the specification of various adjoint equivalences whose components are themselves adjoint equivalences in the doubly-degenerate hom-bicategories.  We are thus interested in adjoint equivalences in doubly degenerate bicategories.

\begin{lem}
Let $B$ be a doubly degenerate bicategory.  Then an adjoint equivalence $(f, g, \eta, \varepsilon)$ in $B$ consists of an invertible element $\eta \in X_{B}$ with $\varepsilon = \eta^{-1}$.
\end{lem}
\begin{proof}
The triangle identities yield the following equation in any bicategory.
\[
\eta * 1_{g} = a^{-1} \circ 1_{g} * \varepsilon^{-1} \circ r^{-1}_{g} \circ l_{g}
\]
Using the fact that $B$ is doubly degenerate, we see that in the commutative monoid $X_{B}$ (with unit written as 1) $a=1, 1_{g} = 1$, and $r = l$.  We also note that $* = \circ$, so the above equation reduces to the fact that $\eta$ and $\varepsilon$ are inverse to each other.
\end{proof}

{\it A priori}, a triply degenerate tricategory $T$ consists of the following data, which we will need to try to ``reduce'':
\begin{itemize}
\item a single object $\star$;
\item a doubly degenerate bicategory $T(\star, \star)$, which will be considered as a commutative monoid with distinguished invertible element, $(T, d_{T})$;
\item a weak functor $T(\star, \star) \times T(\star, \star) \rightarrow T(\star, \star)$, which will be considered as a monoid homomorphism together with a distinguished invertible element, $(\otimes, m_T)$;
\item a weak functor $I:1 \rightarrow T(\star, \star)$, which will be considered as the unique monoid homomorphism $1 \rightarrow T$ together with a distinguished invertible element $u_{T}$;
\item an adjoint equivalence $\mathbf{a}: \otimes \circ \otimes \times 1 \Rightarrow \otimes \circ 1 \times \otimes$, which is the assertion that $\otimes$ is strictly associative as a binary operation on $T$ together with a distinguished invertible element $a_{T}$;
\item adjoint equivalences $\mathbf{l}: \otimes \circ I \times 1 \Rightarrow 1, \mathbf{r}: \otimes \circ 1 \times I \Rightarrow 1$, which is the assertion that $1$ is a unit for $\otimes$ as a binary operation on $T$, together with distinguished invertible elements $l_{T}, r_{T}$;
\item and four distinguished invertible elements $\pi_{T}, \mu_{T}, \lambda_{T}, \rho_{T}$.
\end{itemize}

Thus we have a commutative monoid $T$, a monoid homomorphism 
	\[\otimes:T \times T \rightarrow T,\] 
and distinguished invertible elements $d_{T}, m_T, u_{T}, a_{T}, l_{T}, r_{T}, \pi_{T}, \mu_{T}, \lambda_{T}, \rho_{T}$.  The fact that $\otimes$ is a monoid homomorphism is expressed in the following equation, where we have written the monoid structure on $T$ as concatenation.
\[
(ab) \otimes (cd) = (a \otimes c)(b \otimes d)
\]
The adjoint equivalences $\mathbf{l}, \mathbf{r}$ each imply that 1 is a unit for $\otimes$.  Using this and the equation above, the Eckmann-Hilton argument immediately implies that $a \otimes b = ab$.

We will later need to use the naturality isomorphisms; it is simple to compute that that the naturality isomorphism for the transformation $a$ is 1, and the naturality isomorphisms for $l$ and $r$ are both $m_T$.

There are three tricategory axioms that we must now check to find the dependence between distinguished invertible elements.  Using the above, it is  straightforward to check that the first tricategory axiom is vacuous, the second gives the equation
\[
\lambda \pi = d^{2}m_{T}^{4},
\]
and the third gives the equation
\[
\rho \pi = d^{2}m_{T}^{4}.
\]
Since $\lambda, \rho, \pi$, and $d$ are invertible, 
\[
\lambda = \rho = \pi^{-1} d^{2}m_{T}^{4}.
\]
Thus $\lambda$ and $\rho$ are determined by the remaining data, hence we have the result as summarised above.

\subsection{Weak functors}
In this section we characterise weak functors between triply degenerate tricategories.  {\it A priori} a weak functor $F:S \rightarrow T$ between triply degenerate tricategories consists of the following data, which we will try to simplify:
\begin{itemize}
\item a weak functor $F_{\star, \star}:S(\star, \star) \rightarrow T(\star, \star)$, which by the results of \cite{cg1} is a monoid homomorphism $F:S \rightarrow T$ together with a distinguished invertible element $m_{F} \in T$;
\item an adjoint equivalence $\mathbf{\chi}: \otimes' \circ (F \times F) \Rightarrow F \circ \otimes$, which is the trivial assertion that $F(a \otimes b) = Fa \otimes ' Fb$ together with a distinguished invertible element $\chi \in T$;
\item an adjoint equivalence $\mathbf{\iota}: I'_{\star} \Rightarrow F \circ I_{\star}$, which is the trivial assertion that $F1 = 1$ together with a distinguished invertible element $\iota \in T$;
\item and invertible modifications $\omega, \gamma,$ and $\delta$.
\end{itemize}
Thus we have a monoid homomorphism $F$ and six distinguished invertible elements $m_{F}, \chi, \iota, \omega, \gamma$, and $\delta$.  It is straightforward to compute that the naturality isomorphism for $\chi$ is given by the invertible element $Fm_{S} \cdot (m_{T}m_{F})^{-1}$ and the naturality isomorphism for $\iota$ is given by $m_{F}$.

There are two axioms for weak functors for tricategories.  In the case of triply degenerate tricategories, the first axiom reduces to the equation
\[
\omega \cdot \pi_{T} \cdot  Fm_{S}^{2} \cdot m_{T}^{-2} \cdot Fd_{S}^{2} \cdot d_{T}^{-2} = F\pi_{S}
\]
thus by invertibility $\omega$ is determined by the rest of the data.  The second axiom reduces to the equation
\[
\omega  \cdot \delta  \cdot \gamma  \cdot \mu_{T}  \cdot Fm_{S}^{2}  \cdot m_{T}^{-2}  \cdot Fd_{S}^{2}  \cdot d_{T}^{-2} = F \mu_{S}.
\]
By the previous equation and the invertibility of all terms involved, $\delta$ and $\gamma$ determine each other once the rest of the data is fixed, hence we have the result as summarised above.

\subsection{Tritransformations}
In this section we characterise tritransformations for triply degenerate tricategories.  First we need the following lemma, which is a simple calculation.
\begin{lem}
Let $T$ be a triply degenerate tricategory.  Then the functor
\[
T(1, I_{\star}) = I_{\star} \circ - :T(\star, \star) \rightarrow T(\star, \star)
\]
is given by the identity homomorphism together with the distinguished invertible element $d^{-1}m$.  Additionally, $T(1, I_{\star}) = T(I_{\star}, 1)$.
\end{lem}

{\it A priori}, the data for a tritransformation $\alpha: F \rightarrow G$ of triply degenerate tricategories consists of:
\begin{itemize}
\item an adjoint equivalence $\bs{\alpha}: T(1,I_{\star}) \circ F \Rightarrow T(I_{\star}, 1) \circ G$, which consists of the assertion that $F=G$ as monoid homomorphisms together with a distinguished invertible element $\alpha_{T}$; and
\item distinguished invertible elements $\Pi$ and $M$.
\end{itemize}

It is easy to compute that the naturality isomorphism for the transformation $\alpha$ is $m_{F}^{-1}m_{G}$.  The first transformation axiom reduces to the equation
\[
m_{G} = m_{F},
\]
the second axiom reduces to the equation
\[
\Pi \mu_{T} l_{T} \gamma_{F} = M m_{T}^{4} d_{T}^{2}a_{T}^{-1} \gamma_{G},
\]
and the third to the equation
\[
\Pi \delta_{F} = a_{T}^{-1} l_{T}^{-1} d_{T}^{2}m_{T}^{4} \mu^{-1} M \delta_{G}.
\]
Thus we see that $\Pi$ determines $M$, and that the second and third axioms combine to yield no new information. So we have remaining distinguished invertible elements $\Pi$ and $\alpha_{T}$, giving the results as summarised above.

\subsection{Trimodifications and perturbations}

The data for a trimodification $m:\alpha \Rightarrow \beta$ consists of a single invertible element $m$ in $T$, and there are two axioms.  The first is the equation
\[
m^{2} \cdot \Pi \cdot Gd_{S} = \Pi \cdot Fd_{S} \cdot m
\]
which reduces to $m=1$ since $F=G$ as monoid homomorphisms.  The second axiom also reduces to $m=1$, thus there is a unique trimodification between any two parallel transformations.  Note that this means that \textit{any} diagram of trimodifications in this setting commutes, a fact that will be useful later.

The data for a perturbation $\sigma: m \Rrightarrow n$ consists of an element $\sigma$ in $T$.  The single axiom is vacuous so a perturbation is precisely an element $\sigma \in T$.

\subsection{Overall structure}

We now compare the totalities of, on the one hand triply degenerate tricategories, and on the other hand commutative monoids.  Recall that for the case of doubly degenerate bicategories we were able to attempt comparisons at the level of categories, bicategories and tricategories of such, simply by truncating the full sub-tricategory of \cat{Bicat} to the required dimension.  However, for triply degenerate tricategories we show that truncating the full sub-tetracategory of \cat{Tricat} does not yield a category or a bicategory; truncation does yield a tricategory, and this is the only level that yields an equivalence with commutative monoids.  As in \cite{cg1} we compare with the discrete $j$-categories of commutative monoids obtained by adding higher identity cells to \cat{CMon}.  

Note that we do not actually prove that we have a tetracategory of triply degenerate tricategories; for the comparison, we simply prove that the natural putative functor is not full and faithful therefore cannot be an equivalence.

We have a 4-dimensional structure with

\begin{itemize}
\item 0-cells:  triply degenerate tricategories
\item 1-cells:  weak functors between them
\item 2-cells:  tritransformations between those
\item 3-cells:  trimodifications between those
\item 4-cells:  perturbations between those.
\end{itemize}

We write $\mathbf{Tricat}(3)_{j}$ for the truncation of this structure to a $j$-dimensional structure, and $\mathbf{CMon}_{j}$ for the $j$-category of commutative monoids and their morphisms (and higher identities where necessary).

There are obvious assignments 
\[\begin{tabular}{rcl}
triply degenerate tricategory  & $\mapsto$ & underlying commutative monoid \\
weak functor & $\mapsto$ & underlying homomorphism of monoids
\end{tabular}\]
which, together with the unique maps on higher cells, form the underlying morphism on $j$-globular sets for putative functors 
\[
\xi_{j}: \mathbf{Tricat}(3)_{j} \rightarrow \mathbf{CMon}_{j}.
\]

\begin{thm}\label{triply}\mbox{\hspace{10cm}}
\begin{enumerate}
\item $\mathbf{Tricat}(3)_1$ is not a category.
\item $\mathbf{Tricat}(3)_2$ is not a bicategory.
\item $\mathbf{Tricat}(3)_3$ is a tricategory, and $\xi_3$ defines a functor which is a triequivalence.
\item $\xi_4$ does not give a tetra-equivalence of tetra-categories.
\end{enumerate}
\end{thm}

The rest of this section will constitute a gradual proof of the various parts of this theorem.  We begin by constructing the hom-bicategories for a tricategory structure on $\mathbf{Tricat}(3)_{3}$.

\begin{prop}
Let $X,Y$ be triply degenerate tricategories.  Then there is a bicategory $\mathbf{Tricat}(3)_{3}(X,Y)$ with 0-cells weak functors $F:X \rightarrow Y$, 1-cells tritransformations $\alpha:F \Rightarrow G$, and 2-cells trimodifications $m: \alpha \Rrightarrow \beta$.
\end{prop}

\begin{proof}
To give the bicategory structure, we need only provide unit 1-cells and 1-cell composition since there is a unique trimodification between every pair of parallel tritransformations.  It is simple to read off the required distinguished invertible elements from the corresponding formulae for composites of tritransformations and from the data for the unit tritransformation.
\end{proof}

\begin{rem}
Note that composition of 1-cells in $\mathbf{Tricat}(3)_{3}(X,Y)$ is strictly associative, but is not strictly unital.  In particular, this shows that $\mathbf{Tricat}(3)_{2}$ is not a bicategory, proving Theorem~\ref{triply}, part 2.
\end{rem}

We now construct the composition functor
\[
\otimes: \mathbf{Tricat}(3)_{3}(Y,Z) \times \mathbf{Tricat}(3)_{3}(X,Y) \rightarrow \mathbf{Tricat}(3)_{3}(X,Z).
\]
for any triply degenerate tricategories $X,Y,Z$.  We define the composite $GF$ of functors $F:X \rightarrow Y$, $G:Y \rightarrow Z$ by the following formulae which can be read off directly from the formulae giving the composite of functors between tricategories.
\[
\begin{array}{rcl}
m_{GF} & = & m_{G}Gm_{F} \\
\chi_{GF} & = &  \chi_{G}G(\chi_{F}d_{Y})d_{Z}^{-2} \\
\iota_{GF} & = &  \iota_{G}G(\iota_{F}d_{Y})d_{Z}^{-2} \\
\gamma_{G} & = & d_{Z}^{-2}m_{Z}^{2}m_{G}^{2}\gamma_{G}G(\gamma_{F}d_{Y}m_{Y})
\end{array}
\]
The formulae for the composite $\beta \otimes \alpha$ of two transformations are derived similarly, and thus we have a weak functor $\otimes$ for composition as required.

Similarly, there is a unit functor
\[
I_{X}: 1 \rightarrow  \mathbf{Tricat}(3)_{3}(X,X)
\]
whose value on the unique 0-cell is the identity functor on $X$.

\begin{rem}
The formulae above make it obvious that $\otimes$ is not strictly associative on 0-cells, and that the identity functor is not a strict unit for $\otimes$.  This shows that $\mathbf{Tricat}(3)_{1}$ is not a category, proving Theorem~\ref{triply}, part 1.
\end{rem}

Next we need to specify the required constraint adjoint equivalences.  It is straightforward to find adjoint equivalences
\[
\begin{array}{c}
\bs{\mathcal{A}}: \otimes \circ \otimes \times 1 \Rightarrow \otimes \circ 1 \times \otimes \\
\bs{\mathcal{L}}: \otimes \circ I \times 1 \Rightarrow 1 \\
\bs{\mathcal{R}}: \otimes \circ 1 \times I \Rightarrow 1
\end{array}
\]
in the appropriate functor bicategories; the actual choice of adjoint equivalence is irrelevant, since there is a unique modification between any pair of parallel transformations.



Finally, to finish constructing the tricategory $\mathbf{Tricat}(3)_{3}$ we must define invertible modifications $\pi, \mu, \lambda, \rho$ and check three axioms. However since there are unique trimodifications between parallel tritransformations, these modifications are uniquely determined and the axioms automatically hold.

We now examine the morphism $\xi_3$ of 3-globular sets and show that it defines a functor
\[\mathbf{Tricat}(3)_{3} \lra \mathbf{CMon}_{3};\]
in fact functoriality is trivial as $\mathbf{CMon}_{3}$ has discrete hom-2-categories.  Furthermore we show it is an equivalence as follows.  The functor is clearly surjective on objects, and the functor on hom-bicategories
\[
\mathbf{Tricat}(3)_{3}(X,Y) \rightarrow \mathbf{CMon}_{3}(\xi_{3}X, \xi_{3}Y)
\]
is easily seen to be surjective on objects as well.  This functor on hom-bicategories is also a local equivalence since $\mathbf{CMon}_{3}$ is discrete at dimensions two and three and $\mathbf{Tricat}(3)_{3}$ has unique 3-cells between parallel 2-cells.  This finishes the proof of Theorem~\ref{triply}, part 3.  

For part 4, we observe that the morphism $\xi_{4}$ of 4-globular sets is clearly not locally faithful on 4-cells.  This finishes the proof of Theorem~\ref{triply}.  

\section{Doubly degenerate tricategories}\label{2degen}

We now compare doubly degenerate tricategories with braided monoidal categories.  As described informally in the Introduction the comparison is not straightforward.  Therefore we begin by directly listing the structure that we get on the monoidal category given by the (unique) degenerate hom-bicategory; this is simply a matter of writing out the definitions as nothing simplifies in this case.  Afterwards, we show how to extract a braided monoidal category from this structure.  Essentially, all of the data listed in Section~\ref{doublybasic} can be thought of as ``extra structure'' that arises on the braided monoidal category we will construct.  Finally, we construct some comparison functors.

We will begin with an informal overview of this whole section as we feel that for many readers the ideas will be at least as important as the technical details.  

\subsection{Overview}\label{overview}

It is widely accepted that a doubly degenerate bicategory ``is'' a commutative monoid, and that a doubly degenerate tricategory ``is'' a braided monoidal category.  Moreover, it is widely accepted that the proof of the bicategory case is ``simply'' a question of applying the Eckmann-Hilton argument to the multiplications given by horizontal and vertical composition, and that the tricategory result is proved by doing this process up to isomorphism.  In this section we give an informal overview of the extent to which this is and is not the case.  We believe that this is important because the disparity will increase as dimensions increase, and because this issue seems to lie at the heart of various critical phenomena in higher-dimensional category theory, such as:

\begin{enumerate}

\item why we do not expect every weak $n$-category to be equivalent to a strict one

\item why weak $n$-categories are expected to model homotopy $n$-types while strict ones are known not to do so \cite{gro1, bd3, sim3}

\item why some diagrams of constraints in a tricategory do not in general commute, and why these do not arise in free tricategories \cite{gur1}

\item why strict computads do not form a presheaf category \cite{mz1}

\item why the existing definitions of $n$-categories based on reflexive globular sets fail to be fully weak \cite{cl1}

\item why a notion of semistrict $n$-category with weak units but strict interchange may be weak enough to model homotopy $n$-types and give coherence results \cite{sim1, koc1, jk1}

\item why we need weak $n$-categories at all, and not just strict ones.

\end{enumerate}

A doubly degenerate bicategory \cl{B} has only one 0-cell $\star$ and only one 1-cell $I_\star$.  To show that the 2-cells form a commutative monoid we first use the fact that they are the morphisms of the single hom-category $\cl{B}(\star, \star)$; since this hom-category has only one object $I_\star$ we know it is a monoid, with multiplication given by vertical composition of 2-cells.  To show that it is a {\it commutative} monoid, we apply the Eckmann-Hilton argument to the two multiplications defined on the set of 2-cells: vertical composition and horizontal composition.  

Recall that the Eckmann-Hilton argument says: Let $A$ be a set with two binary operations $\ast$ and $\circ$ such that

\begin{enumerate}
\item  $\ast$ and $\circ$ are unital with the same unit 
\item $\ast$ and $\circ$ distribute over each other i.e. $\forall a,b,c,d \in A$ 
	\[(a \ast b)\circ (c \ast d) = (a \circ c)\ast(b \circ d).\]
\end{enumerate}
Then $\ast$ and $\circ$ are in fact equal and this operation is commutative.  Note that the two binary operations are usually called products (with implied associativity) but in fact associativity is irrelevant to the argument.

However, in our case a difficulty arises because horizontal composition in a bicategory is not strictly unital.  The situation is rescued by the fact that $l_I = r_I$ in any bicategory.  This, together with the naturality of $l$ and $r$, enables us to prove, albeit laboriously, that horizontal composition is strictly unital for 2-cells in a {\it doubly degenerate} bicategory, and moreover that the vertical 2-cell identity also acts as a horizontal identity.  Thus we can in fact apply the Eckmann-Hilton argument.

Generalising this argument to doubly degenerate tricategories directly is tricky.  There are various candidates for a ``categorified Eckmann-Hilton argument'' provided by Joyal and Street \cite{js1, bat2}.  The idea is to replace all the equalities in the argument by isomorphisms, but as usual we need to take some care over {\it specifying} these isomorphisms rather than merely asserting their existence; see Definition~\ref{defmult}.  

However, when we try and apply this result to a doubly degenerate tricategory we have some further difficulties: composition along bounding 0-cells is difficult to manipulate as a multiplication, because we cannot use coherence results for tricategories.  Coherence for tricategories \cite{gur2} tells us that ``every diagram of constraints in a free tricategory (on a category-enriched 2-graph) commutes''.  In particular this means that if we need to use cells that do not arise in a free tricategory, then we cannot use coherence results to check axioms.  This is the case if we attempt to built a multiplication out of composition along 0-cells; we have to use the fact that we only have one 1-cell in our tricategory, and therefore that various composites of 1-cells are all ``accidentally'' the same. This comes down to the fact that the free tricategory on a doubly degenerate tricategory is not itself doubly degenerate; it is not clear how to construct a ``free doubly degenerate tricategory''.

However, to rectify this situation we can look at an alternative way of proving the result for degenerate bicategories, that does not make such identifications.  We still use the Eckmann-Hilton argument but instead of attempting to apply it using horizontal composition of 2-cells, we define a new binary operation on 2-cells that is derived from horizontal composition as follows:

\begin{displaymath}
\beta \odot \alpha = r \circ (\beta * \alpha) \circ l^{-1} 
\end{displaymath}
(Essentially this is what we use to prove that horizontal composition is strictly unital in the previous argument.)  Unlike horizontal composition, this operation does ``categorify correctly'', that is, given a doubly degenerate tricategory we can define a multiplication on its associated monoidal category by using the above formula (this is the content of Theorem~\ref{degenmults}), and we can manipulate it using coherence for tricategories. 

To extract a braiding from this we then have to follow the steps of the Eckmann-Hilton argument and keep track of all the isomorphisms used; this is Proposition~\ref{jsprop}.

We see that we use instances of the following cells, in a lengthy composite:

\begin{itemize}

\item naturality constraints for $l_I$ and $r_I$
\item constraints for weak interchange of 2-cells
\item isomorphisms showing that $l_I \cong r_I$

\end{itemize}

This shows very clearly why a theory with weak units but strict interchange is enough to produce braidings -- the braiding is built from all of the above structure contraints, so if any one of them is weak then braidings will still arise.  Thus if we have strict units then we need weak interchange, but if we maintain weak units we can have strict interchange and still get a braided monoidal category.  As mentioned above we do, however, get a certain amount of extra structure on the braided monoidal category that arises, and there does not seem to be a straightforward way of organising it, or of describing coherently the tricategorical situation in which this extra structure is trivial. 

\subsection{Basic results}\label{doublybasic}

The results in this section are all obtained by simply rewriting the appropriate definitions using the results of \cite{cg1}.    Many of the diagrams needed in the theorems below are excessively large, and have been relegated to the Appendix.

Just as we began the previous section by characterising adjoints in doubly degenerate bicategories, we begin this section by recalling the definition of ``dual pair'' of objects in a monoidal category, since this characterises adjoints for 1-cells in degenerate bicategories; eventually we will of course be interested in adjoint equivalences, not just adjoints.

\begin{Defi}
Let $M$ be a monoidal category.  Then a \emph{dual pair} in $M$ consists of a pair of objects $X, X^{\cdot}$ together with morphisms $\varepsilon: X \otimes X^{\cdot} \rightarrow I, \eta:I \rightarrow X^{\cdot} \otimes X$ satisfying the two equations below, where all unmarked isomorphisms are given by coherence isomorphisms.
\[
\xy
{\ar^{\cong} (0,0)*+{X}; (20,0)*+{XI} };
{\ar^{1 \eta} (20,0)*+{XI}; (40,0)*+{X(X^{\cdot}X)} };
{\ar^{\cong} (40,0)*+{X(X^{\cdot}X)}; (60,0)*+{(XX^{\cdot})X} };
{\ar^{\varepsilon 1} (60,0)*+{(XX^{\cdot})X}; (80,0)*+{IX} };
{\ar^{\cong} (80,0)*+{IX}; (80,-20)*+{X} };
{\ar_{1} (0,0)*+{X}; (80,-20)*+{X} }
\endxy
\]
\[
\xy
{\ar^{\cong} (0,0)*+{X^{\cdot}}; (20,0)*+{IX^{\cdot}} };
{\ar^{\eta 1} (20,0)*+{IX^{\cdot}}; (40,0)*+{(X^{\cdot}X)X^{\cdot}} };
{\ar^{\cong} (40,0)*+{(X^{\cdot}X)X^{\cdot}}; (60,0)*+{X^{\cdot}(XX^{\cdot})} };
{\ar^{1 \varepsilon} (60,0)*+{X^{\cdot}(XX^{\cdot})}; (80,0)*+{X^{\cdot}I} };
{\ar^{\cong} (80,0)*+{X^{\cdot}I}; (80,-20)*+{X^{\cdot}} };
{\ar_{1} (0,0)*+{X^{\cdot}}; (80,-20)*+{X^{\cdot}} }
\endxy
\]
\end{Defi}

\begin{thm}
A doubly degenerate tricategory $B$ is precisely
\begin{itemize}
\item a monoidal category $(B, \otimes, U, a, l, r)$ given by the single degenerate hom-bicategory;
\item a monoidal functor $\boxtimes: B \times B \rightarrow B$ from composition;
\item a monoid $I$ in $B$ and an isomorphism $I \cong U$ as monoids in $B$; this comes from the functor for units $I \lra B(\star, \star)$
\item a dual pair $(A, A^{\cdot}, \varepsilon_{A}, \eta_{A})$ with $\varepsilon_{A}, \eta_{A}$ both invertible, and natural isomorphisms
\[
\begin{array}{c}
A \otimes \Big( (X \boxtimes Y) \boxtimes Z) \big) \cong \Big( X \boxtimes (Y \boxtimes Z) \Big) \otimes A \\
A^{\cdot} \otimes  \Big( X \boxtimes (Y \boxtimes Z) \Big) \cong \Big( (X \boxtimes Y) \boxtimes Z) \big) \otimes A^{\cdot};
\end{array}
\]

subject to diagrams given in the Appendix,

\item a dual pair $(L, L^{\cdot}, \varepsilon_{L}, \eta_{L})$ with with $\varepsilon_{L}, \eta_{L}$ both invertible, and natural isomorphisms

\[
\begin{array}{c}
L \otimes (I \boxtimes X) \cong X \otimes L \\
L^{\cdot} \otimes X \cong (I \boxtimes X) \otimes L^{\cdot}
\end{array}
\]
subject to diagrams given in the Appendix,
\item a dual pair $(R, R^{\cdot}, \varepsilon_{R}, \eta_{R})$ with with $\varepsilon_{L}, \eta_{L}$ both invertible, and natural isomorphisms
\[
\begin{array}{c}
R \otimes (X \boxtimes I) \cong X \otimes R \\
R^{\cdot} \otimes X \cong (X \boxtimes I) \otimes R^{\cdot};
\end{array}
\]
subject to diagrams given in the Appendix,
\item and isomorphisms
\[
\begin{array}{c}
\Bigg( (U \boxtimes A) \otimes \big( A \otimes (A \boxtimes U) \big) \Bigg) \stackrel{\pi}{\cong} A \otimes A \\
\Bigg( (U \boxtimes L) \otimes \big( A \otimes (R^{\cdot} \boxtimes U) \big) \Bigg) \stackrel{\mu}{\cong} U \\
L \boxtimes U \stackrel{\lambda}{\cong} L \otimes A \\
U \boxtimes R^{\cdot} \stackrel{\rho}{\cong} A \otimes R^{\cdot};
\end{array}
\]
\end{itemize}
all subject to three axioms which appear in the Appendix.
\end{thm} 

\begin{rem}
It is important to note that $\boxtimes$ does not {\it a priori} give a monoidal structure on the category $B$; the obstruction is that lax transformations between weak functors of degenerate tricategories are more general than monoidal transformations between the associated monoidal functors (see \cite{cg1}).  As noted in Section~\ref{overview} it may be possible to prove that $\boxtimes$ is a valid monoidal structure, but since we cannot use coherence for tricategories to help us, the proof is not very evident.  Thus to extract a braiding from all this structure, we will not simply apply an Eckmann-Hilton-style argument to $\otimes$ and $\boxtimes$ (see Section~\ref{braidings}).
\end{rem}

We now describe functors, transformations, modifications and perturbations in a similar spirit.

\begin{thm}
A weak functor $F: B \rightarrow B'$ between doubly degenerate tricategories is precisely
\begin{itemize}
\item a monoidal functor $F:B \rightarrow B'$;
\item a dual pair $(\chi, \chi^{\cdot}, \varepsilon_{\chi}, \eta_{\chi})$ in $B'$ with $\varepsilon_{\chi}, \eta_{\chi}$ both invertible, and natural isomorphisms
\[
\begin{array}{c}
\chi \otimes' (FX \boxtimes' FY) \cong F(X \boxtimes Y) \otimes' \chi \\
\chi^{\cdot} \otimes' F(X \boxtimes Y) \cong (FX \boxtimes' FY) \otimes' \chi^{\cdot}
\end{array}
\]
subject to diagrams given in the Appendix,
\item a dual pair $(\iota, \iota^{\cdot}, \varepsilon_{\iota}, \eta_{\iota})$ with $\varepsilon_{\iota}, \eta_{\iota}$ both invertible, and natural isomorphisms
\[
\begin{array}{c}
\iota \otimes' I' \cong FI \otimes' \iota \\
\iota^{\cdot} \otimes' FI \cong I' \otimes' \iota^{\cdot}
\end{array}
\]
subject to diagrams given in the Appendix,
\item and isomorphisms
\[
\begin{array}{c}
FA \otimes' \Big( \chi \otimes' (\chi \boxtimes' U') \Big) \stackrel{ \omega}{\cong} \chi \otimes' \Big( (U' \boxtimes' \chi) \otimes' A' \Big) \\
FL \otimes' \Big( \chi \otimes' (\iota \boxtimes' U') \Big) \stackrel{\gamma}{\cong} L' \\
FR^{\cdot} \stackrel{\delta}{\cong} \chi \otimes' \Big( (U' \boxtimes' \iota) \otimes' (R')^{\cdot} \Big);
\end{array}
\]
\end{itemize}
all subject to axioms given in the Appendix.
\end{thm}

\begin{thm}\label{ddtrans}
A weak transformation $\alpha:F \rightarrow G$ in the above setting is precisely
\begin{itemize}
\item a dual pair $(\alpha, \alpha^{\cdot}, \varepsilon_{\alpha}, \eta_{\alpha})$ with $\varepsilon_{\alpha}, \eta_{\alpha}$ both invertible, and natural isomorphisms
\[
\begin{array}{c}
\alpha \otimes' (U' \boxtimes' FX) \cong (GX \boxtimes' U') \otimes' \alpha \\
\alpha^{\cdot} \otimes' (GX \boxtimes' U') \cong (U' \boxtimes' FX) \otimes' \alpha^{\cdot}
\end{array}
\]
subject to diagrams given in the Appendix,
\item and isomorphisms
\[
\begin{array}{c}
\mbox{\hspace*{-4em}}(\chi_{G} \otimes' U') \otimes' \Bigg( (A')^{\cdot} \otimes' \Big( (U' \boxtimes' \alpha) \otimes' \big(A' \otimes' (\alpha \boxtimes' U') \big) \Big) \Bigg) \stackrel{\Pi}{\cong} \alpha \otimes' \Big( (U' \boxtimes' \chi_{F}) \otimes' A' \Big) \\
\mbox{\hspace*{-4em}}\alpha \otimes' \Big( (U' \boxtimes' \iota_{F}) \otimes' (R')^{\cdot} \Big) \stackrel{M}{\cong} (\iota_{G} \boxtimes' U') \otimes' (L')^{\cdot};
\end{array}
\]
\end{itemize}
all subject to three axioms given in the Appendix.
\end{thm}

The analogous result for lax transformations should be obvious, with dual pair replaced by distinguished object since in the lax case we have a noninvertible morphism instead of an adjoint equivalence.

\begin{thm}
A modification $m:\alpha \Rightarrow \beta$ is precisely
\begin{itemize}
\item an object $m \in B'$ and
\item an isomorphism
\[
(U' \boxtimes m) \otimes' \alpha \cong \beta \otimes' (m \boxtimes' U')
\]
\end{itemize}
subject to two axioms given in the Appendix.
\end{thm}

\begin{thm}
A perturbation $\sigma: m \Rrightarrow n$ is precisely a morphism $\sigma: m \rightarrow n$ in $B'$ satisfying the single axiom in the Appendix. 
\end{thm}

\subsection{Braidings}\label{braidings}

In this section we show that the underlying monoidal category of a doubly degenerate tricategory does have a braiding on it.  To show this, we use the fact that to give a braiding for a monoidal structure, it suffices to give the structure of a multiplication on the monoidal category in question.  We give the relevant definitions below; for additional details, see \cite{js1}.

\begin{Defi}\label{defmult}
Let $M$ be a monoidal category, and equip $M \times M$ with the componentwise monoidal structure.  Then a \textit{multiplication} $\varphi$ on $M$ consists of a monoidal functor $\varphi:M \times M \rightarrow M$ and invertible monoidal transformations $\rho: \varphi \circ (\textrm{id} \times I) \Rightarrow \textrm{id}, \lambda: \varphi \circ (I \times \textrm{id}) \Rightarrow \textrm{id}$ where $I:1 \rightarrow M$ is the canonical monoidal functor whose value on the single object is the unit of $M$ and whose structure constraints are given by unique coherence isomorphisms.
\end{Defi}

The following result, due to Joyal and Street \cite{js1}, says that a multiplication naturally gives rise to a braiding.

\begin{prop}\label{jsprop}
Let $M$ be a monoidal category with multiplication $\varphi$.  Then $M$ is braided with braiding given by the composite below.
\[\mbox{\hspace*{-8em}}
ab \stackrel{\lambda^{-1} \rho^{-1}}{\longrightarrow} \varphi(I, a)\varphi(b,I) \stackrel{\cong}{\longrightarrow} \varphi(Ib, aI) \stackrel{\varphi(l, r)}{\longrightarrow} \varphi(b,a)\] 
\[\mbox{\hspace*{8em}}\stackrel{\varphi(r^{-1}, l^{-1})}{\longrightarrow} \varphi(bI, Ia) \stackrel{\cong}{\longrightarrow} \varphi(b,I)\varphi(I,a) \stackrel{\rho \lambda}{\longrightarrow} ba
\]
\end{prop}

We will use this construction to provide a braiding for the monoidal category associated to a doubly degenerate tricategory.  As can be seen from the above formula, this braiding is ``natural'' but not exactly ``simple''.  

\begin{thm}\label{degenmults}
Let $B$ be a doubly degenerate tricategory, and also denote by $B$ the monoidal category associated to the single (degenerate) hom-bicategory.  Then there is a multiplication $\varphi$ on $B$ with 
	\[\varphi(X,Y) = R \otimes \big( (X \boxtimes Y) \otimes L^{\cdot} \big).\]
\end{thm}
This result is a lengthy calculation that requires repeated use of the coherence theorm for tricategories as well as coherence for bicategories and functors.  We thus omit it, and only record the following crucial corollary.

\begin{cor}\label{degenbraids}
Let $B$ be a doubly degenerate tricategory, and also denote by $B$ the monoidal category associated to the single (degenerate) hom-bicategory.  Then $B$ is a braided monoidal category.
\end{cor}

The situation for functors is similar, with braided monoidal functors arising from ``multiplicative'' functors as follows.  

\begin{Defi}
Let $(M, \varphi)$ and $(N, \psi)$ be monoidal categories equipped with multiplications.  A \textit{multiplicative functor} $F: (M, \varphi) \rightarrow (N, \psi)$ consists of a monoidal functor $F: M \rightarrow N$ and an invertible monoidal transformation $\chi: \psi \circ (F \times F) \Rightarrow F \circ \phi$, satisfying unit axioms.
\end{Defi}

\begin{prop}
Let $(M, \varphi)$ and $(N, \psi)$ be monoidal categories equipped with multiplications, and let  $F: (M, \varphi) \rightarrow (N, \psi)$ be a multiplicative functor between them.  Then the underlying monoidal functor $F$ is braided when $M$ and $N$ are equipped with the braidings induced by their respective multiplications.
\end{prop}

The following theorem says that functors between doubly degenerate tricategories do give rise to multiplicative functors, and as a corollary, braided monoidal functors.  The proof of the theorem is another long calculation involving coherence.

\begin{thm}
Let $B$ and $B'$ be doubly degenerate tricategories, and let $F: B \rightarrow B'$ be a functor between them.  Then the monoidal functor $F$ between the monoidal categories $B$ and $B'$ can be given the structure of a multiplicative functor when we equip $B$ and $B'$ with the multiplications of Theorem \ref{degenmults}.
\end{thm}

\begin{cor}\label{degenmorph}
Let $B$ and $B'$ be doubly degenerate tricategories, and let $F: B \rightarrow B'$ be a functor between them.  Then the monoidal functor $F$ is braided with respect to the braided monoidal categories $B$ and $B'$ as in Corollary \ref{degenbraids}.
\end{cor}

The situation for transformations does not lend itself to the same sort of analysis: a transformation of doubly degenerate tricategories is rather different from a monoidal transformation.  This also occurs in the study of degenerate bicategories, where transformations of degenerate bicategories are rather different from monoidal transformations.  Thus, as discussed in the Introduction, we modify \cat{Tricat} so that the 2-cells between doubly degenerate tricategories {\it do} give rise to monoidal transformations \cite{gg1}.

%
%
%

We construct a bicategory $\widetilde{\mathbf{Tricat}}$ with
\begin{itemize}
\item 0-cells: tricategories, 
\item 1-cells: functors, and
\item 2-cells: ``locally iconic lax transformations'', that is, transformations whose 1- and 2-cell components are identities.
\end{itemize}
Recall that an icon is an oplax transformation of bicategories, all of whose component 1-cells are identities.  Specifying a 2-cell of $\widetilde{\mathbf{Tricat}}$ thus involves specifying, for every 2-cell in the source, a 3-cell component in the target, together with 3-cell constraint data $\Pi$ and $M$ (as exhibited in Theorem~\ref{weaktransform} for example) satisfying the relevant axioms; there is effectively no lower-dimensional data.

Since icons between degenerate bicategories yield monoidal transformations (see \cite{cg1}), the following result is immediate.

\begin{thm}\label{degen2cells}
Let $\alpha:F \Rightarrow G$ be a 2-cell in $\widetilde{\mathbf{Tricat}}$ whose source and target tricategories are doubly degenerate.  Then $\alpha$ gives rise to a monoidal transformation between the braided monoidal functors corresponding to $F$ and $G$.
\end{thm}
This gives us a slightly better comparison with braided monoidal categories---we can at least have functors in both directions---but in the next section we will compare the overall structures and see that we still do not get a biequivalence of bicategories.

\subsection{Overall structure}

In this section we attempt to compare the totalities of structures involved.  That is, on the one hand we have doubly degenerate tricategories, and on the other hand we have braided monoidal categories as predicted by the Periodic Table.  

First observe that the full 4-dimensional structure of tricategories does not yield an equivalence.  We can add higher identity cells to the bicategory \cat{BrMonCat} of braided monoidal categories to form a discrete tetracategory, but it is clear that this cannot be equivalent to the tetracategory $\cat{Tricat}(2)$ of doubly degenerate tricategories, which has too many non-trivial 4-cells.  

Instead we can examine $\widetilde{\cat{Tricat}}(2)$, the full sub-bicategory of $\widetilde{\cat{Tricat}}$ whose 0-cells are the doubly degenerate tricategories. We will now show that there are naturally arising comparison functors to and from \cat{BrMonCat}, but these do not exhibit a biequivalence.

From the results of the previous section we have a morphism of underlying globular sets
	\[U:\widetilde{\cat{Tricat}}(2) \lra \cat{BrMonCat},\]
assigning
\begin{itemize}
\item to each 0-cell associated braided monoidal category as in Corollary~\ref{degenbraids}, 
\item to each 1-cell the associated braided monoidal functor as in Corollary~\ref{degenmorph}, and
\item to each 2-cell the associated monoidal transformation as in Theorem~\ref{degen2cells}.
\end{itemize}
	
\begin{prop}
The morphism $U$ defines a strict functor of bicategories.
\end{prop}
	
\begin{proof}
It is trivial that $U$ strictly preserves composition and units; the functor axioms follow by noting that $U$ sends the constraint isomorphisms $a, l, r$ to identities.
\end{proof}
	
There is also a comparison functor in the opposite direction
	\[F:\cat{BrMonCat} \lra \widetilde{\cat{Tricat}}(2).\]
This is simply a matter of choosing all the extra structure to be given by identities where this makes sense, and using isomorphisms given canonically by coherence constraints elsewhere.  We can then check the following theorem.

\begin{thm}\label{idcomp}
The composite functor $UF:\mathbf{BrMonCat} \rightarrow \mathbf{BrMonCat}$ is the identity 2-functor.
\end{thm}
\begin{proof}
It is obvious that $UF(X)$ has the same monoidal structure as $X$, and following the definition of the braiding carefully and using coherence gives that the braided structures are the same as well.  It follows immediately that $UF$ is the identity on 1- and 2-cells.  It remains to check that the constraints for this composite are also identities; this follows from the definition of the constraints for $F$ and the definition of $U$.
\end{proof}

We can now show that these comparison functors do not exhibit an equivalence.  The problem is at the level of 2-cells; the 2-cells of $\widetilde{\mathbf{Tricat}}(2)$ have an extra choice of the structure constraints $\Pi$ and $M$, and these are forgotten by the functor $U$.  

\begin{thm}
The functor $U:\widetilde{\mathbf{Tricat}}(2) \lra \mathbf{BrMonCat}$ is locally full but not locally faithful.
\end{thm}
\begin{proof}
The first statement follows from Theorem~\ref{idcomp}.  For the second, let $Z$ denote the symmetric monoidal category with only one object $x$ and $Z(x,x) = \mathbb{Z}/2$, with the composition and monoidal structure given on morphisms by addition and all coherence isomorphisms the identity.  The identity functor $1:Z \rightarrow Z$ is a strict braided monoidal functor.  There are two natural transformations $1 \Rightarrow 1$ corresponding to the two different group elements, but only the identity is monoidal.  To show that $U$ is not locally faithful, we will prove that there is more than a single 2-cell $F(1) \Rightarrow F(1)$ in $\widetilde{\mathbf{Tricat}}(2)$; however, there is only one 2-cell $UF(1)=1 \Rightarrow 1 = UF(1)$.

A 2-cell $F(1) \Rightarrow F(1)$ in $\widetilde{\mathbf{Tricat}}(2)$ consists of a monoidal transformation and two group elements satisfying four axioms.  We show that these axioms allow for two different 2-cells.  The first axiom reduces to the equation $\Pi = \Pi$ since all the other cells involved are identities.  The second and third axioms both reduce to $\Pi + M = 0$.  The fourth axiom is then the equation that $\Pi + \Pi = \Pi + \Pi$.  Thus there are two different 2-cells $F(1) \Rightarrow F(1)$ in $\widetilde{\mathbf{Tricat}}(2)$ corresponding to the two different choices of $\Pi$.
\end{proof}

\begin{rem}
The same proof shows that the strict functor $F$ is not locally full.
\end{rem}

Finally note that we could restrict $\widetilde{\mathbf{Tricat}}(2)$ further just for the purposes of getting a biequivalence, as described in the following theorem.  We write $\mathbf{I}$ for the dual pair $(I, I, l, l^{-1})$; this is valid in any monoidal category.

\begin{thm}
There is a 2-category $\mathbf{Tricat}(2)_{2}'$ with
\begin{itemize}
\item 0-cells those doubly degenerate tricategories with $\boxtimes = \otimes$, monoid $I$ with isomorphism $I \cong I$ the identity, all dual pairs $\mathbf{I}$, and all isomorphisms given by unique coherence isomorphisms;
\item 1-cells those functors with all dual pairs $\mathbf{I}$ and all isomorphisms given by unique coherence isomorphisms; and
\item 2-cells those lax transformations with distinguished object $I$ and all constraints given by unique coherence isomorphisms.
\end{itemize}
\end{thm}

The functors $F$ and $U$ then restrict to comparison functors to and from this 2-category, and this does produce a biequivalence, albeit a somewhat tautological one.

\section{Degenerate tricategories}\label{1degen}

We now study degenerate tricategories, and use them to make a definition of monoidal bicategory.  This definition differs from existing definitions \cite{gps1, mcc2} only in that it is fully algebraic.  The difference between these structures becomes more significant at the level of transformation.  In Section~\ref{tricatofmonbicats} we will explore these differences in the process of defining a tricategory of monoidal bicategories.

Since we will {\it define} monoidal bicategories to be degenerate tricategories, a process of ``comparison'' might seem rather circular.  In effect we do little more than observe that our definition of transformation is significantly different from that inherited from \cat{Tricat}, just as in the case of transformations between degenerate bicategories \cite{cg1}.  

\subsection{Basic results}
The results in this section are all obtained by simply rewriting the appropriate definitions using the results of \cite{cg1}.  First we characterise degenerate tricategories and functors between them; this is straightforward.

\begin{thm}
A degenerate tricategory $B$ is precisely
\begin{itemize}
\item a single hom-bicategory which we will also call $B$;
\item a functor $\otimes: B \times B \rightarrow B$;
\item a functor $I: 1 \rightarrow B$;
\item adjoint equivalence $\bs{a}, \bs{l}$, and $\bs{r}$ as in the definition of a tricategory; and
\item invertible modifications $\pi, \mu, \lambda$, and $\rho$ as in the definition of a tricategory
\end{itemize}
all subject to the tricategory axioms.
\end{thm}

\begin{thm}
A weak functor $F:B \rightarrow B'$ between degenerate tricategories is precisely
\begin{itemize}
\item a weak functor $F:B \rightarrow B'$;
\item adjoint equivalences $\bs{\chi}$ and $\bs{\iota}$ as in the definition of weak functor between tricategories; and
\item invertible modifications $\omega, \delta$, and $\gamma$ as in the definition of weak functor, as shown below
\end{itemize}
all subject to axioms which are identical to the functor axioms aside from source and target considerations.
\end{thm}

We now characterise weak transformations, modifications and perturbations.  Here we actually include all the diagrams in the definition, because in Section~\ref{tricatofmonbicats} we will modify these definitions in order to construct a tricategory of monoidal bicategories.

\begin{thm}\label{weaktransform}
A weak transformation $\alpha: F \rightarrow G$ between weak functors of degenerate tricategories is precisely
\begin{itemize}
\item an object $\alpha$ in the target bicategory $B'$, corresponding to the component $\alpha_{\star}$ of the transformation;
\item an adjoint equivalence as displayed below;
\[
\xy 
{\ar_{G} (0,0)*+{B}; (0,-14)*+{B'} };
{\ar^{F} (0,0)*+{B}; (30,0)*+{B'} };
{\ar^{\alpha \otimes' -} (30,0)*+{B'}; (30,-14)*+{B'} };
{\ar_{- \otimes' \alpha} (0,-14)*+{B'}; (30,-14)*+{B'} };
{\ar@{=>}_{\bs{\alpha}} (21, -4)*{}; (9,-10)*{} }
\endxy
\]
\item and invertible modifications as displayed below, where we write $[a,b]$ for $\mbox{{\upshape Hom}}(a,b)$ and $[b,c;a,b]$ for $\mbox{{\upshape Hom}}(b,c) \times \mbox{{\upshape Hom}}(a,b)$;

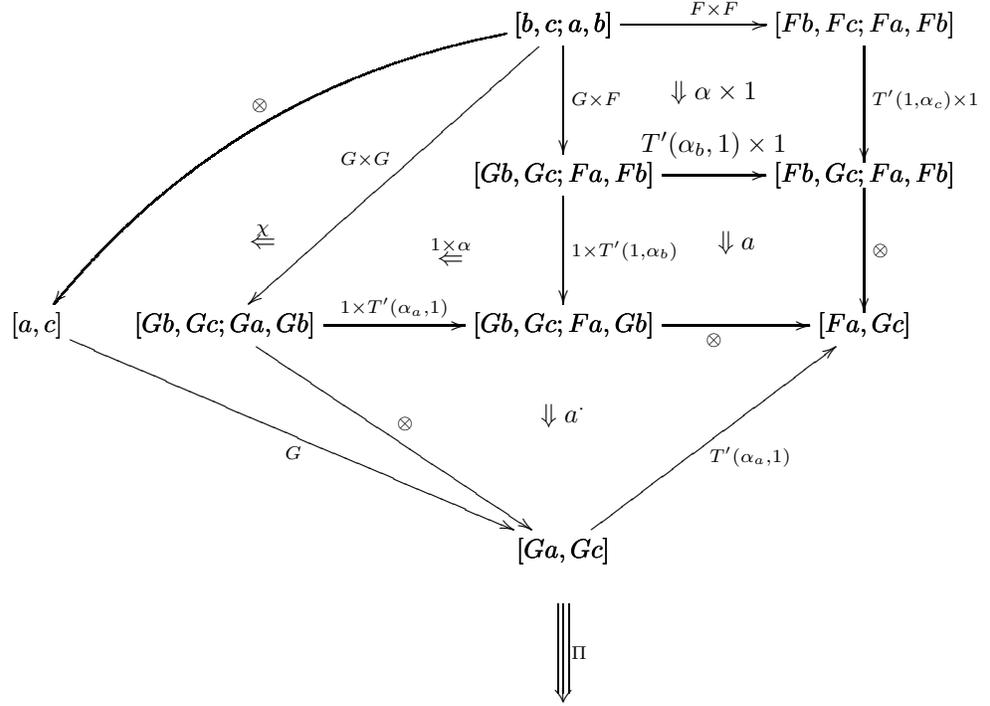
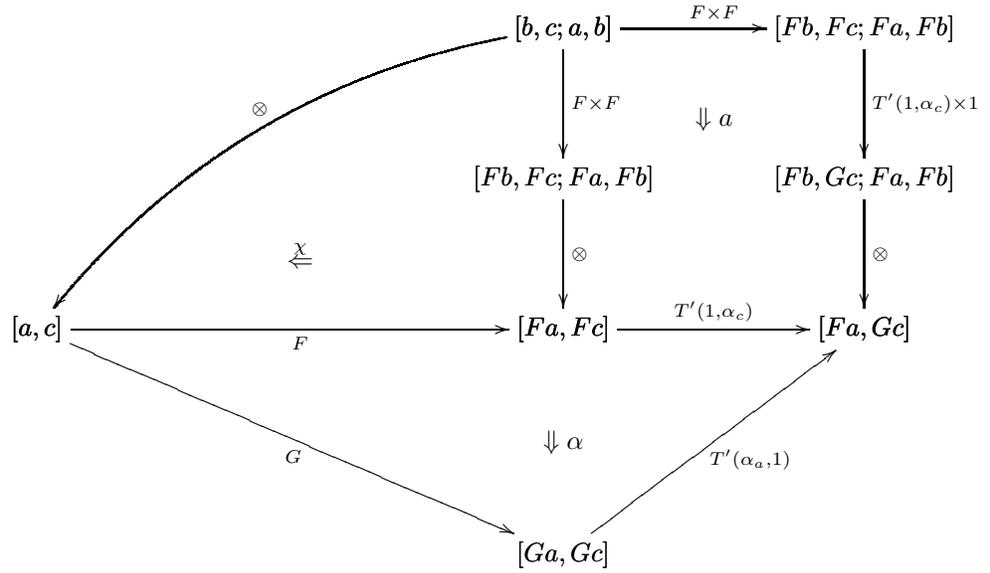
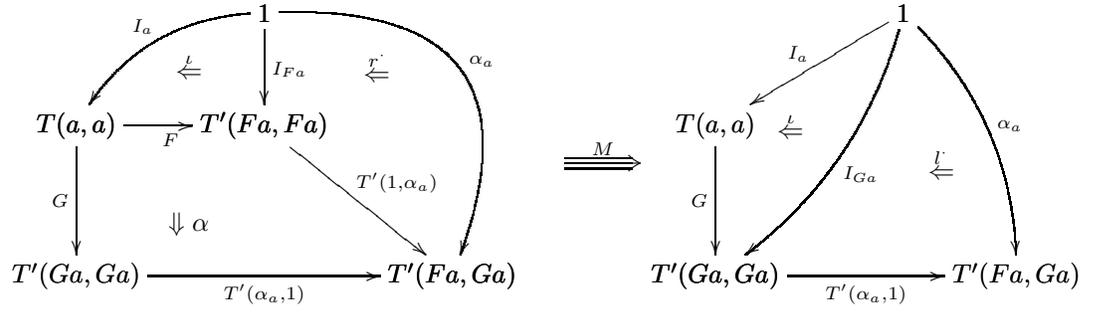
\[
\xy
{\ar@/_2pc/_{\otimes} (0,0)*+{[b,c;a,b]}; (-70,-40)*+{[a,c]} };
{\ar_{G} (-70,-40)*+{[a,c]}; (0,-70)*+{[Ga,Gc]} };
{\ar_{T'(\alpha_{a},1)} (0,-70)*+{[Ga,Gc]}; (40,-40)*+{[Fa,Gc]} };
{\ar^{F \times F} (0,0)*+{[b,c;a,b]}; (40,0)*+{[Fb,Fc; Fa,Fb]} };
{\ar^{T'(1,\alpha_{c}) \times 1} (40,0)*+{[Fb,Fc;Fa,Fb]}; (40,-20)*{[Fb,Gc;Fa,Fb]} };
+{\ar^{\otimes} (40,-20)*{[Fb,Gc;Fa,Fb]}; (40,-40)*{[Fa,Gc]} };
{\ar^{G \times F} (0,0)*+{[b,c;a,b]}; (0,-20)*+{[Gb,Gc; Fa,Fb]} };
{\ar^{1 \times T'(1, \alpha_{b})} (0,-20)*+{[Gb,Gc;Fa,Fb]}; (0,-40)*+{[Gb,Gc;Fa,Gb]} };
{\ar_{\otimes} (0,-40)*+{[Gb,Gc;Fa,Gb]}; (40,-40)*+{[Fa,Gc]} };
{\ar (0,-20)*+{[Gb,Gc;Fa,Fb]}; (40,-20)*+{[Fb,Gc;Fa,Fb]} };
(20,-16)*{T'(\alpha_{b},1) \times 1};
{\ar_{G \times G} (0,0)*+{[b,c;a,b]}; (-45,-40)*+{[Gb,Gc;Ga,Gb]} };
{\ar^{1 \times T'(\alpha_{a}, 1)} (-45,-40)*+{[Gb,Gc;Ga,Gb]}; (0,-40)*+{[Gb,Gc;Fa,Gb]} };
{\ar^{\otimes} (-45,-40)*+{[Gb,Gc;Ga,Gb]}; (0,-70)*+{[Ga,Gc]} };
{\ar@3{->}^{\Pi} (0,-77)*{}; (0,-90)*{} };
{\ar@/_2pc/_{\otimes} (0,-100)*+{[b,c;a,b]}; (-70,-140)*+{[a,c]} };
{\ar_{G} (-70,-140)*+{[a,c]}; (0,-170)*+{[Ga,Gc]} };
{\ar_{T'(\alpha_{a},1)} (0,-170)*+{[Ga,Gc]}; (40,-140)*+{[Fa,Gc]} };
{\ar^{F \times F} (0,-100)*+{[b,c;a,b]}; (40,-100)*+{[Fb,Fc;Fa,Fb]} };
{\ar^{T'(1,\alpha_{c}) \times 1} (40,-100)*+{[Fb,Fc;Fa,Fb]}; (40,-120)*+{[Fb,Gc;Fa,Fb]} };
{\ar^{\otimes} (40,-120)*+{[Fb,Gc;Fa,Fb]}; (40,-140)*+{[Fa,Gc]} };
{\ar^{F \times F} (0,-100)*+{[b,c;a,b]}; (0,-120)*+{[Fb,Fc;Fa,Fb]} };
{\ar^{\otimes} (0,-120)*+{[Fb,Fc;Fa,Fb]}; (0,-140)*+{[Fa,Fc]} };
{\ar^{T'(1, \alpha_{c})} (0,-140)*+{[Fa,Fc]}; (40,-140)*+{[Fa,Gc]} };
{\ar_{F} (-70,-140)*+{[a,c]}; (0,-140)*+{[Fa,Fc]} };
(20,-9)*{\Downarrow \alpha \times 1}; (23,-29)*{\Downarrow a}; (-15,-30)*{\stackrel{1 \times \alpha}{\Leftarrow}}; (0,-52)*{\Downarrow a^{\cdot}}; (-40,-28)*{\stackrel{\chi}{\Leftarrow}};
(20,-112)*{\Downarrow a}; (-35,-130)*{\stackrel{\chi}{\Leftarrow}}; (0,-155)*{\Downarrow \alpha}
\endxy
\]
\[
\xy
{\ar@/_1pc/_{I_{a}} (0,0)*+{1}; (-25,-15)*+{T(a,a)} };
{\ar_{G} (-25,-15)*+{T(a,a)}; (-25,-35)*+{T'(Ga,Ga)} };
{\ar_{T'(\alpha_{a}, 1)} (-25,-35)*+{T'(Ga,Ga)}; (25,-35)*+{T'(Fa,Ga)} };
{\ar@/^4pc/^{\alpha_{a}} (0,0)*+{1}; (25,-35)*+{T'(Fa,Ga)} };
{\ar^{I_{Fa}} (0,0)*+{1}; (0,-15)*+{T'(Fa,Fa)} };
{\ar^{T'(1, \alpha_{a})} (0,-15)*+{T'(Fa,Fa)}; (25,-35)*+{T'(Fa,Ga)} };
{\ar_{F} (-25,-15)*+{T(a,a)}; (0,-15)*+{T'(Fa,Fa)} };
{\ar@3{->}^{M} (40,-20)*{}; (50,-20)*{} };
{\ar_{I_{a}} (85,0)*+{1}; (60,-15)*+{T(a,a)} };
{\ar_{G} (60,-15)*+{T(a,a)}; (60,-35)*+{T'(Ga,Ga)} };
{\ar_{T'(\alpha_{a}, 1)} (60,-35)*+{T'(Ga,Ga)}; (100,-35)*+{T'(Fa,Ga)} };
{\ar@/^1pc/^{\alpha_{a}} (85,0)*+{1}; (100,-35)*+{T'(Fa,Ga)} };
{\ar@/^1pc/^{I_{Ga}} (85,0)*+{1}; (60,-35)*+{T'(Ga,Ga)} };
(-10,-7)*{\stackrel{\iota}{\Leftarrow}}; (-10,-28)*{\Downarrow \alpha}; (15,-7)*{\stackrel{r^{\cdot}}{\Leftarrow}};
(70,-15)*{\stackrel{\iota}{\Leftarrow}}; (90,-20)*{\stackrel{l^{\cdot}}{\Leftarrow}}
\endxy
\]
\end{itemize}
all subject to the following axioms.
\[
\xy
{\ar^{(\alpha1)1} (0,0)*+{((\alpha Ff)Fg)Fh}; (0,20)*+{((Gf \alpha)Fg)Fh} };
{\ar^{a1} (0,20)*+{((Gf \alpha)Fg)Fh}; (5,40)*+{(Gf(\alpha Fg))Fh} };
{\ar^{(1 \alpha) 1} (5,40)*+{(Gf(\alpha Fg))Fh}; (10,60)*+{(Gf(Gg \alpha)) Fh} };
{\ar^{a} (10,60)*+{(Gf(Gg \alpha)) Fh}; (20,80)*+{Gf((Gg \alpha) Fh)} };
{\ar^{1 a} (20,80)*+{Gf((Gg \alpha) Fh)}; (40,100)*+{Gf(Gg(\alpha Fh))} };
{\ar^{1(1 \alpha)} (40,100)*+{Gf(Gg(\alpha Fh))}; (75,100)*+{Gf(Gg(Gh \alpha))} };
{\ar^{a^{\cdot}}  (75,100)*+{Gf(Gg(Gh \alpha))}; (90,80)*+{(GfGg)(Gh \alpha)} };
{\ar^{\chi 1} (90,80)*+{(GfGg)(Gh \alpha)}; (110,60)*+{G(fg)(Gh \alpha)} };
{\ar^{a^{\cdot}} (110,60)*+{G(fg)(Gh \alpha)}; (110,40)*+{(G(fg)Gh)\alpha} };
{\ar^{\chi 1} (110,40)*+{(G(fg)Gh)\alpha}; (115, 20)*+{G((fg)h) \alpha} };
{\ar^{Ga 1} (115, 20)*+{G((fg)h) \alpha}; (115,0)*+{G(f(gh)) \alpha} };
{\ar_{a} (0,0)*+{((\alpha Ff)Fg)Fh}; (15,-20)*+{(\alpha Ff)(Fg Fh)} };
{\ar_{1 \chi} (15,-20)*+{(\alpha Ff)(Fg Fh)}; (40,-50)*+{(\alpha Ff)F(gh)} };
{\ar_{a} (40,-50)*+{(\alpha Ff)F(gh)}; (75,-50)*+{\alpha(FfF(gh))} };
{\ar_{1 \chi} (75,-50)*+{\alpha(FfF(gh))}; (95,-25)*+{\alpha F(f(gh))} };
{\ar_{\alpha} (95,-25)*+{\alpha F(f(gh))}; (115,0)*+{G(f(gh)) \alpha} };
{\ar^{a1} (0,0)*+{((\alpha Ff)Fg)Fh}; (25,16)*+{(\alpha(FfFg))Fh} };
{\ar^{(1 \chi)1} (25,16)*+{(\alpha(FfFg))Fh}; (33,32)*+{(\alpha F(fg))Fh} };
{\ar^{\alpha 1} (33,32)*+{(\alpha F(fg))Fh}; (41, 48)*+{(G(fg)\alpha)Fh} };
{\ar^{a^{\cdot}1} (10,60)*+{(Gf(Gg \alpha)) Fh}; (35,70)*+{((GfGg)\alpha)Fh} };
{\ar_{(\chi 1)1} (35,70)*+{((GfGg)\alpha)Fh}; (41, 48)*+{(G(fg)\alpha)Fh} };
{\ar_{a} (35,70)*+{((GfGg)\alpha)Fh}; (55,85)*+{(GfGg)(\alpha Fh)} };
{\ar^{a^{\cdot}} (40,100)*+{Gf(Gg(\alpha Fh))}; (55,85)*+{(GfGg)(\alpha Fh)} };
{\ar^{1 \alpha} (55,85)*+{(GfGg)(\alpha Fh)}; (90,80)*+{(GfGg)(Gh \alpha)} };
{\ar_{a} (41, 48)*+{(G(fg)\alpha)Fh}; (75,63)*+{G(fg)(\alpha Fh)} };
{\ar^{\chi 1} (55,85)*+{(GfGg)(\alpha Fh)}; (75,63)*+{G(fg)(\alpha Fh)} };
{\ar_{1 \alpha} (75,63)*+{G(fg)(\alpha Fh)}; (110,60)*+{G(fg)(Gh \alpha)} };
{\ar^{a} (33,32)*+{(\alpha F(fg))Fh}; (61,16)*+{\alpha(F(fg)Fh)} };
{\ar^{1 \chi} (61,16)*+{\alpha(F(fg)Fh)}; (88,0)*+{\alpha F((fg)h)} };
{\ar^{\alpha} (88,0)*+{\alpha F((fg)h)}; (115, 20)*+{G((fg)h) \alpha} };
{\ar^{1Fa} (88,0)*+{\alpha F((fg)h)}; (95,-25)*+{\alpha F(f(gh))} };
{\ar_{a} (25,16)*+{(\alpha(FfFg))Fh}; (48,-5)*+{\alpha((FfFg)Fh)} };
{\ar_{1(\chi1)} (48,-5)*+{\alpha((FfFg)Fh)}; (61,16)*+{\alpha(F(fg)Fh)} };
{\ar^{1a} (48,-5)*+{\alpha((FfFg)Fh)}; (60,-25)*+{\alpha(Ff(FgFh))} };
{\ar^{1(1\chi)} (60,-25)*+{\alpha(Ff(FgFh))}; (75,-50)*+{\alpha(FfF(gh))} };
{\ar_{a} (15,-20)*+{(\alpha Ff)(Fg Fh)}; (60,-25)*+{\alpha(Ff(FgFh))} };
(20,50)*{\Downarrow \Pi 1}; (35,85)*{\Downarrow \pi}; (65,93)*{\cong}; (55,66)*{\cong}; (82,70)*{\cong}; (75,38)*{\Downarrow \Pi}; (45,10)*{\cong}; (25,-5)*{\Downarrow \pi}; (48,-37)*{\cong}; (75,-12)*{\Downarrow 1 \omega}; (100,-5)*{\cong};
{\ar@{=} (57,-55)*{}; (57,-65)*{} }
\endxy
\]
\[
\xy
{\ar^{(\alpha1)1} (0,0)*+{((\alpha Ff)Fg)Fh}; (0,20)*+{((Gf \alpha)Fg)Fh} };
{\ar^{a1} (0,20)*+{((Gf \alpha)Fg)Fh}; (5,40)*+{(Gf(\alpha Fg))Fh} };
{\ar^{(1 \alpha) 1} (5,40)*+{(Gf(\alpha Fg))Fh}; (10,60)*+{(Gf(Gg \alpha)) Fh} };
{\ar^{a} (10,60)*+{(Gf(Gg \alpha)) Fh}; (20,80)*+{Gf((Gg \alpha) Fh)} };
{\ar^{1 a} (20,80)*+{Gf((Gg \alpha) Fh)}; (40,100)*+{Gf(Gg(\alpha Fh))} };
{\ar^{1(1 \alpha)} (40,100)*+{Gf(Gg(\alpha Fh))}; (75,100)*+{Gf(Gg(Gh \alpha))} };
{\ar^{a^{\cdot}}  (75,100)*+{Gf(Gg(Gh \alpha))}; (90,80)*+{(GfGg)(Gh \alpha)} };
{\ar^{\chi 1} (90,80)*+{(GfGg)(Gh \alpha)}; (110,60)*+{G(fg)(Gh \alpha)} };
{\ar^{a^{\cdot}} (110,60)*+{G(fg)(Gh \alpha)}; (110,40)*+{(G(fg)Gh)\alpha} };
{\ar^{\chi 1} (110,40)*+{(G(fg)Gh)\alpha}; (115, 20)*+{G((fg)h) \alpha} };
{\ar^{Ga 1} (115, 20)*+{G((fg)h) \alpha}; (115,0)*+{G(f(gh)) \alpha} };
{\ar_{a} (0,0)*+{((\alpha Ff)Fg)Fh}; (15,-20)*+{(\alpha Ff)(Fg Fh)} };
{\ar_{1 \chi} (15,-20)*+{(\alpha Ff)(Fg Fh)}; (40,-50)*+{(\alpha Ff)F(gh)} };
{\ar_{a} (40,-50)*+{(\alpha Ff)F(gh)}; (75,-50)*+{\alpha(FfF(gh))} };
{\ar_{1 \chi} (75,-50)*+{\alpha(FfF(gh))}; (95,-25)*+{\alpha F(f(gh))} };
{\ar_{\alpha} (95,-25)*+{\alpha F(f(gh))}; (115,0)*+{G(f(gh)) \alpha} };
{\ar_{\alpha 1} (40,-50)*+{(\alpha Ff)F(gh)}; (45,-25)*+{(Gf \alpha) F(gh)} };
{\ar_{\alpha 1} (15,-20)*+{(\alpha Ff)(Fg Fh)}; (23,7)*+{(Gf \alpha)(Fg Fh)} };
{\ar^{1 \chi} (23,7)*+{(Gf \alpha)(Fg Fh)}; (45,-25)*+{(Gf \alpha) F(gh)} };
{\ar^{a} (0,20)*+{((Gf \alpha)Fg)Fh}; (23,7)*+{(Gf \alpha)(Fg Fh)} };
{\ar_{a} (5,40)*+{(Gf(\alpha Fg))Fh}; (35,40)*+{Gf((\alpha Fg)Fh)} };
{\ar^{1 a} (35,40)*+{Gf((\alpha Fg)Fh)}; (40,20)*+{Gf(\alpha (FgFh))} };
{\ar^{a} (23,7)*+{(Gf \alpha)(Fg Fh)}; (40,20)*+{Gf(\alpha (FgFh))} };
{\ar^{1(1\chi)} (40,20)*+{Gf(\alpha (FgFh))}; (60,0)*+{Gf(\alpha F(gh))} };
{\ar_{a} (45,-25)*+{(Gf \alpha) F(gh)}; (60,0)*+{Gf(\alpha F(gh))} };
{\ar_{1 \alpha} (60,0)*+{Gf(\alpha F(gh))}; (70,23)*+{Gf(G(gh) \alpha)} };
{\ar_{a^{\cdot}} (70,23)*+{Gf(G(gh) \alpha)}; (90,10)*+{(GfG(gh))\alpha} };
{\ar_{\chi 1} (90,10)*+{(GfG(gh))\alpha}; (115,0)*+{G(f(gh)) \alpha} };
{\ar_{1a^{\cdot}} (75,100)*+{Gf(Gg(Gh \alpha))}; (55,59)*+{Gf((GgGh)\alpha)} };
{\ar_{1(\chi 1)} (55,59)*+{Gf((GgGh)\alpha)}; (70,23)*+{Gf(G(gh) \alpha)} };
{\ar_{a^{\cdot}} (90,80)*+{(GfGg)(Gh \alpha)}; (85,52)*+{((GfGg)Gh)\alpha} };
{\ar_{a1} (85,52)*+{((GfGg)Gh)\alpha}; (90,33)*+{(Gf(GgGh))\alpha} };
{\ar^{(1 \chi)1} (90,33)*+{(Gf(GgGh))\alpha}; (90,10)*+{(GfG(gh))\alpha} };
{\ar^{a^{\cdot}} (55,59)*+{Gf((GgGh)\alpha)}; (90,33)*+{(Gf(GgGh))\alpha} };
{\ar^{(\chi 1)1} (85,52)*+{((GfGg)Gh)\alpha}; (110,40)*+{(G(fg)Gh)\alpha} };
{\ar_{1(\alpha 1)} (35,40)*+{Gf((\alpha Fg)Fh)}; (20,80)*+{Gf((Gg \alpha) Fh)} };
(20,50)*{\cong}; (20,30)*{\Downarrow \pi}; (7,10)*{\cong}; (45,80)*{\Downarrow 1 \Pi}; (40,-1)*{\cong}; (27,-12)*{\cong}; (75,70)*{\stackrel{\pi}{\Leftarrow}}; (71,38)*{\cong}; (95,65)*{\cong}; (105,27)*{\Downarrow \omega 1}; (75,-18)*{\Downarrow \Pi}
\endxy
\]
\[
\xy
{\ar^{r^{\cdot}1} (0,0)*+{\alpha Ff}; (15,15)*+{(\alpha I)Ff} };
{\ar^{(1 \iota)1} (15,15)*+{(\alpha I)Ff}; (30,30)*+{(\alpha FI)Ff} };
{\ar^{\alpha 1} (30,30)*+{(\alpha FI)Ff}; (45,42)*+{(GI \alpha)Ff} };
{\ar^{a} (45,42)*+{(GI \alpha)Ff}; (60,53)*+{GI(\alpha Ff)} };
{\ar^{1 \alpha} (60,53)*+{GI(\alpha Ff)}; (75,42)*+{GI(Gf \alpha)} };
{\ar^{a^{\cdot}} (75,42)*+{GI(Gf \alpha)}; (90,30)*+{(GIGf)\alpha} };
{\ar^{\chi 1} (90,30)*+{(GIGf)\alpha}; (105,15)*+{G(If) \alpha} };
{\ar^{Gl 1}  (105,15)*+{G(If) \alpha}; (120,0)*+{Gf \alpha} };
{\ar_{l^{\cdot} 1} (0,0)*+{\alpha Ff}; (30,-20)*+{(I \alpha) Ff} };
{\ar_{a} (30,-20)*+{(I \alpha) Ff}; (60,-30)*+{I(\alpha Ff)} };
{\ar_{1 \alpha} (60,-30)*+{I(\alpha Ff)}; (90,-20)*+{I(Gf \alpha)} };
{\ar_{l} (90,-20)*+{I(Gf \alpha)}; (120,0)*+{Gf \alpha} };
{\ar^{a} (30,30)*+{(\alpha FI)Ff}; (55,25)*+{\alpha(FIFf)} };
{\ar^{1 \chi} (55,25)*+{\alpha(FIFf)}; (80,20)*+{\alpha F(If)} };
{\ar^{\alpha} (80,20)*+{\alpha F(If)}; (105,15)*+{G(If) \alpha} };
{\ar^{a} (15,15)*+{(\alpha I)Ff}; (47.5,7.5)*+{\alpha(IFf)} };
{\ar^{1l} (47.5,7.5)*+{\alpha(IFf)}; (80,0)*+{\alpha Ff} };
{\ar_{\alpha} (80,0)*+{\alpha Ff}; (120,0)*+{Gf \alpha} };
{\ar^{1 (\iota 1)} (47.5,7.5)*+{\alpha(IFf)}; (55,25)*+{\alpha(FIFf)} };
{\ar^{1 Fl} (80,20)*+{\alpha F(If)}; (80,0)*+{\alpha Ff} };
{\ar_{1} (0,0)*+{\alpha Ff}; (80,0)*+{\alpha Ff} };
{\ar^{l} (60,-30)*+{I(\alpha Ff)}; (80,0)*+{\alpha Ff} };
(60,37)*{\Downarrow \Pi}; (35,22)*{\cong}; (25,8)*{\Downarrow \mu}; (65,17)*{\Downarrow 1 \gamma}; (92,8)*{\cong}; (90,-10)*{\cong}; (40,-10)*{\Downarrow \lambda};
{\ar@{=} (60,-35)*{}; (60,-45)*+{} };
{\ar^{r^{\cdot}1} (0,-103)*+{\alpha Ff}; (15,-88)*+{(\alpha I)Ff} };
{\ar^{(1 \iota)1} (15,-88)*+{(\alpha I)Ff}; (30,-73)*+{(\alpha FI)Ff} };
{\ar^{\alpha 1} (30,-73)*+{(\alpha FI)Ff}; (45,-61)*+{(GI \alpha)Ff} };
{\ar^{a} (45,-61)*+{(GI \alpha)Ff}; (60,-50)*+{GI(\alpha Ff)} };
{\ar^{1 \alpha} (60,-50)*+{GI(\alpha Ff)}; (75,-61)*+{GI(Gf \alpha)} };
{\ar^{a^{\cdot}} (75,-61)*+{GI(Gf \alpha)}; (90,-73)*+{(GIGf)\alpha} };
{\ar^{\chi 1} (90,-73)*+{(GIGf)\alpha}; (105,-88)*+{G(If) \alpha} };
{\ar^{Gl 1}  (105,-88)*+{G(If) \alpha}; (120,-103)*+{Gf \alpha} };
{\ar_{l^{\cdot} 1} (0,-103)*+{\alpha Ff}; (30,-123)*+{(I \alpha) Ff} };
{\ar_{a} (30,-123)*+{(I \alpha) Ff}; (60,-133)*+{I(\alpha Ff)} };
{\ar_{1 \alpha} (60,-133)*+{I(\alpha Ff)}; (90,-123)*+{I(Gf \alpha)} };
{\ar_{l} (90,-123)*+{I(Gf \alpha)}; (120,-103)*+{Gf \alpha} };
{\ar_{(\iota 1)1} (30,-123)*+{(I \alpha) Ff}; (45,-61)*+{(GI \alpha)Ff} };
{\ar_{\iota1} (60,-133)*+{I(\alpha Ff)}; (60,-50)*+{GI(\alpha Ff)} };
{\ar@/^1pc/^{\iota 1} (90,-123)*+{I(Gf \alpha)}; (75,-61)*+{GI(Gf \alpha)} };
{\ar_{a^{\cdot}} (90,-123)*+{I(Gf \alpha)}; (90,-98)*+{(IGf)\alpha} };
{\ar^{(\iota 1)1} (90,-98)*+{(IGf)\alpha}; (90,-73)*+{(GIGf)\alpha} };
{\ar^{l1} (90,-98)*+{(IGf)\alpha}; (120,-103)*+{Gf \alpha} };
(20,-98)*{\Downarrow M 1}; (50,-98)*{\cong}; (70,-98)*{\cong}; (100,-92)*{\Downarrow \gamma 1}; (100,-105)*{\Downarrow \lambda}; (83,-80)*{\cong}
\endxy
\]
\vspace*{-3em}\[
\xy
{\ar^{\alpha} (0,0)*+{\alpha Ff}; (11,20)*+{Gf \alpha} };
{\ar^{1 r^{\cdot}} (11,20)*+{Gf \alpha}; (30,40)*+{Gf( \alpha I)} };
{\ar^{1(1 \iota)} (30,40)*+{Gf( \alpha I)}; (60,60)*+{Gf(\alpha FI)} };
{\ar^{1 \alpha} (60,60)*+{Gf(\alpha FI)}; (80,40)*+{Gf(GI \alpha)} };
{\ar^{a^{\cdot}} (80,40)*+{Gf(GI \alpha)}; (100,20)*+{(GfGI) \alpha} };
{\ar^{\chi 1} (100,20)*+{(GfGI) \alpha}; (120,0)*+{G(fI) \alpha} };
{\ar@/_2pc/_{1}  (0,0)*+{\alpha Ff}; (70,-30)*+{\alpha Ff} };
{\ar_{\alpha}  (70,-30)*+{\alpha Ff};  (95,-20)*+{Gf \alpha} };
{\ar_{Gr^{\cdot} 1} (95,-20)*+{Gf \alpha}; (120,0)*+{G(fI) \alpha} };
{\ar_{r^{\cdot}} (0,0)*+{\alpha Ff}; (38,5)*+{(\alpha Ff)I} };
{\ar^{\alpha 1} (38,5)*+{(\alpha Ff)I}; (30,23)*+{(Gf \alpha)I} };
{\ar_{r^{\cdot}} (11,20)*+{Gf \alpha}; (30,23)*+{(Gf \alpha)I} };
{\ar_{a} (30,23)*+{(Gf \alpha)I}; (30,40)*+{Gf( \alpha I)} };
{\ar_{1 \iota} (38,5)*+{(\alpha Ff)I}; (65,18)*+{(\alpha Ff)FI} };
{\ar_{\alpha 1} (65,18)*+{(\alpha Ff)FI}; (55,35)*+{(Gf \alpha)FI} };
{\ar_{a} (55,35)*+{(Gf \alpha)FI}; (60,60)*+{Gf(\alpha FI)} };
{\ar^{1 \iota} (30,23)*+{(Gf \alpha)I}; (55,35)*+{(Gf \alpha)FI} };
{\ar_{a} (38,5)*+{(\alpha Ff)I}; (55,-13)*+{\alpha(Ff I)} };
{\ar^{1(1 \iota)} (55,-13)*+{\alpha(Ff I)}; (75,5)*+{\alpha (Ff FI)} };
{\ar^{a} (65,18)*+{(\alpha Ff)FI}; (75,5)*+{\alpha (Ff FI)} };
{\ar_{1 r} (55,-13)*+{\alpha(Ff I)}; (70,-30)*+{\alpha Ff} };
{\ar^{1 Fr^{\cdot}} (70,-30)*+{\alpha Ff}; (95,-5)*+{\alpha F(fI)} };
{\ar^{1 \chi} (75,5)*+{\alpha (Ff FI)}; (95,-5)*+{\alpha F(fI)} };
{\ar^{\alpha} (95,-5)*+{\alpha F(fI)}; (120,0)*+{G(fI) \alpha} };
(17,10)*{\cong}; (24,28)*{\Downarrow \rho}; (45,41)*{\cong}; (47,18)*{\cong}; (58,4)*{\cong}; (70,-10)*{\Downarrow 1 \delta}; (100,-10)*{\cong}; (35,-8)*{\Downarrow \rho}; (82,22)*{\Downarrow \Pi};
{\ar@{=} (60,-35)*{}; (60,-45)*{} };
{\ar^{\alpha} (0,-110)*+{\alpha Ff}; (11,-90)*+{Gf \alpha} };
{\ar^{1 r^{\cdot}} (11,-90)*+{Gf \alpha}; (30,-70)*+{Gf( \alpha I)} };
{\ar^{1(1 \iota)} (30,-70)*+{Gf( \alpha I)}; (60,-50)*+{Gf(\alpha FI)} };
{\ar^{1 \alpha} (60,-50)*+{Gf(\alpha FI)}; (80,-70)*+{Gf(GI \alpha)} };
{\ar^{a^{\cdot}} (80,-70)*+{Gf(GI \alpha)}; (100,-90)*+{(GfGI) \alpha} };
{\ar^{\chi 1} (100,-90)*+{(GfGI) \alpha}; (120,-110)*+{G(fI) \alpha} };
{\ar@/_2pc/_{1}  (0,-110)*+{\alpha Ff}; (70,-140)*+{\alpha Ff} };
{\ar_{\alpha}  (70,-140)*+{\alpha Ff};  (95,-130)*+{Gf \alpha} };
{\ar_{Gr^{\cdot} 1} (95,-130)*+{Gf \alpha}; (120,-110)*+{G(fI) \alpha} };
{\ar^{1 l^{\cdot}} (11,-90)*+{Gf \alpha}; (55,-82)*+{Gf (I \alpha)} };
{\ar_{1(\iota 1)} (55,-82)*+{Gf (I \alpha)}; (80,-70)*+{Gf(GI \alpha)} };
{\ar_{1} (11,-90)*+{Gf \alpha}; (95,-130)*+{Gf \alpha} };
{\ar^{a^{\cdot}} (55,-82)*+{Gf (I \alpha)}; (75,-100)*+{(GfI) \alpha} };
{\ar^{r 1} (75,-100)*+{(GfI) \alpha}; (95,-130)*+{Gf \alpha} };
{\ar^{(1 \iota)1} (75,-100)*+{(GfI) \alpha}; (100,-90)*+{(GfGI) \alpha} };
(45,-70)*{\Downarrow 1 M}; (35,-118)*{\cong}; (50,-92)*{\Downarrow \mu}; (80,-85)*{\cong}; (100,-107)*{\Downarrow \delta 1}
\endxy
\]
\end{thm}

\begin{thm}
Let $\alpha$ and $\beta$ be transformations $F \Rightarrow G$ of degenerate tricategories.  Then a modification $m: \alpha \Rightarrow \beta$ is precisely
\begin{itemize}
\item a 1-cell $m: \alpha \rightarrow \beta$ in $B'$ and
\item invertible modifications as shown below,
\[
\xy
{\ar_{G} (0,0)*+{T(a,b)}; (0,-16)*+{T'(Ga,Gb)} };
{\ar^{F} (0,0)*+{T(a,b)}; (36,0)*+{T'(Fa,Fb)} };
{\ar^{T'(1,\alpha_{b})} (36,0)*+{T'(Fa,Fb)}; (36,-16)*+{T'(Fa, Gb)} };
{\ar^{T'(\alpha_{a},1)} (0,-16)*+{T'(Ga,Gb)}; (36,-16)*+{T'(Fa, Gb)} };
{\ar@{=>}_{\alpha} (25, -3)*{}; (10,-10)*{} };
{\ar@/_3pc/_{T'(\beta_{a}, 1)} (0,-16)*+{T'(Ga,Gb)}; (36,-16)*+{T'(Fa, Gb)} };
{\ar@{=>}^{(m_{a})^{*}} (18,-18)*{}; (18,-23)*{} };
{\ar@3{->}^{m} (50,-12)*{}; (60,-12)*{} };
{\ar^{F} (66,0)*+{T(a,b)}; (102,0)*+{T'(Fa,Fb)} };
{\ar_{G} (66,0)*+{T(a,b)}; (66,-28)*+{T'(Ga,Gb)} };
{\ar_{T'(1,\beta_{b})} (102,0)*+{T'(Fa,Fb)}; (102,-28)*+{T'(Fa, Gb)} };
{\ar_{T'(\beta_{a}, 1)} (66,-28)*+{T'(Ga,Gb)}; (102,-28)*+{T'(Fa, Gb)} };
{\ar@/^3pc/^{T'(1, \alpha_{b})} (102,0)*+{T'(Fa,Fb)}; (102,-28)*+{T'(Fa, Gb)} };
{\ar@{=>}_{(m_{b})_{*}} (112,-14)*{}; (105,-14)*{} };
{\ar@{=>}_{\beta} (90,-4)*{}; (75,-21)*{} }
\endxy \]
\end{itemize}
all subject to the following two axioms (unmarked isomorphisms are naturality isomorphisms).
\[
\xy
{\ar^{\alpha  1} (0,0)*+{(\alpha Ff)Fg}; (25,0)*+{(Gf \alpha) Fg} };
{\ar^{a} (25,0)*+{(Gf \alpha) Fg}; (50,0)*+{Gf(\alpha Fg)} };
{\ar^{1  \alpha} (50,0)*+{Gf(\alpha Fg)}; (75,0)*+{Gf(Gg \alpha)} };
{\ar^{a^{\cdot}} (75,0)*+{Gf(Gg \alpha)}; (100,0)*+{(Gf Gg) \alpha} };
{\ar^{\chi  1} (100,0)*+{(Gf Gg) \alpha}; (125,0)*+{G(fg) \alpha} };
{\ar^{G1   m} (125,0)*+{G(fg) \alpha}; (125,-20)*+{G(fg) \beta} };
{\ar_{(m  F1)  F1} (0,0)*+{(\alpha Ff)Fg}; (0,-20)*+{(\beta Ff)Fg} };
{\ar_{a} (0,-20)*+{(\beta Ff)Fg}; (42,-40)*+{\beta(FfFg)} };
{\ar_{1  \chi} (42,-40)*+{\beta(FfFg)}; (83,-40)*+{\beta F(fg)} };
{\ar_{\beta} (83,-40)*+{\beta F(fg)}; (125,-20)*+{G(fg) \beta} };
{\ar^{\beta  1} (0,-20)*+{(\beta Ff)Fg}; (25,-20)*+{(Gf \beta) Fg} };
{\ar_{a} (25,-20)*+{(Gf \beta) Fg}; (50,-20)*{Gf(\beta Fg)} };
{\ar_{1  \beta} (50,-20)*{Gf(\beta Fg)}; (75,-20)*+{Gf(Gg \beta)} };
{\ar_{a^{\cdot}} (75,-20)*+{Gf(Gg \beta)}; (100,-20)*+{(GfGg) \beta} };
{\ar_{\chi  1} (100,-20)*+{(GfGg) \beta}; (125,-20)*+{G(fg) \beta} };
{\ar|{(G1  m)  F1} (25,0)*+{(Gf \alpha) Fg}; (25,-20)*+{(Gf \beta) Fg} };
{\ar|{G1  (m  F1)} (50,0)*+{Gf(\alpha Fg)}; (50,-20)*{Gf(\beta Fg)} };
{\ar|{G1  (G1  m)} (75,0)*+{Gf(Gg \alpha)}; (75,-20)*+{Gf(Gg \beta)} };
{\ar|{(G1  G1)  m} (100,0)*+{(Gf Gg) \alpha}; (100,-20)*+{(GfGg) \beta} };
{\ar^{\alpha  1} (0,-55)*+{(\alpha Ff)Fg}; (25,-55)*+{(Gf \alpha) Fg} };
{\ar^{a} (25,-55)*+{(Gf \alpha) Fg}; (50,-55)*+{Gf(\alpha Fg)} };
{\ar^{1  \alpha} (50,-55)*+{Gf(\alpha Fg)}; (75,-55)*+{Gf(Gg \alpha)} };
{\ar^{a^{\cdot}} (75,-55)*+{Gf(Gg \alpha)}; (100,-55)*+{(Gf Gg) \alpha} };
{\ar^{\chi  1} (100,-55)*+{(Gf Gg) \alpha}; (125,-55)*+{G(fg) \alpha} };
{\ar^{G1   m} (125,-55)*+{G(fg) \alpha}; (125,-75)*+{G(fg) \beta} };
{\ar_{(m  F1)  F1} (0,-55)*+{(\alpha Ff)Fg}; (0,-75)*+{(\beta Ff)Fg} };
{\ar_{a} (0,-75)*+{(\beta Ff)Fg}; (42,-95)*+{\beta(FfFg)} };
{\ar_{1  \chi} (42,-95)*+{\beta(FfFg)}; (83,-95)*+{\beta F(fg)} };
{\ar_{\beta} (83,-95)*+{\beta F(fg)}; (125,-75)*+{G(fg) \beta} };
{\ar^{a} (0,-55)*+{(\alpha Ff)Fg}; (42,-75)*+{\alpha(Ff Fg)} };
{\ar_{1  \chi} (42,-75)*+{\alpha(Ff Fg)}; (83,-75)*+{\alpha F(fg)} };
{\ar_{\alpha} (83,-75)*+{\alpha F(fg)}; (125,-55)*+{G(fg) \alpha} };
{\ar^{m  (F1  F1)} (42,-75)*+{\alpha(Ff Fg)}; (42,-95)*+{\beta(FfFg)} };
{\ar_{m  F(1  1)} (83,-75)*+{\alpha F(fg)}; (83,-95)*+{\beta F(fg)} };
(9,-9)*{\Downarrow m1}; (37.5,-9)*{\cong}; (62.5,-9)*{\stackrel{1m}{\Downarrow}}; (87.5,-9)*{\cong}; (112.5,-9)*{\cong}; (62.5,-31)*{\Downarrow \Pi};
{\ar@{=} (62.5,-45)*{}; (62.5,-50)*{} };
(62.5,-65)*{\Downarrow \Pi}; (22,-75)*{\cong}; (102,-75)*{\Downarrow m}; (62.5,-85)*{\cong}
\endxy
\]
\[
\xy
{\ar^{r^{\cdot}} (0,0)*+{\alpha}; (35,0)*+{\alpha I} };
{\ar^{1 \iota} (35,0)*+{\alpha I}; (70,0)*+{\alpha FI} };
{\ar^{\alpha} (70,0)*+{\alpha FI}; (105,0)*+{GI \alpha} };
{\ar^{1m} (105,0)*+{GI \alpha}; (105,-20)*+{GI \beta} };
{\ar_{m} (0,0)*+{\alpha}; (0,-20)*+{\beta} };
{\ar_{l^{\cdot}} (0,-20)*+{\beta}; (52,-40)*+{I \beta} };
{\ar_{\iota 1} (52,-40)*+{I \beta}; (105,-20)*+{GI \beta} };
{\ar^{r^{\cdot}} (0,-20)*+{\beta}; (35,-20)*+{\beta I} };
{\ar_{1 \iota} (35,-20)*+{\beta I}; (70,-20)*+{1 \iota} };
{\ar^{\beta} (70,-20)*+{1 \iota}; (105,-20)*+{GI \beta} };
{\ar^{m 1} (35,0)*+{\alpha I}; (35,-20)*+{\beta I} };
{\ar^{m 1} (70,0)*+{\alpha FI}; (70,-20)*+{1 \iota} };
{\ar@{=} (52,-47)*{}; (52,-54)*{} };
{\ar^{r^{\cdot}} (0,-60)*+{\alpha}; (35,-60)*+{\alpha I} };
{\ar^{1 \iota} (35,-60)*+{\alpha I}; (70,-60)*+{\alpha FI} };
{\ar^{\alpha} (70,-60)*+{\alpha FI}; (105,-60)*+{GI \alpha} };
{\ar^{1m} (105,-60)*+{GI \alpha}; (105,-80)*+{GI \beta} };
{\ar_{m} (0,-60)*+{\alpha}; (0,-80)*+{\beta} };
{\ar_{l^{\cdot}} (0,-80)*+{\beta}; (52,-100)*+{I \beta} };
{\ar_{\iota 1} (52,-100)*+{I \beta}; (105,-80)*+{GI \beta} };
{\ar_{l^{\cdot}} (0,-60)*+{\alpha}; (52,-78)*+{I \alpha} };
{\ar_{\iota 1} (52,-78)*+{I \alpha}; (105,-60)*+{GI \alpha} };
{\ar^{1m} (52,-78)*+{I \alpha}; (52,-100)*+{I \beta} };
(17,-10)*{\cong}; (52,-10)*{\cong}; (87,-10)*{\Downarrow m};
(52,-28)*{\Downarrow M};
(52,-68)*{\Downarrow M}; (25,-80)*{\cong}; (80,-80)*{\cong}
\endxy
\]
\end{thm}

\begin{thm}
A perturbation $\sigma: m \Rrightarrow n$ of degenerate tricategories is precisely a 2-cell $\sigma: m \Rightarrow n$ in $B'$ such that the following axiom holds.
\[
\xy
{\ar^{\alpha} (0,0)*+{\alpha \otimes Ff}; (35,0)*+{Gf \otimes \alpha} };
{\ar^{1 \otimes m} (35,0)*+{Gf \otimes \alpha}; (35,-20)*+{Gf \otimes \beta} };
{\ar^{m \otimes 1} (0,0)*+{\alpha \otimes Ff}; (0,-20)*+{\beta \otimes Ff} };
{\ar_{\beta} (0,-20)*+{\beta \otimes Ff}; (35,-20)*+{Gf \otimes \beta} };
{\ar@/_3pc/_{n \otimes 1} (0,0)*+{\alpha \otimes Ff}; (0,-20)*+{\beta \otimes Ff} };
{\ar@{=} (45,-10)*{}; (51,-10)*+{} };
{\ar^{\alpha} (60,0)*+{\alpha \otimes Ff}; (95,0)*+{Gf \otimes \alpha} };
{\ar_{1 \otimes n} (95,0)*+{Gf \otimes \alpha}; (95,-20)*+{Gf \otimes \beta} };
{\ar@/^3pc/^{1 \otimes m} (95,0)*+{Gf \otimes \alpha}; (95,-20)*+{Gf \otimes \beta} };
{\ar_{n \otimes 1} (60,0)*+{\alpha \otimes Ff}; (60,-20)*+{\beta \otimes Ff} };
{\ar_{\beta} (60,-20)*+{\beta \otimes Ff}; (95,-20)*+{Gf \otimes \beta} };
(-5,-10)*{\stackrel{\sigma \otimes 1}{\Leftarrow}}; (100,-10)*+{\stackrel{1 \otimes \sigma}{\Leftarrow}};
{\ar@{=>}_{m} (28,-5)*{}; (7,-15)*{} };
{\ar@{=>}_{n} (88,-5)*{}; (67,-15)*{} }
\endxy
\]
\end{thm}

\subsection{Overall structure}\label{tricatofmonbicats}

In this section, we will construct a tricategory of monoidal bicategories.  The objects and 1-cells will be given by degenerate tricategories and functors between them, but the higher cells will be given only by special transformations and modifications which have their components at the lowest dimension chosen to be the identity, as discussed in the Introduction.  This is similar to the case of doubly degenerate tricategories above, and is in direct analogy with the 2-dimensional version in which the bicategory of monoidal categories, monoidal functors and monoidal transformations can be found as a full sub-bicategory of the bicategory of icons.

\begin{Defi} \mbox{\hspace{10cm}}
\begin{enumerate}
\item  A \emph{monoidal bicategory} is a degenerate tricategory. 
\item  A \emph{weak monoidal functor}, which we now shorten to monoidal functor, is a weak functor between the corresponding degenerate tricategories.
\end{enumerate}
\end{Defi}

We now define monoidal transformations as a special case of lax transformations where the single object component is the identity, the lax transformation $\alpha$ is actually weak, and where the two modifications $\Pi$ and $M$ are invertible.  The data and axioms presented here use collapsed versions of the transformation diagrams, making use of the left and right unit adjoint equivalences to simplify the diagrams involved.

\begin{Defi}
Let $B,B'$ be monoidal bicategories and $F,G: B \rightarrow B'$ be monoidal functors between them.  A \emph{monoidal transformation} $\alpha:F \Rightarrow G$ consists of
\begin{itemize}
\item a weak transformation $\alpha: F \Rightarrow G$ between the underlying weak functors,
\item an invertible modification as displayed below,
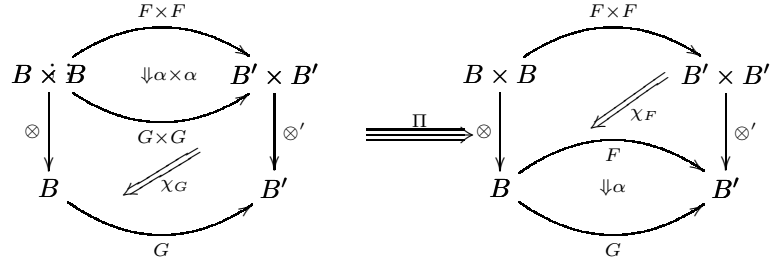
\[
\xy
{\ar@/^1.5pc/^{F \times F} (0,0)*+{B \times B}; (30,0)*+{B' \times B'} };
{\ar@/_1.5pc/_{G \times G} (0,0)*+{B \times B}; (30,0)*+{B' \times B'} };
(16,0)*{\scriptstyle \Downarrow \alpha \times \alpha};
{\ar_{\otimes} (0,0)*+{B \times B}; (0,-15)*+{B} };
{\ar@/_1.5pc/_{G} (0,-15)*+{B}; (30,-15)*+{B'} };
{\ar^{\otimes'} (30,0)*+{B' \times B'}; (30,-15)*+{B'} };
{\ar@{=>}^{\chi_{G}} (20,-10)*{}; (10,-16)*{} };
{\ar@/^1.5pc/^{F \times F} (60,0)*+{B \times B}; (90,0)*+{B' \times B'} };
{\ar_{\otimes} (60,0)*+{B \times B}; (60,-15)*+{B} };
{\ar@/^1.5pc/_{F} (60,-15)*+{B}; (90,-15)*+{B'} };
{\ar@/_1.5pc/_{G} (60,-15)*+{B}; (90,-15)*+{B'} };
(75,-15)*{\scriptstyle \Downarrow \alpha}
{\ar^{\otimes'} (90,0)*+{B' \times B'}; (90,-15)*+{B'} };
{\ar@{=>}^{\chi_{F}} (80,0)*{}; (70,-7)*{} };
{\ar@3{->}^{\Pi} (38,-8)*{}; (52,-8)*{} }
\endxy
\]
\item and an invertible modification as displayed below,
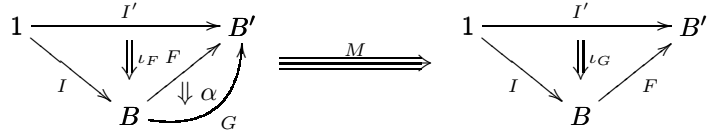
\[
\xy
{\ar^{I'} (0,0)*+{1}; (30,0)*+{B'} };
{\ar_{I} (0,0)*+{1}; (15,-12)*+{B} };
{\ar^{F} (15,-12)*+{B}; (30,0)*+{B'} };
{\ar@/_1.5pc/_{G} (15,-12)*+{B}; (30,0)*+{B'} };
{\ar@{=>}^{\iota_{F}} (15,-2)*{}; (15,-7)*+{} };
(24,-9)*{\Downarrow \alpha};
{\ar^{I'} (60,0)*+{1}; (90,0)*+{B'} };
{\ar_{I} (60,0)*+{1}; (75,-12)*+{B} };
{\ar_{F} (75,-12)*+{B}; (90,0)*+{B'} };
{\ar@{=>}^{\iota_{G}} (75,-2)*{}; (75,-7)*+{} };
{\ar@3{->}^{M} (35,-5)*{}; (55,-5)*{}}
\endxy
\]
\end{itemize}
all subject to the following three axioms.
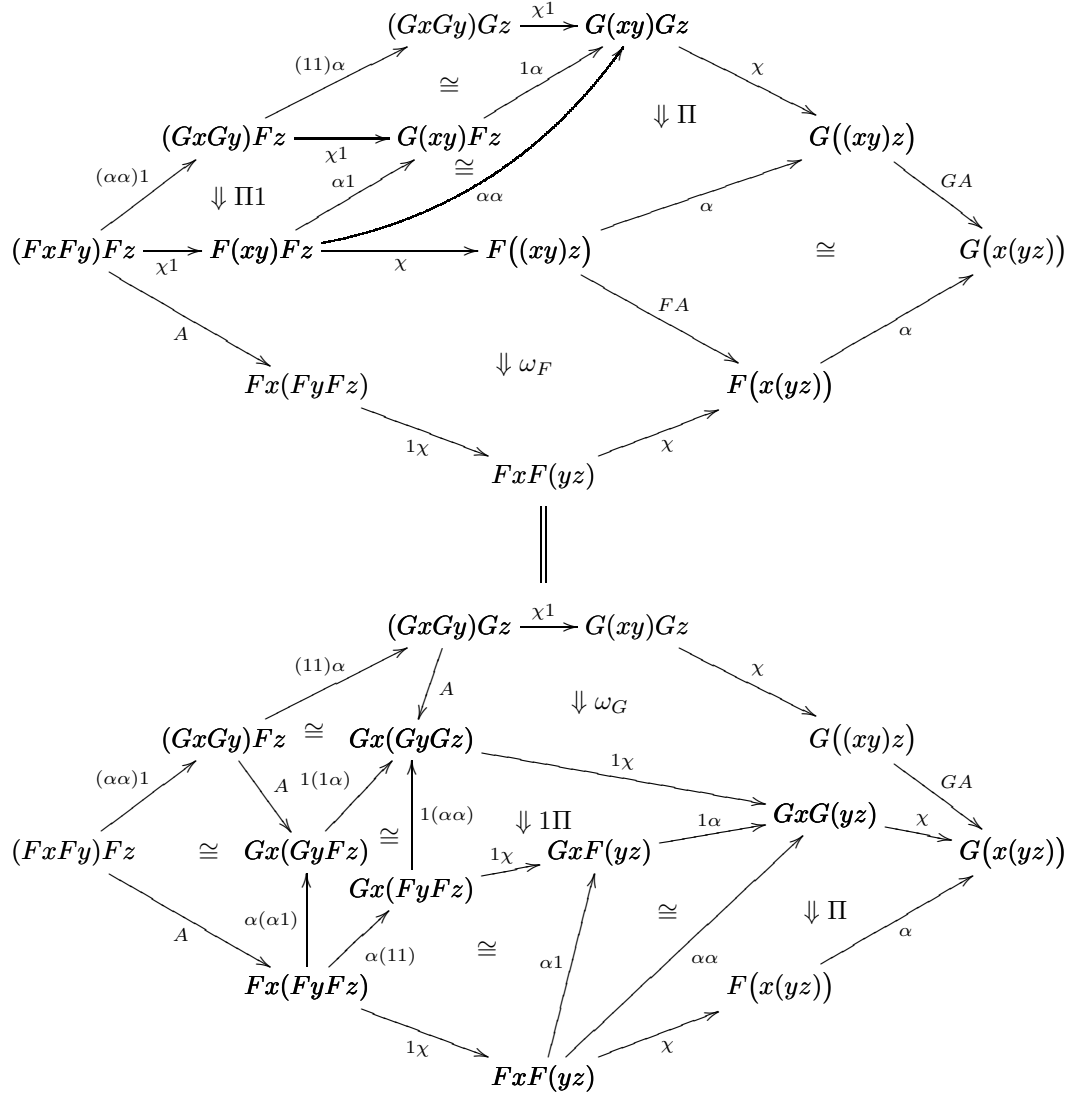
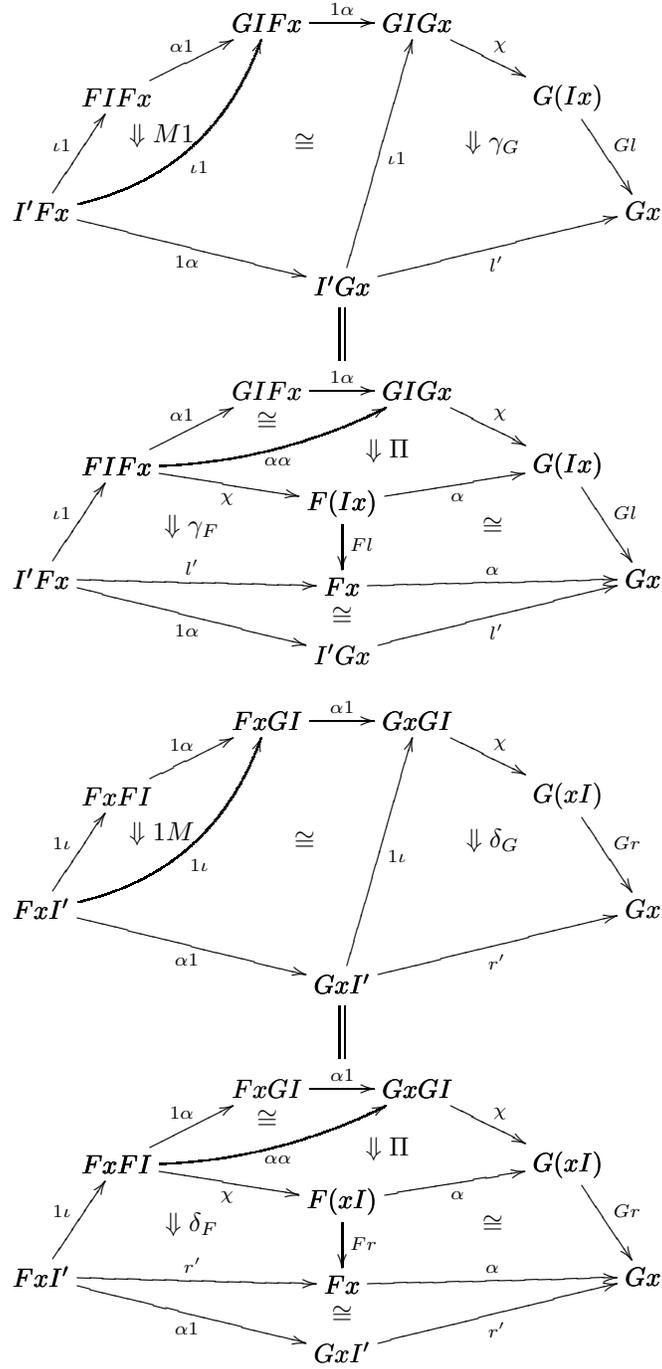
\[
\xy
{\ar^{(\alpha \alpha) 1} (0,0)*+{(FxFy)Fz}; (20,15)*+{(GxGy)Fz} };
{\ar^{(11)\alpha} (20,15)*+{(GxGy)Fz}; (50,30)*+{(GxGy)Gz} };
{\ar^{\chi 1} (50,30)*+{(GxGy)Gz}; (75,30)*+{G(xy)Gz} };
{\ar^{\chi} (75,30)*+{G(xy)Gz}; (105,15)*+{G \big( (xy)z \big)} };
{\ar^{GA} (105,15)*+{G \big( (xy)z \big)}; (125,0)*+{G \big( x(yz) \big) } };
{\ar_{A} (0,0)*+{(FxFy)Fz}; (31,-18)*+{Fx(FyFz)} };
{\ar_{1 \chi} (31,-18)*+{Fx(FyFz)}; (62.5,-30)*+{FxF(yz)} };
{\ar_{\chi} (62.5,-30)*+{FxF(yz)}; (94,-18)*+{F \big( x(yz) \big) } };
{\ar_{\alpha} (94,-18)*+{F \big( x(yz) \big) }; (125,0)*+{G \big( x(yz) \big)} }; 
{\ar_{\chi 1} (0,0)*+{(FxFy)Fz}; (25,0)*+{F(xy)Fz} };
{\ar_{\chi 1} (20,15)*+{(GxGy)Fz}; (50,15)*+{G(xy)Fz} };
{\ar^{\alpha 1} (25,0)*+{F(xy)Fz}; (50,15)*+{G(xy)Fz} };
{\ar^{1 \alpha} (50,15)*+{G(xy)Fz}; (75,30)*+{G(xy)Gz} };
{\ar@/_1.7pc/_{\alpha \alpha} (25,0)*+{F(xy)Fz}; (75,30)*+{G(xy)Gz} };
{\ar_{\chi} (25,0)*+{F(xy)Fz}; (62,0)*+{F \big( (xy)z \big)} };
{\ar^{FA} (62,0)*+{F \big( (xy)z \big)}; (94,-18)*+{F \big( x(yz) \big) } };
{\ar_{\alpha} (62,0)*+{F \big( (xy)z \big)}; (105,15)*+{G \big( (xy)z \big)} };
(50, 22)*{\cong}; (22, 7)*{\Downarrow \Pi 1}; (52,11)*{\cong}; (80, 18)*{\Downarrow \Pi}; (100,0)*{\cong}; (60, -15)*{\Downarrow \omega_{F}};
{\ar^{(\alpha \alpha) 1} (0,-80)*+{(FxFy)Fz}; (20,-65)*+{(GxGy)Fz} };
{\ar^{(11)\alpha} (20,-65)*+{(GxGy)Fz}; (50,-50)*+{(GxGy)Gz} };
{\ar^{\chi 1} (50,-50)*+{(GxGy)Gz}; (75,-50)*+{G(xy)Gz} };
{\ar^{\chi} (75,-50)*+{G(xy)Gz}; (105,-65)*+{G \big( (xy)z \big)} };
{\ar^{GA} (105,-65)*+{G \big( (xy)z \big)}; (125,-80)*+{G \big( x(yz) \big) } };
{\ar_{A} (0,-80)*+{(FxFy)Fz}; (31,-98)*+{Fx(FyFz)} };
{\ar_{1 \chi} (31,-98)*+{Fx(FyFz)}; (62.5,-110)*+{FxF(yz)} };
{\ar_{\chi} (62.5,-110)*+{FxF(yz)}; (94,-98)*+{F \big( x(yz) \big) } };
{\ar_{\alpha} (94,-98)*+{F \big( x(yz) \big) }; (125,-80)*+{G \big( x(yz) \big)} }; 
{\ar^{A} (20,-65)*+{(GxGy)Fz}; (31,-80)*+{Gx(GyFz)} };
{\ar^{\alpha (\alpha 1)} (31,-98)*+{Fx(FyFz)}; (31,-80)*+{Gx(GyFz)} };
{\ar^{A} (50,-50)*+{(GxGy)Gz}; (45,-65)*+{Gx(GyGz)} };
{\ar^{1(1 \alpha)} (31,-80)*+{Gx(GyFz)}; (45,-65)*+{Gx(GyGz)} };
{\ar^{1 \chi} (45,-65)*+{Gx(GyGz)}; (100,-75)*+{GxG(yz)} };
{\ar^{\chi} (100,-75)*+{GxG(yz)}; (125,-80)*+{G \big( x(yz) \big)} }; 
{\ar_{\alpha \alpha} (62.5,-110)*+{FxF(yz)}; (100,-75)*+{GxG(yz)} };
{\ar_{\alpha (11)} (31,-98)*+{Fx(FyFz)}; (45,-85)*+{Gx(FyFz)} };
{\ar_{1(\alpha \alpha)} (45,-85)*+{Gx(FyFz)}; (45,-65)*+{Gx(GyGz)} };
{\ar^{1 \chi} (45,-85)*+{Gx(FyFz)}; (70,-80)*+{GxF(yz)} };
{\ar^{1 \alpha} (70,-80)*+{GxF(yz)}; (100,-75)*+{GxG(yz)} };
{\ar^{\alpha 1} (62.5,-110)*+{FxF(yz)}; (70,-80)*+{GxF(yz)} };
(18,-80)*{\cong}; (32, -64)*{\cong}; (42,-78)*{\cong}; (70,-60)*{\Downarrow \omega_{G}}; (62.5,-76)*{\Downarrow 1 \Pi}; (55,-93)*{\cong}; (79,-88)*{\cong}; (100,-88)*{\Downarrow \Pi};
{\ar@{=} (62.5,-34)*{}; (62.5,-44)*{} }
\endxy
\]
\[
\xy
{\ar^{\iota 1} (0,0)*+{I'Fx}; (10, 15)*+{FIFx} };
{\ar^{\alpha 1} (10, 15)*+{FIFx}; (30, 25)*+{GIFx} };
{\ar^{1 \alpha} (30, 25)*+{GIFx}; (50, 25)*+{GIGx} };
{\ar^{\chi} (50, 25)*+{GIGx}; (70, 15)*+{G(Ix)} };
{\ar^{Gl} (70, 15)*+{G(Ix)}; (80,0)*+{Gx} };
{\ar_{1 \alpha} (0,0)*+{I'Fx}; (40,-10)*+{I'Gx} };
{\ar_{l'} (40,-10)*+{I'Gx}; (80,0)*+{Gx} };
{\ar@/_1.5pc/_{\iota 1} (0,0)*+{I'Fx}; (30, 25)*+{GIFx} };
{\ar_{\iota 1} (40,-10)*+{I'Gx}; (50, 25)*+{GIGx} };
(16, 10)*{\Downarrow M1}; (35,9)*{\cong}; (60,9)*{\Downarrow \gamma_{G}};
{\ar^{\iota 1} (0,-49)*+{I'Fx}; (10, -34)*+{FIFx} };
{\ar^{\alpha 1} (10, -34)*+{FIFx}; (30, -24)*+{GIFx} };
{\ar^{1 \alpha} (30, -24)*+{GIFx}; (50, -24)*+{GIGx} };
{\ar^{\chi} (50, -24)*+{GIGx}; (70, -34)*+{G(Ix)} };
{\ar^{Gl} (70, -34)*+{G(Ix)}; (80,-49)*+{Gx} };
{\ar_{1 \alpha} (0,-49)*+{I'Fx}; (40,-59)*+{I'Gx} };
{\ar_{l'} (40,-59)*+{I'Gx}; (80,-49)*+{Gx} };
{\ar@/_0.7pc/_{\alpha \alpha} (10, -34)*+{FIFx}; (50, -24)*+{GIGx} };
{\ar_{\chi} (10, -34)*+{FIFx}; (40,-39)*+{F(Ix)} };
{\ar_{\alpha}  (40,-39)*+{F(Ix)}; (70, -34)*+{G(Ix)} };
{\ar^{l'} (0,-49)*+{I'Fx}; (40, -50)*+{Fx} };
{\ar^{\alpha} (40, -50)*+{Fx}; (80,-49)*+{Gx} };
{\ar^{Fl} (40,-39)*+{F(Ix)}; (40, -50)*+{Fx} };
(30,-28)*{\cong}; (46,-32)*{\Downarrow \Pi}; (20, -42)*{\Downarrow \gamma_{F}}; (60,-42)*{\cong}; (40,-54)*{\cong};
{\ar@{=} (40,-13)*{}; (40,-20)*{} }
\endxy
\]

\[
\xy
{\ar^{1 \iota} (0,0)*+{Fx I'}; (10, 15)*+{FxFI} };
{\ar^{1\alpha} (10, 15)*+{FxFI}; (30, 25)*+{FxGI} };
{\ar^{\alpha 1} (30, 25)*+{FxGI}; (50, 25)*+{GxGI} };
{\ar^{\chi} (50, 25)*+{GxGI}; (70, 15)*+{G(xI)} };
{\ar^{Gr} (70, 15)*+{G(xI)}; (80,0)*+{Gx} };
{\ar_{ \alpha 1} (0,0)*+{Fx I'}; (40,-10)*+{Gx I'} };
{\ar_{r'} (40,-10)*+{GxI'}; (80,0)*+{Gx} };
{\ar@/_1.5pc/_{1 \iota} (0,0)*+{Fx I'}; (30, 25)*+{FxGI} };
{\ar_{1 \iota} (40,-10)*+{GxI'}; (50, 25)*+{GxGI} };
(16, 10)*{\Downarrow 1 M}; (35,9)*{\cong}; (60,9)*{\Downarrow \delta_{G}};
{\ar^{1 \iota} (0,-49)*+{Fx I'}; (10, -34)*+{FxFI} };
{\ar^{1 \alpha} (10, -34)*+{FxFI}; (30, -24)*+{FxGI} };
{\ar^{\alpha 1} (30, -24)*+{FxGI}; (50, -24)*+{GxGI} };
{\ar^{\chi} (50, -24)*+{GxGI}; (70, -34)*+{G(xI)} };
{\ar^{Gr} (70, -34)*+{G(xI)}; (80,-49)*+{Gx} };
{\ar_{\alpha 1} (0,-49)*+{FxI'}; (40,-59)*+{GxI'} };
{\ar_{r'} (40,-59)*+{GxI'}; (80,-49)*+{Gx} };
{\ar@/_0.7pc/_{\alpha \alpha} (10, -34)*+{FxFI}; (50, -24)*+{GxGI} };
{\ar_{\chi} (10, -34)*+{FxFI}; (40,-39)*+{F(xI)} };
{\ar_{\alpha}  (40,-39)*+{F(xI)}; (70, -34)*+{G(xI)} };
{\ar^{r'} (0,-49)*+{FxI'}; (40, -50)*+{Fx} };
{\ar^{\alpha} (40, -50)*+{Fx}; (80,-49)*+{Gx} };
{\ar^{Fr} (40,-39)*+{F(xI)}; (40, -50)*+{Fx} };
(30,-28)*{\cong}; (46,-32)*{\Downarrow \Pi}; (20, -42)*{\Downarrow \delta_{F}}; (60,-42)*{\cong}; (40,-54)*{\cong};
{\ar@{=} (40,-13)*{}; (40,-20)*{} }
\endxy
\]
Note that in the previous diagram we have written $\delta_{F}$ and $\delta_{G}$ when in fact their mates are used.
\end{Defi}

We now define monoidal modifications between monoidal bicategories in a similar fashion, as a special case of lax modifications with the component at the single object being given by an identity.  Using the left and right unit adjoint equivalences, we are then able to simplify the diagrams to those given below.  

\begin{Defi}
Let $\alpha, \beta:F \Rightarrow G$ be monoidal transformations between monoidal functors.  A \emph{monoidal modification} $m: \alpha \Rrightarrow \beta$ consists of a modification $m: \alpha \Rrightarrow \beta$ between the underlying transformations such that the following two axioms hold.
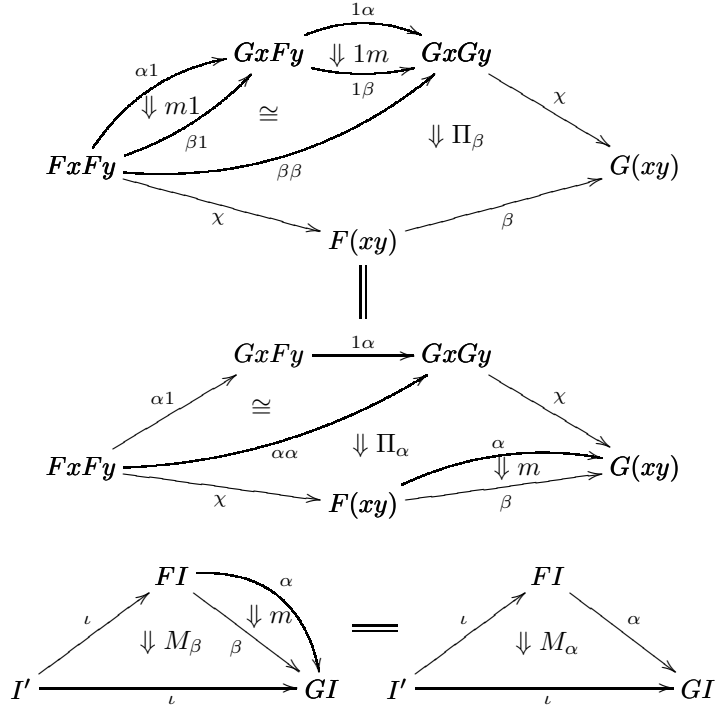
\[
\xy
{\ar@/^1pc/^{\alpha 1} (0,0)*+{FxFy}; (25,15)*+{GxFy} };
{\ar@/^1pc/^{1 \alpha} (25,15)*+{GxFy}; (50, 15)*+{GxGy} };
{\ar^{\chi} (50, 15)*+{GxGy};(75,0)*+{G(xy)} };
{\ar_{\chi} (0,0)*+{FxFy}; (37.5, -10)*+{F(xy)} };
{\ar_{\beta} (37.5, -10)*+{F(xy)}; (75,0)*+{G(xy)} };
{\ar@/_0.6pc/_{\beta 1} (0,0)*+{FxFy}; (25,15)*+{GxFy} };
{\ar@/_0.6pc/_{1 \beta} (25,15)*+{GxFy}; (50, 15)*+{GxGy} };
{\ar@/_1.5pc/_{\beta \beta} (0,0)*+{FxFy}; (50, 15)*+{GxGy} };
(12, 8)*{\Downarrow m 1}; (37,15)*{\Downarrow 1 m}; (25, 7)*{\cong}; (50,4)*{\Downarrow \Pi_{\beta}};
{\ar^{\alpha 1} (0,-40)*+{FxFy}; (25,-25)*+{GxFy} };
{\ar^{1 \alpha} (25,-25)*+{GxFy}; (50, -25)*+{GxGy} };
{\ar^{\chi} (50, -25)*+{GxGy};(75,-40)*+{G(xy)} };
{\ar_{\chi} (0,-40)*+{FxFy}; (37.5, -45)*+{F(xy)} };
{\ar_{\beta} (37.5, -45)*+{F(xy)}; (75,-40)*+{G(xy)} };
{\ar@/_1pc/_{\alpha \alpha} (0,-40)*+{FxFy}; (50, -25)*+{GxGy} };
{\ar@/^1pc/^{\alpha} (37.5, -45)*+{F(xy)}; (75,-40)*+{G(xy)} };
(24,-32)*{\cong}; (40,-37)*{\Downarrow \Pi_{\alpha}}; (58,-40)*{\Downarrow m};
{\ar@{=} (37.5, -13)*{}; (37.5,-20)*{} }
\endxy
\]

\[
\xy
{\ar^{\iota} (0,0)*+{I'}; (20,15)*+{FI} };
{\ar@/^1.5pc/^{\alpha} (20,15)*+{FI}; (40,0)*+{GI} };
{\ar_{\beta} (20,15)*+{FI}; (40,0)*+{GI} };
{\ar_{\iota} (0,0)*+{I'}; (40,0)*+{GI} };
(20, 6)*{\Downarrow M_{\beta}}; (33,9.5)*{\Downarrow m};
{\ar^{\iota} (50,0)*+{I'}; (70,15)*+{FI} };
{\ar^{\alpha} (70,15)*+{FI}; (90,0)*+{GI} };
{\ar_{\iota} (50,0)*+{I'}; (90,0)*+{GI} };
(70, 6)*{\Downarrow M_{\alpha}};
{\ar@{=} (44,8)*{}; (50,8)*{} } 
\endxy
\]
\end{Defi}

The rest of this section will be devoted to defining the structure of the tricategory $\mathbf{MonBicat}$ whose 0-cells are monoidal bicategories, 1-cells are monoidal functors, 2-cells are monoidal transformations, and 3-cells are monoidal modifications.  We begin by defining the hom-bicategories for this tricategory; note that composition is not inherited directly from $\mathbf{Tricat}$ but can be thought of as a ``hybrid'' of the respective structures of \cat{Tricat} and \cat{Bicat}.  

For 1-cell composition, consider monoidal transformations $\alpha:F \Rightarrow G$ and $\beta:G \Rightarrow H$.  We define a monoidal transformation $\beta \alpha$ as follows:
\begin{itemize}
\item its underlying transformation is the composite $\beta \alpha$,
\item the invertible modification $\Pi_{\beta \alpha}$ has component at $(X,Y)$ given by the diagram below,
\[
\xy
{\ar^{(\beta \alpha) \otimes (\beta \alpha)} (0,0)*+{\scriptstyle FX \otimes FY}; (40,0)*+{\scriptstyle HX \otimes HY} };
{\ar^{\chi_{H}}  (40,0)*+{\scriptstyle HX \otimes HY}; (40,-30)*+{\scriptstyle H(X \otimes Y)} };
{\ar_{\chi_{F}} (0,0)*+{\scriptstyle FX \otimes FY}; (0,-30)*+{\scriptstyle F(X \otimes Y)} };
{\ar_{\alpha} (0,-30)*+{\scriptstyle F(X \otimes Y)}; (20,-30)*+{\scriptstyle G(X \otimes Y)} };
{\ar_{\beta} (20,-30)*+{\scriptstyle G(X \otimes Y)}; (40,-30)*+{\scriptstyle H(X \otimes Y)} };
{\ar_{\alpha \otimes \alpha} (0,0)*+{\scriptstyle FX \otimes FY}; (20,-15)*+{\scriptstyle GX \otimes GY} };
{\ar_{\beta \otimes \beta} (20,-15)*+{\scriptstyle GX \otimes GY}; (40,0)*+{\scriptstyle HX \otimes HY} };
{\ar_{\chi_{G}} (20,-15)*+{\scriptstyle GX \otimes GY}; (20,-30)*+{\scriptstyle G(X \otimes Y)} };
(20, -7)*{\cong}; (9, -19)*{\Downarrow \Pi_{\alpha}}; (30, -19)*{\Downarrow \Pi_{\beta}}
\endxy
\]
\item and the invertible modification $M_{\beta \alpha}$ is given by the diagram below.
\[
\xy
{\ar^{\iota_{F}} (0,0)*+{I'}; (15,15)*+{FI} };
{\ar^{\alpha} (15,15)*+{FI}; (30, 15)*+{GI} };
{\ar^{\beta} (30, 15)*+{GI}; (45,0)*+{HI} };
{\ar@/_0.6pc/_{\iota_{G}} (0,0)*+{I'}; (30, 15)*+{GI} };
{\ar_{\iota_{H}} (0,0)*+{I'}; (45,0)*+{HI} };
(17, 10)*+{\Downarrow M_{\alpha}}; (30,6)*+{\Downarrow M_{\beta}}
\endxy
\]
\end{itemize}
The three axioms are easily checked by a simple diagram chase.

For identity 1-cells, consider a monoidal functor $F$.  Then the identity transformation $u: F \Rightarrow F$ can be equipped with the structure of a monoidal transformation with both $\Pi_{u}$ and $M_{u}$ being given by unique coherence isomorphisms.  The axioms follow immediately from the coherence theorem for tricategories.

For vertical 2-cell composition, consider monoidal modifications  $m: \alpha \Rrightarrow \beta$ and $n:\beta \Rrightarrow \gamma$.  Then we can check that the composite $nm: \alpha \Rrightarrow \gamma$ in $\mathbf{Bicat}$ is in fact monoidal, and likewise the identity.

For horizontal 2-cell composition, consider monoidal modifications as displayed below.
\[
\xy
{\ar@/^6pc/^{F} (0,0)*+{X}; (50,0)*+{Y} };
{\ar^{G} (0,0)*+{X}; (50,0)*+{Y} };
{\ar@/_6pc/_{H} (0,0)*+{X}; (50,0)*+{Y} };
{\ar@/_1pc/_{\alpha} (14,14)*{}; (14,3)*{} };
{\ar@/^1pc/^{\beta} (36,14)*{}; (36,3)*{} };
{\ar@/_1pc/_{\gamma} (14,-3)*{}; (14,-14)*{} };
{\ar@/^1pc/^{\delta} (36,-3)*{}; (36,-14)*{} };
{\ar@3{->}^{m} (15,8)*{}; (35,8)*{} };
{\ar@3{->}_{n} (15,-8)*{}; (35,-8)*{} }
\endxy
\]
Then we can check that the composite $n*m:\gamma \alpha \Rrightarrow \delta \beta$ in $\mathbf{Bicat}$ is in fact monoidal, and that this composition is functorial.

For coherence isomorphisms in the hom-bicategories, consider monoidal transformations $\alpha: F \Rightarrow G$, $\beta:G \Rightarrow H$, and $\gamma:H \Rightarrow J$.
\begin{itemize}
\item  Let $r: \alpha u_{F} \Rrightarrow \alpha$ be the modification with component at $X$ the right unit isomorphism $r_{\alpha_{X}}$.  It follows from coherence for tricategories that $r$ and $r^{-1}$ are monoidal.
\item   Let $l: u_{G} \alpha \Rrightarrow \alpha$ be the modification with component at $X$ the left unit isomorphism $l_{\alpha_{X}}$.  Observe as above that this modification and its inverse $l^{-1}$ are monoidal.
\item  Let $a: (\gamma \beta) \alpha \Rrightarrow \gamma (\beta \alpha)$ be the modification with component at $X$ the associativity isomorphism $a_{\gamma_{X} \beta_{X} \alpha_{X}}$ is monoidal.  Observe as above that this modification and its inverse $a^{-1}$ are monoidal.
\end{itemize}

\begin{thm}
The above structure defines a bicategory $\mathbf{MonBicat}(X,Y)$.
\end{thm}
\begin{proof}
The axioms follow from the bicategory axioms in $Y$.
\end{proof}

We next define composition along 0-cells for the tricategory $\mathbf{MonBicat}$, which we will denote $\boxtimes$; we simply extend the definition of composition in the tricategory $\mathbf{Bicat}$ which we now recall.  Consider functors, transformations, and modifications as below.
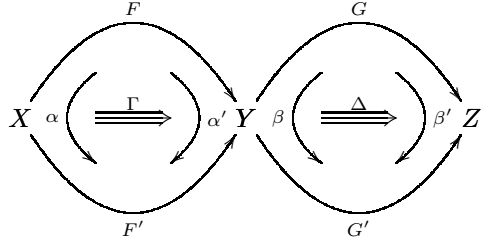
\[
\xy
{\ar@/^3pc/^{F} (0,0)*+{X}; (30,0)*+{Y} };
{\ar@/_3pc/_{F'} (0,0)*+{X}; (30,0)*+{Y} };
{\ar@/_0.9pc/_{\alpha} (10,6)*{}; (10,-6)*{} };
{\ar@/^0.9pc/^{\alpha'} (20,6)*{}; (20,-6)*{} };
{\ar@3{->}^{\Gamma} (10,0)*+{}; (20,0)*+{} };
{\ar@/^3pc/^{G} (30,0)*+{Y}; (60,0)*+{Z} };
{\ar@/_3pc/_{G'} (30,0)*+{Y}; (60,0)*+{Z} };
{\ar@/_0.9pc/_{\beta} (40,6)*{}; (40,-6)*{} };
{\ar@/^0.9pc/^{\beta'} (50,6)*{}; (50,-6)*{} };
{\ar@3{->}^{\Delta} (40,0)*+{}; (50,0)*+{} };
\endxy
\]
Then we have the following formulae.
\[
\begin{array}{c}
G \otimes F : = GF \\
\beta \otimes \alpha := (G' * \alpha) \circ (\beta * F) \\
(\Delta \otimes \Gamma)_{x} := G'\Gamma_{x} * \Delta_{Fx}
\end{array}
\]
Now suppose all of the above data are monoidal.
\begin{enumerate}
\item The composite $G \boxtimes F$ is the composite of the functors of the underlying degenerate tricategories.
\item The composite $\beta \boxtimes \alpha$ has underlying transformation $\beta \otimes \alpha$ as above together with
\begin{itemize}
\item invertible modification $\Pi$ given by the diagram below, and
\[
\xy
{\ar^{\beta \otimes \beta} (0,0)*+{GFX \otimes GFY}; (35,0)*+{G'FX \otimes G'FY} };
{\ar^{G'\alpha \otimes G'\alpha}  (35,0)*+{G'FX \otimes G'FY}; (70,0)*+{G'F'X \otimes G'F'Y} };
{\ar@/^2pc/^{(\beta \boxtimes \alpha) \otimes (\beta \boxtimes \alpha)} (0,0)*+{GFX \otimes GFY}; (70,0)*+{G'F'X \otimes G'F'Y} };
{\ar_{\chi_{G}} (0,0)*+{GFX \otimes GFY}; (0,-15)*+{G(FX \otimes FY)} };
{\ar_{G \chi_{F}} (0,-15)*+{G(FX \otimes FY)}; (0,-30)*+{GF(X \otimes Y)} };
{\ar_{\beta} (0,-15)*+{G(FX \otimes FY)}; (35,-15)*+{G'(FX \otimes FY)} };
{\ar_{\beta} (0,-30)*+{GF(X \otimes Y)}; (35,-30)*+{G'F(X \otimes Y)} };
{\ar_{\chi_{G'}} (35,0)*+{G'FX \otimes G'FY}; (35,-15)*+{G'(FX \otimes FY)} };
{\ar_{G'\chi_{F}} (35,-15)*+{G'(FX \otimes FY)}; (35,-30)*+{G'F(X \otimes Y)} };
{\ar_{G'(\alpha \otimes \alpha)} (35,-15)*+{G'(FX \otimes FY)}; (70,-15)*+{G'(F'X \otimes F'Y)} };
{\ar_{G' \alpha} (35,-30)*+{G'F(X \otimes Y)}; (70,-30)*+{G'F'(X \otimes Y)} };
{\ar^{\chi_{G'}} (70,0)*+{G'F'X \otimes G'F'Y}; (70,-15)*+{G'(F'X \otimes F'Y)} };
{\ar^{G'\chi_{F'}} (70,-15)*+{G'(F'X \otimes F'Y)}; (70,-30)*+{G'F'(X \otimes Y)} };
(35,4)*{\cong}; (17,-8)*{\Downarrow \Pi_{\beta}}; (17,-23)*{\cong}; (52,-8)*{\cong}; (52,-23)*{\Downarrow G'\Pi_{\alpha}}
\endxy
\]
\item invertible modification $M$ given by the diagram below.
\[
\xy
{\ar^{\iota_{G}} (0,0)*+{I''}; (22,0)*+{GI'} };
{\ar^{\G\iota_{F}} (22,0)*+{GI'}; (44, 0)*+{GFI} };
{\ar^{\beta_{FI}} (44, 0)*+{GFI}; (66,0)*+{G'FI} };
{\ar^{G'\alpha_{I}} (66,0)*+{G'FI}; (88,0)*+{G'F'I} };
{\ar_{\beta_{I'}} (22,0)*+{GI'}; (44,-13)*+{G'I'} };
{\ar_{G\iota_{F}} (44,-13)*+{G'I'}; (66,0)*+{G'FI} };
{\ar@/_1pc/_{\iota_{G'}} (0,0)*+{I''}; (44,-13)*+{G'I'} };
{\ar@/_1pc/_{G'\iota_{F}} (44,-13)*+{G'I'}; (88,0)*+{G'F'I} };
(44,-6)*{\cong}; (22,-6)*{\Downarrow M_{\beta}}; (67,-6)*{\Downarrow G'M_{\alpha}}
\endxy
\]
\end{itemize}
\item The modification $\Delta \otimes \Gamma$ is a monoidal modification, so we can put $\Delta \boxtimes \Gamma = \Delta \otimes \Gamma$.
\end{enumerate}

\begin{thm}
The assignments above extend to a functor
\[
\boxtimes: \mathbf{MonBicat}(Y,Z) \times \mathbf{MonBicat}(X,Y) \rightarrow \mathbf{MonBicat}(X,Z).
\]
\end{thm}
\begin{proof}
The constraint modifications are the same as those given in my \cite{gur2}; we need only check that they are monoidal modifications, which is accomplished by a lengthy, but routine, diagram chase.  The functor axioms follow from coherence and the transformation axioms.
\end{proof}

\noindent We now define units for the composition $\boxtimes$.  


\begin{prop}
Let $X$ be a monoidal bicategory.  There is a functor $I_{X}: 1 \rightarrow \mathbf{MonBicat}(X,X)$ whose value on the single object is the identity monoidal functor and whose value on the single 1-cell is the identity monoidal transformation.
\end{prop}
\begin{proof}
Functoriality determines that the value on the single 2-cell is the identity.  The unit constraint is the identity, and the composition constraint is given by the left (or right) unit isomorphism in $X$, which we have already determined is a monoidal modification.  The axioms then follow from coherence.
\end{proof}

We now define the adjoint equivalences
\[
\begin{array}{c}
\mathbf{a}: \boxtimes \circ (\boxtimes \times 1) \Rightarrow \boxtimes \circ (1 \times \boxtimes) \\
\mathbf{l}: \boxtimes \circ (I_{X} \times 1) \Rightarrow 1 \\
\mathbf{r}: \boxtimes \circ (1 \times I_{X}) \Rightarrow 1.
\end{array}
\]
The underlying adjoint equivalences of transformations are all the same as the relevant adjoint equivalences in $\mathbf{Bicat}$.  It remains to provide the component modifications, check that these choices give monoidal transformations, check that the unit and counit modifications are monoidal, and check the triangle identities.  All the cells involved are coherence cells, and we can use coherence for tricategories to check that all necessary diagrams commute.

\begin{thm}
There is a tricategory $\mathbf{MonBicat}$ with
\begin{itemize}
\item 0-cells monoidal bicategories;
\item hom-bicategories given by the bicategories $\mathbf{MonBicat}(X,Y)$ defined above;
\item composition functor given by $\boxtimes$;
\item unit given by the functor $I_{X}:1 \rightarrow \mathbf{MonBicat}(X,X)$;
\item adjoint equivalences $\mathbf{a}, \mathbf{l}, \mathbf{r}$ as above; and
\item invertible modifications $\pi, \lambda, \rho, \mu$ with each modification having components given by unique coherence cells in the target bicategory.
\end{itemize}
Furthermore, the obvious forgetful functor $\mathbf{MonBicat} \rightarrow \mathbf{Bicat}$ is a strict functor between tricategories.
\end{thm}
\begin{proof}
The tricategory axioms follow from coherence for bicategories.  The fact that the modifications above are monoidal follows from coherence for tricategories.
\end{proof}

\appendix

\section{Appendix: diagrams}

Here we include all the diagrams that were omitted from the text itself.

\subsection{Doubly degenerate tricategories}
We will write the monoidal structure $\otimes$ as concatenation, invoking coherence in order to ignore association.  The monoidal functor $\boxtimes$ will be written as $\cdot$ to save space, and we will enclose terms involving $\boxtimes$ with square brackets when necessary to avoid an excess of parentheses.  Isomorphisms will remain largely unmarked as there will always be an obvious choice.

The dual pair $\bs{A}$ and its attendant natural isomorphisms must satisfy the following axioms. 
\[
\xy
{\ar (0,0)*+{\scriptscriptstyle (A[(X_{1} \cdot Y_{1}) \cdot Z_{1}])[(X_{2} \cdot Y_{2}) \cdot Z_{2}]}; (40,0)*+{\scriptscriptstyle ([X_{1} \cdot (Y_{1} \cdot Z_{1})]A)[(X_{2} \cdot Y_{2}) \cdot Z_{2}]} };
{\ar (40,0)*+{\scriptscriptstyle ([X_{1} \cdot (Y_{1} \cdot Z_{1})]A)[(X_{2} \cdot Y_{2}) \cdot Z_{2}]}; (80,0)*+{\scriptscriptstyle [X_{1} \cdot (Y_{1} \cdot Z_{1})](A[(X_{2} \cdot Y_{2}) \cdot Z_{2}])} };
{\ar (80,0)*+{\scriptscriptstyle [X_{1} \cdot (Y_{1} \cdot Z_{1})](A[(X_{2} \cdot Y_{2}) \cdot Z_{2}])}; (80,-18)*+{\scriptscriptstyle [X_{1} \cdot (Y_{1} \cdot Z_{1})]([X_{2} \cdot (Y_{2} \cdot Z_{2})]A)} };
{\ar (80,-18)*+{\scriptscriptstyle [X_{1} \cdot (Y_{1} \cdot Z_{1})]([X_{2} \cdot (Y_{2} \cdot Z_{2})]A)}; (80,-36)*+{\scriptscriptstyle ([X_{1} \cdot (Y_{1} \cdot Z_{1})][X_{2} \cdot (Y_{2} \cdot Z_{2})])A} };
{\ar (80,-36)*+{\scriptscriptstyle ([X_{1} \cdot (Y_{1} \cdot Z_{1})][X_{2} \cdot (Y_{2} \cdot Z_{2})])A}; (80,-54)*+{\scriptscriptstyle [(X_{1}X_{2}) \cdot \Big( (Y_{1}Y_{2}) \cdot (Z_{1}Z_{2}) \Big) ]A} };
{\ar (0,0)*+{\scriptscriptstyle (A[(X_{1} \cdot Y_{1}) \cdot Z_{1}])[(X_{2} \cdot Y_{2}) \cdot Z_{2}]}; (0,-27)*+{\scriptscriptstyle A \Big( [(X_{1} \cdot Y_{1}) \cdot Z_{1}][(X_{2} \cdot Y_{2}) \cdot Z_{2}] \Big) } };
{\ar (0,-27)*+{\scriptscriptstyle A \Big( [(X_{1} \cdot Y_{1}) \cdot Z_{1}][(X_{2} \cdot Y_{2}) \cdot Z_{2}] \Big) }; (0,-54)*+{\scriptscriptstyle  A[\Big( (X_{1}X_{2}) \cdot (Y_{1}Y_{2}) \Big) \cdot (Z_{1}Z_{2})]} };
{\ar  (0,-54)*+{\scriptscriptstyle  A[\Big( (X_{1}X_{2}) \cdot (Y_{1}Y_{2}) \Big) \cdot (Z_{1}Z_{2})]}; (80,-54)*+{\scriptscriptstyle [(X_{1}X_{2}) \cdot \Big( (Y_{1}Y_{2}) \cdot (Z_{1}Z_{2}) \Big) ]A} }
\endxy
\]
\[
\xy
{\ar (0,0)*+{A}; (30,0)*+{AU} };
{\ar (30,0)*+{AU}; (60, 0)*+{A[(U \cdot U) \cdot U]} };
{\ar (60, 0)*+{A[(U \cdot U) \cdot U]}; (60,-18)*+{[U \cdot (U \cdot U)]A} };
{\ar (0,0)*+{A}; (0,-18)*+{UA} };
{\ar (0,-18)*+{UA}; (60,-18)*+{[U \cdot (U \cdot U)]A} };
\endxy
\]

Similar axioms hold for $A^{\cdot}$, and the dual pairs $\bs{L}$ and $\bs{R}$.

The three tricategory axioms now take the following form, where the labels $\pi, \mu, \lambda, \rho$ indicate that the arrow is built up from that constraint using $\otimes$ and $\boxtimes$; unmarked arrows are given by naturality or unique coherence isomorphisms.
\[
\xy
{\ar (0,0)*+{\scriptscriptstyle [U \cdot (U \cdot A)][U \cdot A] A [(U \cdot A) \cdot U][A \cdot U][(A \cdot U) \cdot U]}; (60,0)*+{\scriptscriptstyle [U \cdot (U \cdot A)][U \cdot A][U \cdot (A \cdot U)]A [A \cdot U][(A \cdot U) \cdot U]} };
{\ar^{\pi} (60,0)*+{\scriptscriptstyle [U \cdot (U \cdot A)][U \cdot A][U \cdot (A \cdot U)]A [A \cdot U][(A \cdot U) \cdot U]}; (60,-12)*+{\scriptscriptstyle [U \cdot A][U \cdot A]A[A \cdot U][(A \cdot U) \cdot U]} };
{\ar^{\pi} (60,-12)*+{\scriptscriptstyle [U \cdot A][U \cdot A]A[A \cdot U][(A \cdot U) \cdot U]}; (60,-24)*+{\scriptscriptstyle [U \cdot A]AA[(A \cdot U) \cdot U]} };
{\ar (60,-24)*+{\scriptscriptstyle [U \cdot A]AA[(A \cdot U) \cdot U]}; (60,-36)*+{\scriptscriptstyle [U \cdot A]A[A \cdot (U \cdot U)]A} };
{\ar (60,-36)*+{\scriptscriptstyle [U \cdot A]A[A \cdot (U \cdot U)]A}; (60,-48)*+{\scriptscriptstyle [U \cdot A]A[A \cdot U]A} };
{\ar^{\pi} (60,-48)*+{\scriptscriptstyle [U \cdot A]A[A \cdot U]A}; (60,-60)*+{\scriptscriptstyle AAA} };
{\ar_{\pi} (0,0)*+{\scriptscriptstyle [U \cdot (U \cdot A)][U \cdot A] A [(U \cdot A) \cdot U][A \cdot U][(A \cdot U) \cdot U]}; (0,-15)*+{\scriptscriptstyle [U \cdot (U \cdot A)][U \cdot A]A [A \cdot U][A \cdot U]} };
{\ar_{\pi} (0,-15)*+{\scriptscriptstyle [U \cdot (U \cdot A)][U \cdot A]A [A \cdot U][A \cdot U]}; (0,-30)*+{\scriptscriptstyle [U \cdot (U \cdot A)]AA[A \cdot U]} };
{\ar (0,-30)*+{\scriptscriptstyle [U \cdot (U \cdot A)]AA[A \cdot U]}; (0,-45)*+{\scriptscriptstyle A[(U \cdot U) \cdot A]A[A \cdot U]} };
{\ar (0,-45)*+{\scriptscriptstyle A[(U \cdot U) \cdot A]A[A \cdot U]}; (0,-60)*+{\scriptscriptstyle A [U \cdot A] A [A \cdot U]} };
{\ar_{\pi} (0,-60)*+{\scriptscriptstyle A [U \cdot A] A [A \cdot U]}; (60,-60)*+{\scriptscriptstyle AAA} };
\endxy
\]

\[
\xy
{\ar (0,0)*+{\scriptscriptstyle A[(U \cdot L) \cdot U][A \cdot U][(R^{\cdot} \cdot U) \cdot U]}; (7,10)*+{\scriptscriptstyle [U \cdot (L \cdot U)]A[A \cdot U][(R^{\cdot} \cdot U) \cdot U]} };
{\ar^{\lambda} (7,10)*+{\scriptscriptstyle [U \cdot (L \cdot U)]A[A \cdot U][(R^{\cdot} \cdot U) \cdot U]}; (14,20)*+{\scriptscriptstyle [U \cdot (LA)]A[A \cdot U][(R^{\cdot} \cdot U) \cdot U]} };
{\ar (14,20)*+{\scriptscriptstyle [U \cdot (LA)]A[A \cdot U][(R^{\cdot} \cdot U) \cdot U]}; (56,20)*+{\scriptscriptstyle [U \cdot L][U \cdot A]A[A \cdot U][(R^{\cdot} \cdot U) \cdot U]} };
{\ar^{\pi} (56,20)*+{\scriptscriptstyle [U \cdot L][U \cdot A]A[A \cdot U][(R^{\cdot} \cdot U) \cdot U]}; (63,10)*+{\scriptscriptstyle [U \cdot L]AA[(R^{\cdot} \cdot U) \cdot U]} };
{\ar (63,10)*+{\scriptscriptstyle [U \cdot L]AA[(R^{\cdot} \cdot U) \cdot U]}; (70,0)*+{\scriptscriptstyle [U \cdot L]A[R^{\cdot} \cdot (U \cdot U)]A} };
{\ar (70,0)*+{\scriptscriptstyle [U \cdot L]A[R^{\cdot} \cdot (U \cdot U)]A}; (70,-15)*+{\scriptscriptstyle [U \cdot L]A[R^{\cdot} \cdot U]A} };
{\ar^{\mu} (70,-15)*+{\scriptscriptstyle [U \cdot L]A[R^{\cdot} \cdot U]A}; (70,-30)*+{\scriptscriptstyle UA} };
{\ar (70,-30)*+{\scriptscriptstyle UA}; (70,-45)*+{\scriptscriptstyle A} };
{\ar (0,0)*+{\scriptscriptstyle A[(U \cdot L) \cdot U][A \cdot U][(R^{\cdot} \cdot U) \cdot U]}; (0,-15)*+{\scriptscriptstyle A[ \big( (U \cdot L)A(R^{\cdot} \cdot U) \big) \cdot (UUU)]} };
{\ar (0,-15)*+{\scriptscriptstyle A[ \big( (U \cdot L)A(R^{\cdot} \cdot U) \big) \cdot (UUU)]}; (0,-30)*+{\scriptscriptstyle A[ \big( (U \cdot L)A(R^{\cdot} \cdot U) \big) \cdot U]} };
{\ar_{\mu} (0,-30)*+{\scriptscriptstyle A[ \big( (U \cdot L)A(R^{\cdot} \cdot U) \big) \cdot U]}; (0,-45)*+{\scriptscriptstyle A[U \cdot U]} };
{\ar (0,-45)*+{\scriptscriptstyle A[U \cdot U]}; (35,-45)*+{\scriptscriptstyle AU} };
{\ar (35,-45)*+{\scriptscriptstyle AU}; (70,-45)*+{\scriptscriptstyle A} };
\endxy
\]
\[
\xy
{\ar (0,0)*+{\scriptscriptstyle [U \cdot (U \cdot L)][U \cdot A][U \cdot (R^{\cdot} \cdot U)]A}; (7,10)*+{\scriptscriptstyle [U \cdot (U \cdot L)][U \cdot A]A[(U \cdot R^{\cdot}) \cdot U]} };
{\ar^{\rho} (7,10)*+{\scriptscriptstyle [U \cdot (U \cdot L)][U \cdot A]A[(U \cdot R^{\cdot}) \cdot U]}; (14,20)*+{\scriptscriptstyle [U \cdot (U \cdot L)][U \cdot A]A[(AR^{\cdot}) \cdot U]} };
{\ar (14,20)*+{\scriptscriptstyle [U \cdot (U \cdot L)][U \cdot A]A[(AR^{\cdot}) \cdot U]}; (56,20)*+{\scriptscriptstyle [U \cdot (U \cdot L)][U \cdot A]A[A \cdot U][R^{\cdot} \cdot U]} };
{\ar^{\pi} (56,20)*+{\scriptscriptstyle [U \cdot (U \cdot L)][U \cdot A]A[A \cdot U][R^{\cdot} \cdot U]}; (63,10)*+{\scriptscriptstyle [U \cdot (U \cdot L)]AA [R^{\cdot} \cdot U]} };
{\ar (63,10)*+{\scriptscriptstyle [U \cdot (U \cdot L)]AA [R^{\cdot} \cdot U]}; (70,0)*+{\scriptscriptstyle A [(U \cdot U) \cdot L] A [R^{\cdot} \cdot U]} };
{\ar (70,0)*+{\scriptscriptstyle A [(U \cdot U) \cdot L] A [R^{\cdot} \cdot U]}; (70,-15)*+{\scriptscriptstyle A [U \cdot L] A [R^{\cdot} \cdot U]} };
{\ar^{\mu} (70,-15)*+{\scriptscriptstyle A [U \cdot L] A [R^{\cdot} \cdot U]}; (70,-30)*+{\scriptscriptstyle AU} };
{\ar (70,-30)*+{\scriptscriptstyle AU}; (70,-45)*+{\scriptscriptstyle A} };
{\ar (0,0)*+{\scriptscriptstyle [U \cdot (U \cdot L)][U \cdot A][U \cdot (R^{\cdot} \cdot U)]A}; (0,-15)*+{\scriptscriptstyle [(UUU) \cdot \big( (U \cdot L)A (R^{\cdot} \cdot U) \big)]A} };
{\ar (0,-15)*+{\scriptscriptstyle [(UUU) \cdot \big( (U \cdot L)A (R^{\cdot} \cdot U) \big)]A}; (0,-30)*+{\scriptscriptstyle [U \cdot \big( (U \cdot L)A (R^{\cdot} \cdot U) \big)]A} };
{\ar_{\mu}  (0,-30)*+{\scriptscriptstyle [U \cdot \big( (U \cdot L)A (R^{\cdot} \cdot U) \big)]A}; (0,-45)*+{\scriptscriptstyle [U \cdot U]A} };
{\ar (0,-45)*+{\scriptscriptstyle [U \cdot U]A}; (35,-45)*+{\scriptscriptstyle UA} };
{\ar (35,-45)*+{\scriptscriptstyle UA}; (70,-45)*+{\scriptscriptstyle A} };
\endxy
\]

\subsection{Functors between doubly degenerate tricategories}
Here we provide the diagrams omitted in the characterisation of weak functors between doubly degenerate tricategories.  First note that the dual pairs $\bs{\chi}$ and $\bs{\iota}$ must satisfy axioms similar to those for $\bs{A}, \bs{L},$ and $\bs{R}$.  Then, the two functor axioms become the two diagrams below.
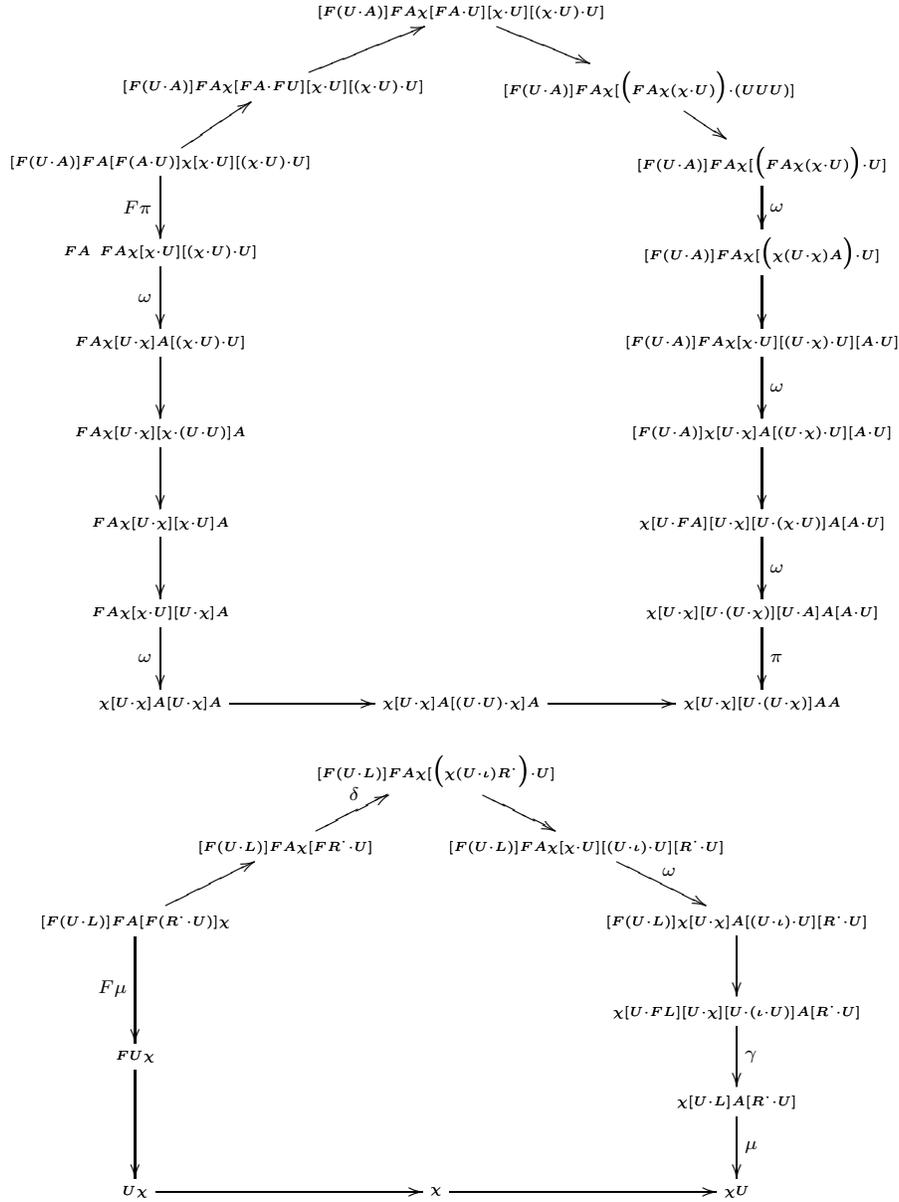
\[
\xy
{\ar (0,0)*+{\scriptscriptstyle [F(U\cdot A)]FA [F(A \cdot U)] \chi [\chi \cdot U][(\chi \cdot U) \cdot U]}; (15,10)*+{\scriptscriptstyle [F(U \cdot A)]FA \chi [FA \cdot FU][\chi \cdot U][(\chi \cdot U) \cdot U]} };
{\ar (15,10)*+{\scriptscriptstyle [F(U \cdot A)]FA \chi [FA \cdot FU][\chi \cdot U][(\chi \cdot U) \cdot U]}; (40,20)*+{\scriptscriptstyle [F(U \cdot A)]FA \chi [FA \cdot U][\chi \cdot U][(\chi \cdot U) \cdot U]} };
{\ar (40,20)*+{\scriptscriptstyle [F(U \cdot A)]FA \chi [FA \cdot U][\chi \cdot U][(\chi \cdot U) \cdot U]}; (65,10)*+{\scriptscriptstyle [F(U \cdot A)] FA \chi [ \big( FA \chi (\chi \cdot U) \big) \cdot (UUU)]} };
{\ar (65,10)*+{\scriptscriptstyle [F(U \cdot A)] FA \chi [ \big( FA \chi (\chi \cdot U) \big) \cdot (UUU)]}; (80,0)*+{\scriptscriptstyle [F(U \cdot A)] FA \chi [ \big( FA \chi (\chi \cdot U) \big) \cdot U]} };
{\ar^{\omega} (80,0)*+{\scriptscriptstyle [F(U \cdot A)] FA \chi [ \big( FA \chi (\chi \cdot U) \big) \cdot U]}; (80,-12)*+{\scriptscriptstyle [F(U \cdot A)]FA \chi [\big( \chi (U \cdot \chi) A \big) \cdot U]} };
{\ar (80,-12)*+{\scriptscriptstyle [F(U \cdot A)]FA \chi [\big( \chi (U \cdot \chi) A \big) \cdot U]}; (80,-24)*+{\scriptscriptstyle [F(U \cdot A)]FA \chi [\chi \cdot U][(U \cdot \chi) \cdot U][A \cdot U]} };
{\ar^{\omega} (80,-24)*+{\scriptscriptstyle [F(U \cdot A)]FA \chi [\chi \cdot U][(U \cdot \chi) \cdot U][A \cdot U]}; (80,-36)*+{\scriptscriptstyle [F(U \cdot A)]\chi [U \cdot \chi]A [(U \cdot \chi) \cdot U][A \cdot U]} };
{\ar (80,-36)*+{\scriptscriptstyle [F(U \cdot A)]\chi [U \cdot \chi]A [(U \cdot \chi) \cdot U][A \cdot U]}; (80,-48)*+{\scriptscriptstyle \chi [U \cdot FA][U \cdot \chi][U \cdot (\chi \cdot U)]A[A \cdot U]} };
{\ar^{\omega} (80,-48)*+{\scriptscriptstyle \chi [U \cdot FA][U \cdot \chi][U \cdot (\chi \cdot U)]A[A \cdot U]}; (80,-60)*+{\scriptscriptstyle \chi [U \cdot \chi][U \cdot (U \cdot \chi)][U \cdot A]A [A \cdot U]} };
{\ar^{\pi}  (80,-60)*+{\scriptscriptstyle \chi [U \cdot \chi][U \cdot (U \cdot \chi)][U \cdot A]A [A \cdot U]}; (80,-72)*+{\scriptscriptstyle \chi [U \cdot \chi][U \cdot (U \cdot \chi)]AA} };
{\ar_{F\pi} (0,0)*+{\scriptscriptstyle [F(U\cdot A)]FA [F(A \cdot U)] \chi [\chi \cdot U][(\chi \cdot U) \cdot U]}; (0,-12)*+{\scriptscriptstyle FA \phantom{1} FA \chi [\chi \cdot U][(\chi \cdot U) \cdot U]} };
{\ar_{\omega} (0,-12)*+{\scriptscriptstyle FA \phantom{1} FA \chi [\chi \cdot U][(\chi \cdot U) \cdot U]}; (0,-24)*+{\scriptscriptstyle FA \chi [U \cdot \chi]A [(\chi \cdot U) \cdot U]} };
{\ar (0,-24)*+{\scriptscriptstyle FA \chi [U \cdot \chi]A [(\chi \cdot U) \cdot U]}; (0,-36)*+{\scriptscriptstyle FA \chi [U \cdot \chi][\chi \cdot (U \cdot U)]A} };
{\ar (0,-36)*+{\scriptscriptstyle FA \chi [U \cdot \chi][\chi \cdot (U \cdot U)]A}; (0,-48)*+{\scriptscriptstyle FA \chi [U \cdot \chi][\chi \cdot U]A} };
{\ar (0,-48)*+{\scriptscriptstyle FA \chi [U \cdot \chi][\chi \cdot U]A}; (0,-60)*+{\scriptscriptstyle FA \chi [\chi \cdot U][U \cdot \chi]A} };
{\ar_{\omega} (0,-60)*+{\scriptscriptstyle FA \chi [\chi \cdot U][U \cdot \chi]A}; (0,-72)*+{\scriptscriptstyle \chi [U \cdot \chi]A[U \cdot \chi]A} };
{\ar (0,-72)*+{\scriptscriptstyle \chi [U \cdot \chi]A[U \cdot \chi]A}; (40,-72)*+{\scriptscriptstyle \chi [U \cdot \chi]A[(U\cdot U) \cdot \chi]A} };
{\ar (40,-72)*+{\scriptscriptstyle \chi [U \cdot \chi]A[(U\cdot U) \cdot \chi]A}; (80,-72)*+{\scriptscriptstyle \chi [U \cdot \chi][U \cdot (U \cdot \chi)]AA} }
\endxy
\]
\[
\xy
{\ar (0,0)*+{\scriptscriptstyle [F(U \cdot L)] FA [F(R^{\cdot} \cdot U)] \chi}; (20,10)*+{\scriptscriptstyle [F(U \cdot L)] FA \chi [FR^{\cdot} \cdot U]} };
{\ar^{\delta} (20,10)*+{\scriptscriptstyle [F(U \cdot L)] FA \chi [FR^{\cdot} \cdot U]}; (40,20)*+{\scriptscriptstyle [F(U \cdot L)] FA \chi [\big( \chi (U \cdot \iota) R^{\cdot} \big) \cdot U]} };
{\ar (40,20)*+{\scriptscriptstyle [F(U \cdot L)] FA \chi [\big( \chi (U \cdot \iota) R^{\cdot} \big) \cdot U]}; (60,10)*+{\scriptscriptstyle [F(U \cdot L)] FA \chi [\chi \cdot U][(U \cdot \iota) \cdot U][R^{\cdot} \cdot U]} };
{\ar^{\omega} (60,10)*+{\scriptscriptstyle [F(U \cdot L)] FA \chi [\chi \cdot U][(U \cdot \iota) \cdot U][R^{\cdot} \cdot U]}; (80,0)*+{\scriptscriptstyle [F(U \cdot L)] \chi [U \cdot \chi] A [(U \cdot \iota) \cdot U] [R^{\cdot} \cdot U]} };
{\ar (80,0)*+{\scriptscriptstyle [F(U \cdot L)] \chi [U \cdot \chi] A [(U \cdot \iota) \cdot U] [R^{\cdot} \cdot U]}; (80,-12)*+{\scriptscriptstyle \chi [U \cdot FL][U \cdot \chi][U \cdot (\iota \cdot U)]A [R^{\cdot} \cdot U]} };
{\ar^{\gamma}  (80,-12)*+{\scriptscriptstyle \chi [U \cdot FL][U \cdot \chi][U \cdot (\iota \cdot U)]A [R^{\cdot} \cdot U]}; (80,-24)*+{\scriptscriptstyle \chi [U \cdot L]A [R^{\cdot} \cdot U]} };
{\ar^{\mu} (80,-24)*+{\scriptscriptstyle \chi [U \cdot L]A [R^{\cdot} \cdot U]}; (80,-36)*+{\scriptscriptstyle \chi U} };
{\ar_{F\mu} (0,0)*+{\scriptscriptstyle [F(U \cdot L)] FA [F(R^{\cdot} \cdot U)] \chi}; (0,-18)*+{\scriptscriptstyle FU \chi} };
{\ar (0,-18)*+{\scriptscriptstyle FU \chi}; (0,-36)*+{\scriptscriptstyle  U \chi} };
{\ar (0,-36)*+{\scriptscriptstyle  U \chi}; (40,-36)*+{\scriptscriptstyle \chi} };
{\ar (40,-36)*+{\scriptscriptstyle \chi}; (80,-36)*+{\scriptscriptstyle \chi U} }
\endxy
\]

\subsection{Transformations for doubly degenerate tricategories}
Now we provide the diagrams omitted from the characterisation of weak transformations in the context of doubly degenerate tricategories.  As before, the dual pair $\bs{\alpha}$ must satisfy two axioms similar to those for $\bs{A}, \bs{L}$, and $\bs{R}$.  The three transformations axioms become the diagrams below.  Since we need lax transformations as well as weak ones, it should be noted that the diagrams below do not change in the lax case, except that we have a distinguished object $\alpha$ instead of the dual pair, and thus $\alpha$ satisfies axioms similar to those for $A, L,$ and $R$.
\[
\xy
{\ar (0,0)*+{\scriptscriptstyle [GA \cdot U][\chi_{G} \cdot U]A^{\cdot}[\chi_{G} \cdot U]A^{\cdot}[U \cdot (U \cdot \alpha)][U \cdot A]A[(U \cdot \alpha) \cdot U][A \cdot U][(\alpha \cdot U) \cdot U]}; (60,-10)*+{\scriptscriptstyle [GA \cdot U][\chi_{G} \cdot U]A^{\cdot} [\chi_{G} \cdot U][U \cdot \alpha]A^{\cdot}[U \cdot A]A[(U \cdot \alpha) \cdot U][A \cdot U][(\alpha \cdot U) \cdot U]} };
{\ar (60,-10)*+{\scriptscriptstyle [GA \cdot U][\chi_{G} \cdot U]A^{\cdot} [\chi_{G} \cdot U][U \cdot \alpha]A^{\cdot}[U \cdot A]A[(U \cdot \alpha) \cdot U][A \cdot U][(\alpha \cdot U) \cdot U]}; (60,-25)*+{\scriptscriptstyle [GA \cdot U][\chi_{G} \cdot U]A^{\cdot}[U \cdot \alpha] [\chi_{G} \cdot U]A^{\cdot}[U \cdot A]A[(U \cdot \alpha) \cdot U][A \cdot U][(\alpha \cdot U) \cdot U]} };
{\ar^{\pi} (60,-25)*+{\scriptscriptstyle [GA \cdot U][\chi_{G} \cdot U]A^{\cdot}[U \cdot \alpha] [\chi_{G} \cdot U]A^{\cdot}[U \cdot A]A[(U \cdot \alpha) \cdot U][A \cdot U][(\alpha \cdot U) \cdot U]}; (60,-37)*+{\scriptscriptstyle [GA \cdot U][\chi_{G} \cdot U]A^{\cdot}[U \cdot \alpha] [\chi_{G} \cdot U]A[A^{\cdot} \cdot U][(U \cdot \alpha) \cdot U][A \cdot U][(\alpha \cdot U) \cdot U]} };
{\ar (60,-37)*+{\scriptscriptstyle [GA \cdot U][\chi_{G} \cdot U]A^{\cdot}[U \cdot \alpha] [\chi_{G} \cdot U]A[A^{\cdot} \cdot U][(U \cdot \alpha) \cdot U][A \cdot U][(\alpha \cdot U) \cdot U]}; (60,-55)*+{\scriptscriptstyle [GA \cdot U][\chi_{G} \cdot U]A^{\cdot}[U \cdot \alpha]A [(\chi_{G} \cdot U) \cdot U][A^{\cdot} \cdot U][(U \cdot \alpha) \cdot U][A \cdot U][(\alpha \cdot U) \cdot U]} };
{\ar^{\Pi} (60,-55)*+{\scriptscriptstyle [GA \cdot U][\chi_{G} \cdot U]A^{\cdot}[U \cdot \alpha]A [(\chi_{G} \cdot U) \cdot U][A^{\cdot} \cdot U][(U \cdot \alpha) \cdot U][A \cdot U][(\alpha \cdot U) \cdot U]}; (60,-70)*+{\scriptscriptstyle [GA \cdot U][\chi_{G} \cdot U]A^{\cdot}[U \cdot \alpha]A [\alpha \cdot U][(U \cdot \chi_{F}) \cdot U][A \cdot U]} };
{\ar^{\Pi} (60,-70)*+{\scriptscriptstyle [GA \cdot U][\chi_{G} \cdot U]A^{\cdot}[U \cdot \alpha]A [\alpha \cdot U][(U \cdot \chi_{F}) \cdot U][A \cdot U]}; (60,-85)*+{\scriptscriptstyle [GA \cdot U]\alpha[U \cdot \chi_{F}]A [(U \cdot \chi_{F}) \cdot U][A \cdot U]} };
{\ar (60,-85)*+{\scriptscriptstyle [GA \cdot U]\alpha[U \cdot \chi_{F}]A [(U \cdot \chi_{F}) \cdot U][A \cdot U]}; (60,-100)*+{\scriptscriptstyle \alpha [U \cdot FA][U \cdot \chi_{F}][U \cdot (\chi_{F} \cdot U)]A [A \cdot U]} };
{\ar^{\omega} (60,-100)*+{\scriptscriptstyle \alpha [U \cdot FA][U \cdot \chi_{F}][U \cdot (\chi_{F} \cdot U)]A [A \cdot U]}; (60, -115)*+{\scriptscriptstyle \alpha [U \cdot \chi_{F}][U \cdot (U \cdot \chi_{F})][U \cdot A]A [A \cdot U]} };
{\ar^{\pi} (60, -115)*+{\scriptscriptstyle \alpha [U \cdot \chi_{F}][U \cdot (U \cdot \chi_{F})][U \cdot A]A [A \cdot U]}; (60, -130)*+{\scriptscriptstyle \alpha [U \cdot \chi_{F}][U \cdot (U \cdot \chi_{F})]AA} };
{\ar (60, -130)*+{\scriptscriptstyle \alpha [U \cdot \chi_{F}][U \cdot (U \cdot \chi_{F})]AA}; (60, -145)*+{\scriptscriptstyle \alpha [U \cdot \chi_{F}]A[(U \cdot U) \cdot \chi_{F}]A} };
{\ar (0,0)*+{\scriptscriptstyle [GA \cdot U][\chi_{G} \cdot U]A^{\cdot}[\chi_{G} \cdot U]A^{\cdot}[U \cdot (U \cdot \alpha)][U \cdot A]A[(U \cdot \alpha) \cdot U][A \cdot U][(\alpha \cdot U) \cdot U]}; (0,-20)*+{\scriptscriptstyle [GA \cdot U][\chi_{G} \cdot U][(\chi_{G} \cdot U) \cdot U]A^{\cdot}A^{\cdot}[U \cdot (U \cdot \alpha)][U \cdot A][U \cdot (\alpha \cdot U)]A[A \cdot U][(\alpha \cdot U) \cdot U]} };
{\ar_{\omega} (0,-20)*+{\scriptscriptstyle [GA \cdot U][\chi_{G} \cdot U][(\chi_{G} \cdot U) \cdot U]A^{\cdot}A^{\cdot}[U \cdot (U \cdot \alpha)][U \cdot A][U \cdot (\alpha \cdot U)]A[A \cdot U][(\alpha \cdot U) \cdot U]}; (0,-41)*+{\scriptscriptstyle [\chi_{G} \cdot U][(U \cdot \chi_{G}) \cdot U][A \cdot U]A^{\cdot}A^{\cdot}[U \cdot (U \cdot \alpha)][U \cdot A][U \cdot (\alpha \cdot U)]A[A \cdot U][(\alpha \cdot U) \cdot U]} };
{\ar_{\pi} (0,-41)*+{\scriptscriptstyle [\chi_{G} \cdot U][(U \cdot \chi_{G}) \cdot U][A \cdot U]A^{\cdot}A^{\cdot}[U \cdot (U \cdot \alpha)][U \cdot A][U \cdot (\alpha \cdot U)]A[A \cdot U][(\alpha \cdot U) \cdot U]}; (0,-62)*+{\scriptscriptstyle [\chi_{G} \cdot U][(U \cdot \chi_{G}) \cdot U]A^{\cdot}[U \cdot A^{\cdot}][U \cdot (U \cdot \alpha)][U \cdot A][U \cdot (\alpha \cdot U)]A[A \cdot U][(\alpha \cdot U) \cdot U]} };
{\ar (0,-62)*+{\scriptscriptstyle [\chi_{G} \cdot U][(U \cdot \chi_{G}) \cdot U]A^{\cdot}[U \cdot A^{\cdot}][U \cdot (U \cdot \alpha)][U \cdot A][U \cdot (\alpha \cdot U)]A[A \cdot U][(\alpha \cdot U) \cdot U]}; (0,-83)*+{\scriptscriptstyle [\chi_{G} \cdot U]A^{\cdot}[U \cdot (\chi_{G} \cdot U)][U \cdot A^{\cdot}][U \cdot (U \cdot \alpha)][U \cdot A][U \cdot (\alpha \cdot U)]A[A \cdot U][(\alpha \cdot U) \cdot U]} };
{\ar_{\Pi} (0,-83)*+{\scriptscriptstyle [\chi_{G} \cdot U]A^{\cdot}[U \cdot (\chi_{G} \cdot U)][U \cdot A^{\cdot}][U \cdot (U \cdot \alpha)][U \cdot A][U \cdot (\alpha \cdot U)]A[A \cdot U][(\alpha \cdot U) \cdot U]}; (0,-104)*+{\scriptscriptstyle [\chi_{G} \cdot U]A^{\cdot}[U \cdot \alpha][U \cdot (U \cdot \chi_{F})][U \cdot A]A[A\cdot U][(A \cdot U) \cdot U]} };
{\ar_{\pi} (0,-104)*+{\scriptscriptstyle [\chi_{G} \cdot U]A^{\cdot}[U \cdot \alpha][U \cdot (U \cdot \chi_{F})][U \cdot A]A[A\cdot U][(A \cdot U) \cdot U]}; (0,-125)*+{\scriptscriptstyle [\chi_{G} \cdot U]A^{\cdot}[U \cdot \alpha][U \cdot (U \cdot \chi_{F})]AA[(A \cdot U) \cdot U]} };
{\ar (0,-125)*+{\scriptscriptstyle [\chi_{G} \cdot U]A^{\cdot}[U \cdot \alpha][U \cdot (U \cdot \chi_{F})]AA[(A \cdot U) \cdot U]}; (0,-135)*+{\scriptscriptstyle [\chi_{G} \cdot U]A^{\cdot}[U \cdot \alpha]A[U \cdot \chi_{F}][\alpha \cdot U]A} };
{\ar (0,-135)*+{\scriptscriptstyle [\chi_{G} \cdot U]A^{\cdot}[U \cdot \alpha]A[U \cdot \chi_{F}][\alpha \cdot U]A}; (20,-145)*+{\scriptscriptstyle [\chi_{G} \cdot U]A^{\cdot}[U \cdot \alpha]A[\alpha \cdot U][U \cdot \chi_{F}]A} };
{\ar (20,-145)*+{\scriptscriptstyle [\chi_{G} \cdot U]A^{\cdot}[U \cdot \alpha]A[\alpha \cdot U][U \cdot \chi_{F}]A}; (60, -145)*+{\scriptscriptstyle \alpha [U \cdot \chi_{F}]A[(U \cdot U) \cdot \chi_{F}]A} };
\endxy
\]

\[
\xy
{\ar^{M} (0,0)*+{\scriptscriptstyle [GL \cdot U][\chi_{G} \cdot U]A^{\cdot}[U \cdot \alpha]A[\alpha \cdot U][(U \cdot \iota_{F}) \cdot U][R^{\cdot} \cdot U]}; (60,0)*+{\scriptscriptstyle [GL \cdot U][\chi_{G} \cdot U]A^{\cdot}[U \cdot \alpha]A [(\iota_{G} \cdot U) \cdot U][L^{\cdot} \cdot U]} };
{\ar (60,0)*+{\scriptscriptstyle [GL \cdot U][\chi_{G} \cdot U]A^{\cdot}[U \cdot \alpha]A [(\iota_{G} \cdot U) \cdot U][L^{\cdot} \cdot U]}; (60,-12)*+{\scriptscriptstyle [GL \cdot U][\chi_{G} \cdot U]A^{\cdot}[U \cdot \alpha][\iota_{G} \cdot U]A [L^{\cdot} \cdot U]} };
{\ar (60,-12)*+{\scriptscriptstyle [GL \cdot U][\chi_{G} \cdot U]A^{\cdot}[U \cdot \alpha][\iota_{G} \cdot U]A [L^{\cdot} \cdot U]}; (60,-24)*+{\scriptscriptstyle [GL \cdot U][\chi_{G} \cdot U]A^{\cdot}[\iota_{G} \cdot U][U \cdot \alpha]A [L^{\cdot} \cdot U]} };
{\ar (60,-24)*+{\scriptscriptstyle [GL \cdot U][\chi_{G} \cdot U]A^{\cdot}[\iota_{G} \cdot U][U \cdot \alpha]A [L^{\cdot} \cdot U]}; (60,-36)*+{\scriptscriptstyle [GL \cdot U][\chi_{G} \cdot U][(\iota_{G} \cdot U) \cdot U]A^{\cdot}[U \cdot \alpha]A [L^{\cdot} \cdot U]} };
{\ar^{\gamma_{G}} (60,-36)*+{\scriptscriptstyle [GL \cdot U][\chi_{G} \cdot U][(\iota_{G} \cdot U) \cdot U]A^{\cdot}[U \cdot \alpha]A [L^{\cdot} \cdot U]}; (60,-48)*+{\scriptscriptstyle [L \cdot U]A^{\cdot} [U \cdot \alpha]A[L^{\cdot} \cdot U]} };
{\ar^{\lambda} (60,-48)*+{\scriptscriptstyle [L \cdot U]A^{\cdot} [U \cdot \alpha]A[L^{\cdot} \cdot U]}; (60,-60)*+{\scriptscriptstyle L [U \cdot \alpha]A [L^{\cdot} \cdot U]} };
{\ar_{\Pi} (0,0)*+{\scriptscriptstyle [GL \cdot U][\chi_{G} \cdot U]A^{\cdot}[U \cdot \alpha]A[\alpha \cdot U][(U \cdot \iota_{F}) \cdot U][R^{\cdot} \cdot U]}; (0,-15)*+{\scriptscriptstyle [GL \cdot U]\alpha[U \cdot \chi_{F}]A[(U \cdot \iota_{F}) \cdot U][R^{\cdot} \cdot U]} };
{\ar (0,-15)*+{\scriptscriptstyle [GL \cdot U]\alpha[U \cdot \chi_{F}]A[(U \cdot \iota_{F}) \cdot U][R^{\cdot} \cdot U]}; (0,-30)*+{\scriptscriptstyle [GL \cdot U]\alpha[U \cdot \chi_{F}][U \cdot (\iota_{F} \cdot U)]A[R^{\cdot} \cdot U]} };
{\ar (0,-30)*+{\scriptscriptstyle [GL \cdot U]\alpha[U \cdot \chi_{F}][U \cdot (\iota_{F} \cdot U)]A[R^{\cdot} \cdot U]}; (0,-45)*+{\scriptscriptstyle \alpha [U \cdot FL] [U \cdot \chi_{F}][U \cdot (\iota_{F} \cdot U)]A[R^{\cdot} \cdot U]} };
{\ar_{\gamma_{F}} (0,-45)*+{\scriptscriptstyle \alpha [U \cdot FL] [U \cdot \chi_{F}][U \cdot (\iota_{F} \cdot U)]A[R^{\cdot} \cdot U]}; (0,-60)*+{\scriptscriptstyle \alpha[U \cdot L]A[R^{\cdot} \cdot U] } };
{\ar_{\mu} (0,-60)*+{\scriptscriptstyle \alpha[U \cdot L]A[R^{\cdot} \cdot U] }; (20,-60)*+{\scriptscriptstyle \alpha U} };
{\ar_{\lambda} (20,-60)*+{\scriptscriptstyle \alpha U}; (40,-60)*+{\scriptscriptstyle \alpha LA [L^{\cdot} \cdot U]} };
{\ar (40,-60)*+{\scriptscriptstyle \alpha LA [L^{\cdot} \cdot U]}; (60,-60)*+{\scriptscriptstyle L [U \cdot \alpha]A [L^{\cdot} \cdot U]} }
\endxy
\]

We will now drop the subscripts for $\chi$ and $\iota$ as they are determined by the rest of the diagram.

\[
\xy
{\ar^{\rho} (0,0)*+{\scriptscriptstyle [\chi \cdot U]A^{\cdot}[U \cdot \alpha][U \cdot (U \cdot \iota)][U \cdot R^{\cdot}]\alpha}; (20,10)*+{\scriptscriptstyle [\chi \cdot U]A^{\cdot}[U \cdot \alpha][U \cdot (U \cdot \iota)]AR^{\cdot} \alpha} };
{\ar (20,10)*+{\scriptscriptstyle [\chi \cdot U]A^{\cdot}[U \cdot \alpha][U \cdot (U \cdot \iota)]AR^{\cdot} \alpha}; (40,20)*+{\scriptscriptstyle [\chi \cdot U]A^{\cdot}[U \cdot \alpha]A[\alpha \cdot U][U \cdot \iota]R^{\cdot}} };
{\ar^{\Pi} (40,20)*+{\scriptscriptstyle [\chi \cdot U]A^{\cdot}[U \cdot \alpha]A[\alpha \cdot U][U \cdot \iota]R^{\cdot}}; (60,10)*+{\scriptscriptstyle \alpha [U \cdot \chi]A [U \cdot \iota]R^{\cdot}} };
{\ar (60,10)*+{\scriptscriptstyle \alpha [U \cdot \chi]A [U \cdot \iota]R^{\cdot}}; (80,0)*+{\scriptscriptstyle \alpha [U \cdot \chi][U \cdot (U \cdot \iota)]AR^{\cdot} } };
{\ar^{\delta} (80,0)*+{\scriptscriptstyle \alpha [U \cdot \chi][U \cdot (U \cdot \iota)]AR^{\cdot} }; (80,-15)*+{\scriptscriptstyle \alpha [U \cdot FR^{\cdot}][U \cdot R]AR^{\cdot}} };
{\ar^{\rho} (80,-15)*+{\scriptscriptstyle \alpha [U \cdot FR^{\cdot}][U \cdot R]AR^{\cdot}}; (80,-30)*+{\scriptscriptstyle \alpha [U \cdot FR^{\cdot}]U} };
{\ar (80,-30)*+{\scriptscriptstyle \alpha [U \cdot FR^{\cdot}]U}; (80,-45)*+{\scriptscriptstyle [GR^{\cdot} \cdot U]\alpha} };
{\ar_{M} (0,0)*+{\scriptscriptstyle [\chi \cdot U]A^{\cdot}[U \cdot \alpha][U \cdot (U \cdot \iota)][U \cdot R^{\cdot}]\alpha}; (0,-15)*+{\scriptscriptstyle [\chi \cdot U]A^{\cdot}[U \cdot (\iota \cdot U)][U \cdot L^{\cdot}]\alpha} };
{\ar (0,-15)*+{\scriptscriptstyle [\chi \cdot U]A^{\cdot}[U \cdot (\iota \cdot U)][U \cdot L^{\cdot}]\alpha}; (0,-30)*+{\scriptscriptstyle [\chi \cdot U][(U \cdot \iota) \cdot U]A^{\cdot}[U \cdot L^{\cdot}]\alpha} };
{\ar_{\delta} (0,-30)*+{\scriptscriptstyle [\chi \cdot U][(U \cdot \iota) \cdot U]A^{\cdot}[U \cdot L^{\cdot}]\alpha}; (0,-45)*+{\scriptscriptstyle [GR^{\cdot} \cdot U][R \cdot U]A^{\cdot}[U \cdot L^{\cdot}]\alpha} };
{\ar_{\mu} (0,-45)*+{\scriptscriptstyle [GR^{\cdot} \cdot U][R \cdot U]A^{\cdot}[U \cdot L^{\cdot}]\alpha}; (40,-45)*+{\scriptscriptstyle [GR^{\cdot} \cdot U]U \alpha} };
{\ar (40,-45)*+{\scriptscriptstyle [GR^{\cdot} \cdot U]U \alpha}; (80,-45)*+{\scriptscriptstyle [GR^{\cdot} \cdot U]\alpha} }
\endxy
\]

\subsection{Modifications for doubly degenerate tricategories}
Here we give the two axioms for a modification $m: \alpha \Rightarrow \beta$ in the context of doubly degenerate tricategories.  Once again we will only mark some of the maps, and the maps marked by $m$ are obtained from the single isomorphism in the definition of the modification.
\[
\xy
{\ar (0,0)*+{\scriptstyle [GU \cdot m][\chi \cdot U]A^{\cdot}[U \cdot \alpha]A[\alpha \cdot U]}; (60,0)*+{\scriptstyle [\chi \cdot U][(GU \cdot GU) \cdot m]A^{\cdot}[U \cdot \alpha]A[\alpha \cdot U]} };
{\ar (60,0)*+{\scriptstyle [\chi \cdot U][(GU \cdot GU) \cdot m]A^{\cdot}[U \cdot \alpha]A[\alpha \cdot U]}; (60,-15)*+{\scriptstyle [\chi \cdot U]A^{\cdot}[GU \cdot (GU \cdot m)][U \cdot \alpha]A[\alpha \cdot U]} };
{\ar^{m} (60,-15)*+{\scriptstyle [\chi \cdot U]A^{\cdot}[GU \cdot (GU \cdot m)][U \cdot \alpha]A[\alpha \cdot U]}; (60,-30)*+{\scriptstyle [\chi \cdot U]A^{\cdot}[U \cdot \beta][GU \cdot (m \cdot FU)]A [\alpha \cdot U]} };
{\ar (60,-30)*+{\scriptstyle [\chi \cdot U]A^{\cdot}[U \cdot \beta][GU \cdot (m \cdot FU)]A [\alpha \cdot U]}; (60,-45)*+{\scriptstyle [\chi \cdot U]A^{\cdot}[U \cdot \beta]A[(GU \cdot m) \cdot FU] [\alpha \cdot U]} };
{\ar^{m} (60,-45)*+{\scriptstyle [\chi \cdot U]A^{\cdot}[U \cdot \beta]A[(GU \cdot m) \cdot FU] [\alpha \cdot U]}; (60,-60)*+{\scriptstyle [\chi \cdot U]A^{\cdot}[U \cdot \beta]A[\beta \cdot U][(m \cdot FU) \cdot FU]} };
{\ar^{\Pi} (60,-60)*+{\scriptstyle [\chi \cdot U]A^{\cdot}[U \cdot \beta]A[\beta \cdot U][(m \cdot FU) \cdot FU]}; (60,-75)*+{\scriptstyle \beta [U \cdot \chi]A[(m \cdot FU) \cdot FU]} };
{\ar_{\Pi} (0,0)*+{\scriptstyle [GU \cdot m][\chi \cdot U]A^{\cdot}[U \cdot \alpha]A[\alpha \cdot U]}; (0,-25)*+{\scriptstyle [GU \cdot m]\alpha [U \cdot \chi]A} };
{\ar_{m} (0,-25)*+{\scriptstyle [GU \cdot m]\alpha [U \cdot \chi]A}; (0,-50)*+{\scriptstyle \beta [m \cdot F(U \cdot U)][U \cdot \chi]A} };
{\ar (0,-50)*+{\scriptstyle \beta [m \cdot F(U \cdot U)][U \cdot \chi]A}; (0,-75)*+{\scriptstyle \beta [U \cdot \chi][m \cdot (FU \cdot FU)]A} };
{\ar (0,-75)*+{\scriptstyle \beta [U \cdot \chi][m \cdot (FU \cdot FU)]A}; (60,-75)*+{\scriptstyle \beta [U \cdot \chi]A[(m \cdot FU) \cdot FU]} }
\endxy
\]

\[
\xy
{\ar^{m} (0,0)*+{\scriptstyle [U \cdot m]\alpha[U \cdot \iota]R^{\cdot}}; (40,0)*+{\scriptstyle \beta [m \cdot U][U \cdot \iota]R^{\cdot}} };
{\ar (40,0)*+{\scriptstyle \beta [m \cdot U][U \cdot \iota]R^{\cdot}}; (40,-15)*+{\scriptstyle \beta [U \cdot \iota][m \cdot U]R^{\cdot}} };
{\ar (40,-15)*+{\scriptstyle \beta [U \cdot \iota][m \cdot U]R^{\cdot}}; (40,-30)*+{\scriptstyle \beta [U \cdot \iota]R^{\cdot} m} };
{\ar^{M} (40,-30)*+{\scriptstyle \beta [U \cdot \iota]R^{\cdot} m}; (40,-45)*+{\scriptstyle  [\iota \cdot U]L^{\cdot}m} };
{\ar_{M} (0,0)*+{\scriptstyle [U \cdot m]\alpha[U \cdot \iota]R^{\cdot}}; (0,-22.5)*+{\scriptstyle [U \cdot m][\iota \cdot U] L^{\cdot}} };
{\ar (0,-22.5)*+{\scriptstyle [U \cdot m][\iota \cdot U] L^{\cdot}}; (0,-45)*+{\scriptstyle [\iota \cdot U][U \cdot m]L^{\cdot} } };
{\ar (0,-45)*+{\scriptstyle [\iota \cdot U][U \cdot m]L^{\cdot} }; (40,-45)*+{\scriptstyle  [\iota \cdot U]L^{\cdot}m} }
\endxy
\]

\subsection{Perturbations for doubly degenerate tricategories}
Here we give the single perturbation axiom in the context of doubly degenerate tricategories.  The axiom for the perturbation $\sigma: m \Rrightarrow n$ becomes the diagram below.

\[
\xy
{\ar^{m} (0,0)*+{[U \cdot m]\alpha}; (40,0)*+{\beta [m \cdot U]} };
{\ar^{1[\sigma \cdot 1]} (40,0)*+{\beta [m \cdot U]}; (40,-18)*+{\beta[n \cdot U]} };
{\ar_{[1 \cdot \sigma]1} (0,0)*+{[U \cdot m]\alpha}; (0,-18)*+{[U \cdot n]\alpha} };
{\ar_{n} (0,-18)*+{[U \cdot n]\alpha}; (40,-18)*+{\beta[n \cdot U]} };
\endxy
\]

\bibliography{../bib/bib0706}

\end{document}